\documentclass[12pt]{article}
\usepackage[a4paper, total={7in, 8.7in}]{geometry}
\usepackage{verbatim}
\usepackage{dsfont}			
\usepackage{multirow}
\usepackage[mathscr]{euscript}
\usepackage{fancyhdr}
\usepackage{xcolor}
\usepackage[section]{placeins}

\usepackage{csquotes}
\usepackage[style=ext-authoryear, maxcitenames=1, giveninits=true, maxbibnames=99, sortcites=false, natbib, uniquename=init]{biblatex} 
\DeclareNameAlias{author}{family-given}
\DeclareNameAlias{editor}{family-given}

\renewbibmacro{in:}{
  \ifentrytype{article}
    {}
    {\bibstring{in}
     \printunit{\intitlepunct}}}
\DeclareFieldFormat[article]{pages}{#1}

\usepackage{float}
\usepackage{wrapfig}
\newcommand{\blind}{1}

\usepackage[english]{babel}

\usepackage{amsthm, amsmath, amssymb, mathrsfs, mathtools, bm, dutchcal}
\usepackage{cancel}
\usepackage{array}
\usepackage{graphicx}
\usepackage{subcaption}
\usepackage{bigints}
\usepackage[colorlinks=true, allcolors=blue]{hyperref}
\usepackage{booktabs}
\usepackage[capitalise]{cleveref}
\usepackage{listings}
\usepackage{enumerate}

\newtheorem{theorem}{\underline{\bf Theorem}}
\newtheorem{corollary}{\underline{\bf Corollary}}

\newtheorem{proposition}{\underline{\bf Proposition}}

\theoremstyle{definition}
\newtheorem{definition}{Definition}
\newtheorem{remark}{\underline{\bf Remark}}

\newtheorem{example}{Example}

\makeatletter
\newtheorem*{rep@theorem}{\rep@title}
\newcommand{\newreptheorem}[2]{%
	\newenvironment{rep#1}[1]{%
		\def\rep@title{#2 \ref{##1} (Continued)}%
		\begin{rep@theorem}}%
		{\end{rep@theorem}}}
\makeatother

\newreptheorem{example}{Example}

\begin{document}
\newrefsection

\definecolor{amethyst}{rgb}{0.6, 0.4, 0.8}
\newcommand{\gio}{\color{red}}
\newcommand{\bea}{\color{amethyst}}
\newcommand{\ech}{\rm \color{black}}
\newcommand*{\swap}[2]{#2#1} 


\def\Cor{\mathbb{C}\mathrm{or}}
\def\Cov{\mathbb{C}\mathrm{ov}}
\def\Var{\mathbb{V}\mathrm{ar}}
\def\bbE{\mathbb{E}}
\def\bbP{\mathbb{P}}
\def\bbN{\mathbb{N}}
\def\bbX{\mathbb{X}}
\def\bbR{\mathbb{R}}
\def\eqas{\stackrel{\mbox{\scriptsize{a.s.}}}{=}}
\def\eqd{\stackrel{\mbox{\scriptsize{d}}}{=}}
\def\simind{\stackrel{\mbox{\scriptsize{ind}}}{\sim}}
\def\simiid{\stackrel{\mbox{\scriptsize{iid}}}{\sim}}

\def\d{ \mathrm{d} }

\def\A{ {\mathcal{A}} }
\def\B{ {\mathcal{B}} }
\def\D{ {\mathcal{D}} }
\def\L{ \mathcal{L} }
\def\H{ \mathcal{H} }
\def\I{ \mathcal{I} }
\def\J{ {\mathcal{J}} }
\def\P{ {\mathcal{P}} }

\def\tP{ {\tilde {P}} }
\def\tpi{ {\tilde {\pi}} }

\def\Bern{\mathrm{Bernoulli}}
\def\Beta{\mathrm{Beta}}
\def\DE{\mathrm{DE}}
\def\Dir{\mathrm{Dir}}
\def\DM{\mathrm{DM}}
\def\DP{\mathrm{DP}}
\def\Exp{\mathrm{Exp}}
\def\EPPF{\mathrm{EPPF}}
\def\Ga{\mathrm{Gamma}}
\def\GN{\mathrm{GN}}
\def\GP{\mathrm{GP}}
\def\GEM{\mathrm{GEM}}
\def\HDP{\mathrm{HDP}}
\def\HHDP{\mathrm{HHDP}}
\def\HHPYP{\mathrm{HHPYP}}
\def\HOHMM{\mathrm{HOHMM}}
\def\HPYP{\mathrm{HPYP}}
\def\HSSP{\mathrm{HSSP}}
\def\IGs{\mathrm{Inv-Gs}}
\def\IG{\mathrm{Inv-Ga}}
\def\IW{\mathrm{IW}}
\def\LN{\mathrm{Log-Normal}}
\def\Mult{\mathrm{Mult}}
\def\MVN{\mathrm{MVN}}
\def\mSSP{\mathrm{mSSP}}
\def\Normal{\mathrm{Normal}}
\def\NDP{\mathrm{NDP}}
\def\NPYP{\mathrm{NPYP}}
\def\pEPPF{\mathrm{pEPPF}}
\def\Poi{\mathrm{Poi}}
\def\PYP{\mathrm{PYP}}
\def\SB{\mathrm{SB}}
\def\SSP{\mathrm{SSP}}
\def\Unif{\mathrm{Uniform}}
\def\Wish{\mathrm{W}}

\def\bn{ \bm{n} }
\def\bq{ \bm{q} }
\def\bX{ \mathbf{X}}
\def\bpi{ \bm{\pi} }

\date{}

\title{\bf Multivariate Species Sampling Models}

\if0\blind
\author{}
\else
\author{
Beatrice Franzolini,\\
DEMS Department,
University of Milano-Bicocca
\vspace*{0.2cm}\\
Antonio Lijoi, Igor Pr\"unster\\
Bocconi Institute for Data Science and Analytics, Bocconi University
\vspace*{0.2cm}\\
Giovanni Rebaudo \\
ESOMAS Department, University of Torino and Collegio Carlo Alberto
}
\fi

\maketitle

\bigskip
\abstract{Species sampling processes have long served as the fundamental framework for modeling random discrete distributions and exchangeable sequences.
However, data arising from distinct but related sources require a broader notion of probabilistic invariance, making partial exchangeability a natural choice.
Countless models for partially exchangeable data, collectively known as dependent nonparametric priors, have been proposed.
These include hierarchical, nested and additive processes, widely used in statistics and machine learning.
Still, a unifying framework is lacking and key questions about their underlying learning mechanisms remain unanswered.\\
We fill this gap by introducing multivariate species sampling models, a new general class of nonparametric priors that encompasses most existing finite- and infinite-dimensional dependent processes.
They are characterized by the induced partially exchangeable partition probability function encoding their multivariate clustering structure.
We establish their core distributional properties and analyze their dependence structure, demonstrating that borrowing of information across groups is entirely determined by shared ties.
This provides new insights into the underlying learning mechanisms, offering, for instance, a principled rationale for the previously unexplained correlation structure observed in existing models.
Beyond providing a cohesive theoretical foundation, our approach serves as a constructive tool for developing new models and opens novel research directions for capturing richer dependence structures beyond the framework of multivariate species sampling processes.
}

\vfill 
\vfill

\section{Introduction}
\label{sec: intro}

A fundamental homogeneity assumption in Bayesian inference is the (infinite) exchangeability of observables, which corresponds to distributional invariance with respect to permutations of the data.
According to de Finetti's representation theorem, there is a one-to-one correspondence between the law of an exchangeable sequence $(X_{i})_{i\geq 1}$ and a random probability measure, conditionally on which the $X_{i}$'s are independent and identically distributed (i.i.d.).
This foundational result supports the Bayesian framework, based on likelihood and prior, via a probabilistic symmetry assumption on the data.
At the same time, it establishes a conceptual bridge to the classical i.i.d.\ assumption in frequentist inference.
In a parametric setup, the random probability measure associated with $(X_{i})_{i\geq 1}$ is indexed by a finite-dimensional parameter, whereas in a nonparametric setting, no such restriction is imposed.
A cornerstone of the latter is represented by the Dirichlet process (DP) introduced in \textcite{ferguson1973bayesian}.
Its full weak support property implies remarkable flexibility compared to parametric counterparts, and several popular Bayesian nonparametric (BNP) models can be seen as extensions of the DP itself.

In a seminal work, \textcite{pitman1996some} introduced a unifying framework for studying almost surely discrete random probability measures under the assumption that weights and locations are independent.
The resulting class, known as species sampling processes (SSPs) \citep{ghosal2017fundamentals}, includes the DP as a special case and satisfies several structural properties that have been pivotal to understanding and constructing discrete priors for modeling exchangeable data.
Entire classes of popular nonparametric and parametric priors, such as homogeneous normalized random measures with independent increments \citep{regazzini2003distributional, james2009posterior}, Gibbs-type priors \citep{gnedin2006exchangeable, deblasi2015gibbs}, and stick-breaking processes \citep{ishwaran2001gibbs, gil2023stick} fall within the framework of SSPs, and hence also their notable special cases, which include Pitman--Yor \citep{pitman1997two}, normalized inverse Gaussian \citep{lijoi2005hierarchical} and normalized generalized gamma \citep{lijoi2007controlling} processes.
Moreover, relevant special cases also include finite-dimensional processes such as the finite Dirichlet multinomial \citep{green2001modelling} and mixture of finite-dimensional processes with a prior on the number of locations \citep{nobile1994bayesian, richardson1997bayesian, nobile2007bayesian, gnedin2010species, deblasi2013asymptotic, miller2018mixture}.
A summary of J.~Pitman's theory on univariate SSP \citep{pitman1996some} can be found in Section~\ref{sec: SSP supp} of the Supplementary Material.

However, exchangeability is often too restrictive an assumption in applied settings.
The pioneering contributions of \textcite{maceachern1999dependent, maceachern2000dependent} opened a new research line in both the statistics and machine learning literature, whose goal is to develop models that accommodate heterogeneity across data sources or experimental conditions.
More specifically, when data originate from $J$ distinct populations, such as in meta-analysis, topic modeling or multi-center studies, the exchangeability assumption becomes overly restrictive, since it fails to account for heterogeneity across distinct groups.
Conversely, assuming independence between populations precludes information sharing across experiments, which is often a key goal in multi-sample studies \citep[see, for instance,][]{woodcock2017master, chen2019bayesian, ouma2022bayesian, su2022comparative}.
A natural compromise between these extremes is the probabilistic framework of partial exchangeability \citep{finetti1938condition}, which entails exchangeability within but not across different populations, while still allowing for dependence among them.
Consider a random array $\bX$ with $J$ rows and infinite columns.
Then $\bX$ is partially exchangeable if and only if its distribution is invariant with respect to finite permutations within each row but not necessarily across rows.
This means that elements within each population, i.e., belonging to the same row, are exchangeable, but permuting elements across populations, i.e., belonging to different rows, would alter the distribution of $\bX$.
For instance, suppose we have partially exchangeable binary data from $J=2$ groups and that we observe $X_{1,1}$ and $X_{1,2}$ from the first group, and $X_{2,1}$ from the second.
Then, $\bbP(X_{1,1}=1,\,X_{1,2}=0,\, X_{2,1}= 1)=\bbP(X_{1,2}=1,\, X_{1,1}=0, \, X_{2,1}=1)$ as this is a permutation within group $1$.
However, it is possible that $\bbP(X_{1,1}=1,\, X_{1,2}=0, \, X_{2,1}=1) \ne \bbP(X_{1,1}= 1, \, X_{2,1}=0, \, X_{1,2}=1)$, since invariance with respect to permutations across groups is not preserved.
Similarly to the exchangeable case, partial exchangeability implies the existence of a vector of (dependent) random probability measures $(P_{1}, \ldots, P_{J})$, such that $X_{j, i}\mid P_{1}, \ldots, P_{J} \simind P_{j}$, for $i\geq1$ and $j=1,\ldots, J$.
From a Bayesian perspective, modeling a partially exchangeable array is equivalent to defining a prior distribution for a vector of dependent probability measures.
Countless approaches have been proposed in the literature and several success stories have been recorded.
Notable instances are hierarchical DPs \citep{teh2006hierarchical}, hierarchical normalized completely random measures \citep{camerlenghi2019distribution}, hierarchical species sampling models \citep{bassetti2020hierarchical}, nested constructions \citep{rodriguez2008nested, camerlenghi2019latent}, additive constructions \citep{muller2004method,lijoi2014bayesian}, copula constructions \citep{leisen2011vectors}, normalized compound random measures \citep{griffin2017compound}, normalized completely random vectors \citep{lijoi2014bayesian,catalano2021measuring}, single-atoms dependent processes \citep{maceachern1999dependent, maceachern2000dependent, quintana2022dependent}, compositions of some of the previous \citep{camerlenghi2019latent, beraha2021semi, lijoi2023flexible, balocchi2023clustering, denti2023common}, and many others \citep[e.g.][]{horiguchi2025tree, yan2023bayesian, bi2023class, lee2025logistic}.

\begin{figure}[tb]
	\begin{subfigure}{0.32\linewidth}
		\centering
		\includegraphics[ width=\linewidth]{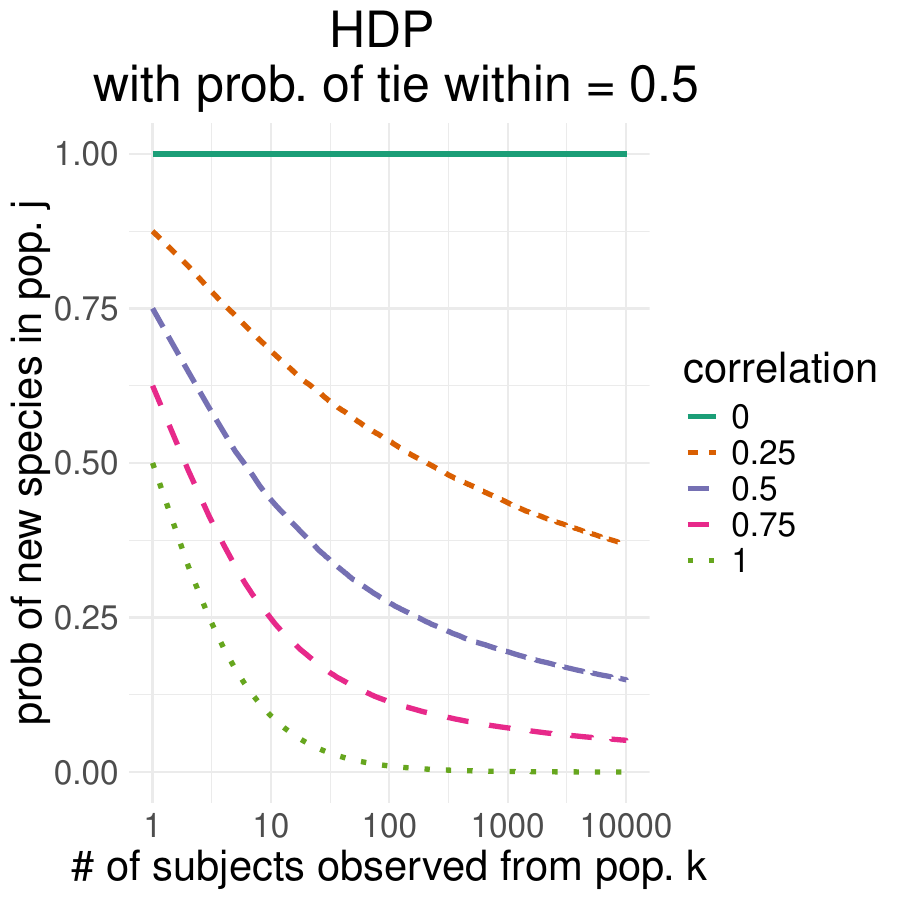}
	\end{subfigure}
 \hfill
 \begin{subfigure}{0.32\linewidth}
		\centering
		\includegraphics[ width=\linewidth]{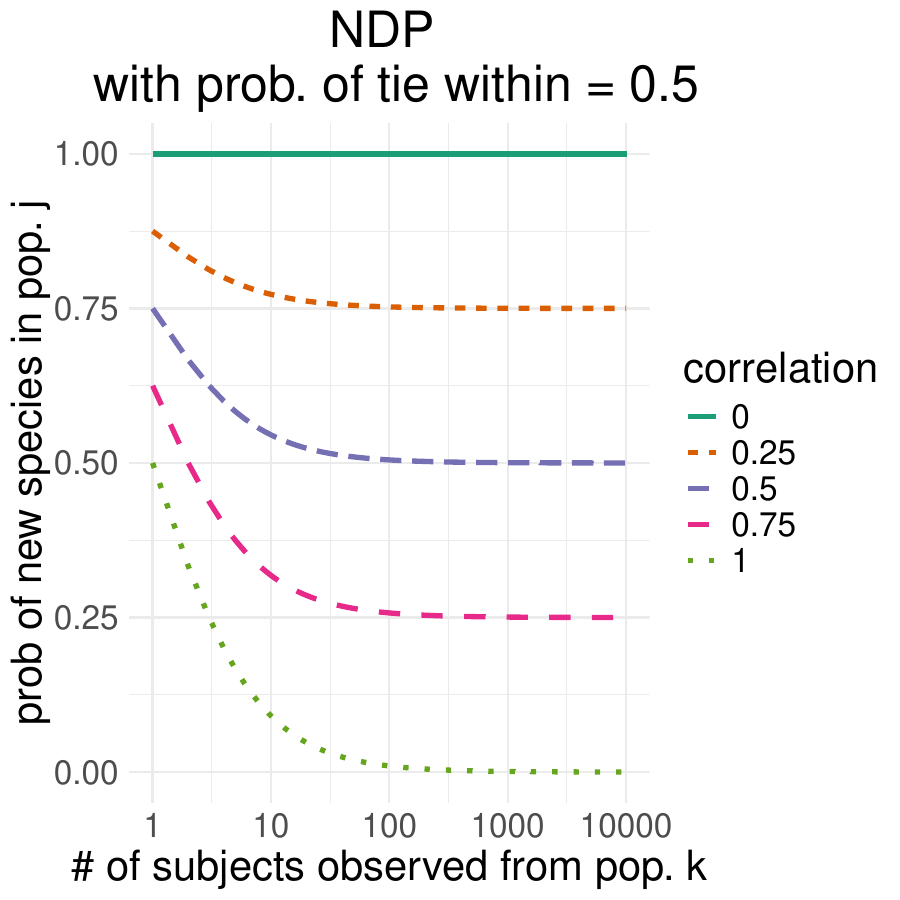}
	\end{subfigure}
 \hfill
		\begin{subfigure}{0.32\linewidth}
		\centering
		\includegraphics[ width=\linewidth]{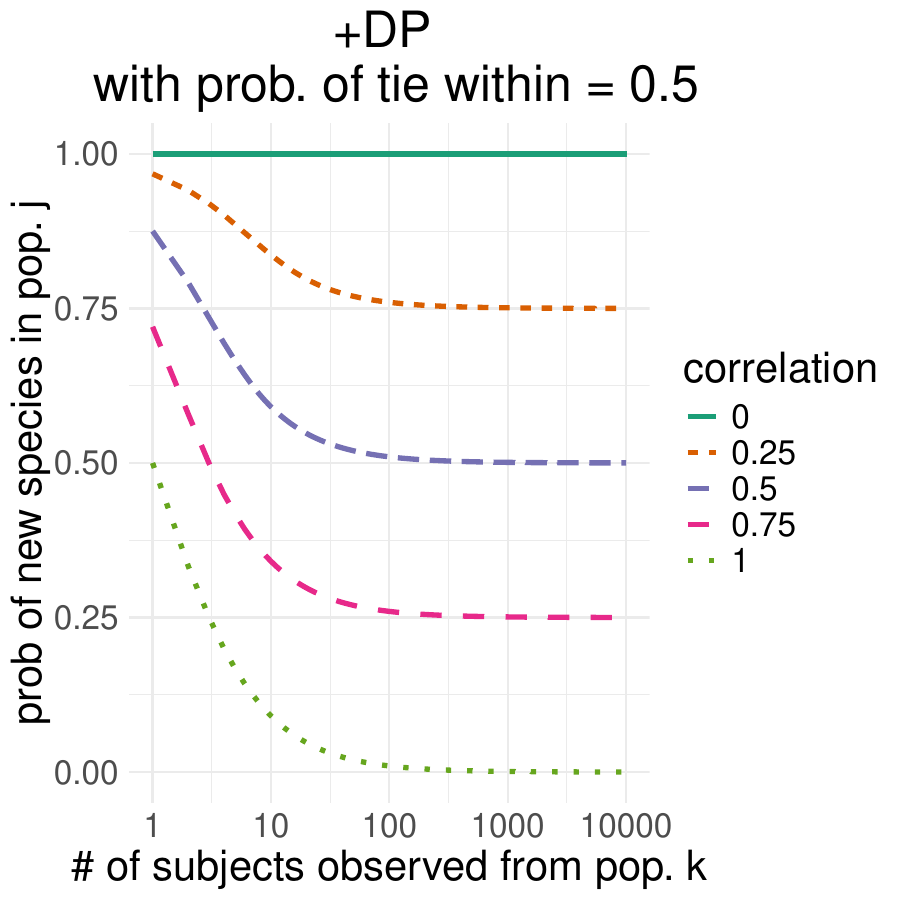}
	\end{subfigure}
 \caption{\label{fig: intro} {Probability of observing a new species at the first draw from population $j$, after $n$ observations from population $k$, i.e., $\bbP(X_{j,1} \notin \{X_{k,1},\ldots, X_{k,n}\})$, plotted as a function of $n$.
 Results are shown (left to right) for hierarchical Dirichlet processes \citep[HDP,][]{teh2006hierarchical}, nested Dirichlet processes \citep[NDP,][]{rodriguez2008nested}, and additive Dirichlet processes \citep[+DP,][]{muller2004method} for varying values of the correlation measure $\Cor[P_{j}(A), P_{k}(A)]$.}} 
\end{figure}

The main goal of these approaches is to flexibly model dependence across populations to facilitate the sharing of information.
This is achieved by defining a collection of dependent latent random probabilities $(P_{1}, \ldots, P_{J})$, which in turn induce dependence among the observables $\bX$.
Clearly, having a way to quantify the dependence between these random probability measures is crucial to understanding and guiding such modeling strategies.
The most widely adopted measure of inter-population dependence is the pairwise correlation between $P_{j}$ and $P_{k}$, for $j\neq k$, evaluated on the same set $A$, namely 
\begin{equation} \label{eq:corr}
	\Cor[P_{j}(A), P_{k}(A)].
\end{equation}
The main reasons this measure has become the benchmark for quantifying dependence are threefold.
First, for most models, it does not depend on the choice of the set $A$, allowing it to be interpreted as a global measure of dependence.
Second, it is often computable either analytically or numerically, which greatly enhances its practical appeal.
Third, although the correlation in \eqref{eq:corr} relies only on the first two moments, it effectively summarizes the dependence structure.
On the one hand, it tends to agree with alternative, more complex measures of dependence that account for the infinite-dimensional nature of the $P_{i}$'s whenever these alternatives can be computed (see \citealp{catalano2021measuring, catalano2024wasserstein}).
On the other hand, different correlation values correspond to markedly different behaviours of the distribution of observable quantities; this is illustrated in Figure~\ref{fig: intro}, which displays the probability of observing a new species for three popular dependent processes as the correlation varies.
However, several fundamental questions remain open, and their resolution is crucial for establishing a rigorous foundation for correlation as a measure of dependence in this context.
First, it is still obscure why \eqref{eq:corr} typically does not depend on the choice of the set $A$; related to this, one would like to understand what conditions on the prior ensure that \eqref{eq:corr} does not depend on the set $A$.
Second, the broader issue of how the properties of the latent random probability measures $P_{i}$'s translate to the observable quantities is still unexplored.
This can be decomposed into two key open questions in terms of correlation: (a) How does the correlation among the latent $P_{i}$'s translate into dependence among observations?
(b) Are correlations among the $P_{i}$'s or observable quantities $\bX$ reliable indicators of dependence?
Specifically, are there models for which $\Cor[P_{j}(A), P_{k}(A)]=1$ if and only if $ P_{j}= P_{k}$ almost surely (i.e., observations are exchangeable), and $\Cor[P_{j}(A), P_{k}(A)]=0$ if and only if $ P_{j} \perp P_{k}$?\\
These questions point to an even more general theme: many structural properties appear to be shared across several classes of dependent processes (e.g., additive, hierarchical, nested, or combinations thereof), and this naturally raises the question of whether it is possible to identify a unifying framework, one that would: (i) encompass most existing models; (ii) allow their structural properties to be studied in a unified way; (iii) provide a foundation for the principled development of new models.

In this work, we provide comprehensive answers to the open problems outlined above.
Specifically, we introduce a new unifying framework for partially exchangeable data that encompasses most existing dependent nonparametric processes.
The resulting class of nonparametric prior processes is termed multivariate species sampling processes (mSSPs), as they play the same foundational role for vectors $(P_{1}, \ldots, P_{J})$ as species sampling processes (SSPs) do in the univariate case.
We characterize mSSPs through their partially exchangeable partition probability function, which encodes the induced multivariate clustering structure.
We show that BNP models currently used for partially exchangeable data belong to a notable subclass of mSSPs, which we refer to as \emph{regular}, which enjoys additional appealing properties.
We also analyze the dependence structure of these processes and prove that borrowing of information across groups is entirely determined by shared ties.
This leads to new insights into the learning mechanisms, offering a principled explanation for the correlation structure discussed above.
Beyond providing a cohesive theoretical foundation, our approach serves as a constructive tool for developing new models and opens new research directions aimed at capturing even richer dependence structures beyond the mSSP framework.\\
Finally, while mSSPs generalize SSPs, their essence lies in their multivariate nature: in particular, in the dependence induced across populations and, consequently, across elements of the vector $(P_{1}, \ldots, P_{J})$.
This feature is obviously absent in classical SSPs, and the structural unification of multivariate, infinite-dimensional objects represents a key innovation of this work.

The paper is structured as follows.
mSSPs and the notable subclass of \emph{regular} mSSPs are defined in Section~\ref{sec: mSSP}.
In Section~\ref{sec: depend}, we derive expressions for marginal and mixed moments of these latent processes in terms of observable quantities.
A substantial part of this section is devoted to analyzing the correlation between the random measures, and we provide the theoretical foundations for it to be regarded as the benchmark measure of dependence within the class of regular mSSPs, where uncorrelation implies independence.
Sections~\ref{sec: pEPPF} and \ref{sec: pred} are devoted to the study of the random partitions induced by mSSPs and their corresponding predictive distributions, respectively.
In Section~\ref{sec: illustration}, we compare the performance of different regular mSSPs in the context of a multi-armed bandit problem aimed at maximizing species discoveries, when sampling sequentially across multiple sites.
Finally, Section~\ref{sec: conclusion} outlines future research directions.
The Supplementary Material contains a review of univariate species sampling processes, all proofs, and further details on the application.
Code to reproduce the experiments is available at \if0\blind the GitHub page of the authors.
\else \url{https://github.com/GiovanniRebaudo/MSSP}.\fi

\section{Multivariate species sampling processes}\label{sec: mSSP}
\subsection{General multivariate species sampling processes}
When extending a univariate random probability $P$ to the multivariate setting involving a vector $(P_{1}, \ldots, P_{J})$ of random probabilities, the species sampling framework of \citet{pitman1996some} can be naturally generalized to multiple populations according to the following definition.
Recall that $\bpi \coloneqq (\pi_{h})_{h\geq1}$ is a sub-probability sequence if $\pi_{h}\in [0,1]$, for all $h$, and $\sum_{h \ge 1} \pi_{h}\leq 1 $.

\begin{definition}
	\label{def: mSSP1}
	A vector of random probability measures $(P_{1}, \ldots, P_{J})$ is a \emph{multivariate species sampling process} (mSSP) if 
	\begin{equation*}
		P_{j} \eqas \sum_{h \ge 1} \pi_{j,h} \delta_{\theta_{h}} + \bigg(1-\sum_{h \ge 1} \pi_{j,h}\bigg) P_{0}, \quad \quad \text{for }j=1,\ldots,J,
	\end{equation*}
	where $P_{0}$ is a non-atomic (deterministic) probability measure on a space $\bbX$, $\bpi_{j} = (\pi_{j,h})_{h\geq1}$ is a random sub-probability sequence, for all $j$, and $\bm{\theta} = (\theta_{h})_{h\geq1}$ are i.i.d.\ from $P_{0}$ and independent of $\bpi = (\bpi_{1}, \ldots, \bpi_{J})\sim \L_{\bpi}$.
	We write $(P_{1}, \ldots, P_{J}) \sim \mSSP(\L_{\bpi}, P_{0})$.
	Moreover, if $\sum_{h \geq 1} \pi_{j,h}\eqas 1$, for any $j$, $(P_{1}, \ldots, P_{J})$ is said to be \emph{proper}.
    By convention $\sum_{h=1}^{0} x_{h} = 0$, for any $(x_{h})$.
\end{definition} 
\noindent According to the standard terminology \citep[see, for instance][]{ghosal2017fundamentals}, we refer to the elements in $\bm{\theta}$ as \emph{atoms}, \emph{labels}, \emph{locations}, or \emph{species}, interchangeably, and to the elements in $\bpi$ as the \emph{weights} of the mSSP.
The structural independence assumption underlying SSPs of independence between the locations and a single sequence of weights, is replaced by independence between the locations $(\theta_{h})_{h\geq1}$ and an array of weights $(\bpi_{1}, \ldots, \bpi_{J})$; the dependence among the $\bpi_{j}$'s then induces a multivariate structure across populations.
From Definition~\ref{def: mSSP1}, the link between mSSPs and SSPs is apparent: it is indeed straightforward to prove that each coordinate of an mSSP is marginally an SSP.
More generally, the class of mSSPs is closed under marginalization.
From now on, for any $n \in \bbN$, let $[n] \coloneqq \{1,\ldots,n\}$.

\begin{proposition}
	\label{prop: marginal mSSP}
	If $(P_{1}, \ldots, P_{J}) \sim \mSSP(\L_{\bpi}, P_{0})$ and $\emptyset \neq \J \subseteq [J]$, then $(P_{j})_{j \in \J} \sim \mSSP(\L_{\bpi_{\J}}, P_{0})$, where $\bpi_{\J} = (\bpi_{j}: j \in \J) \sim \L_{\bpi_{\J}}$ and $\L_{\bpi_{\J}}$ is the marginal distribution obtained from the joint $\L_{\bpi}$.
    Moreover, if $(P_{1}, \ldots, P_{J})$ is proper, then $(P_{j})_{j \in \J}$ is also proper.
\end{proposition}

\noindent At first glance, Definition~\ref{def: mSSP1} might lead one to think of an mSSP as arising from several univariate SSPs sharing the same atoms, i.e., the atoms are given by the same sequence, $(\theta_{h})_{h\geq1}$ for each group.
While this is certainly a possibility, the class of mSSPs is much more general.
In fact, each of the $\pi_{j,h}$'s might be zero almost surely, resulting in mSSPs where the random probabilities $P_{1}, \ldots, P_{J}$ share only a handful or even none of the atoms with positive probability.
Hence, from an intuitive point of view, one should think of an mSSP as arising from several SSPs, which \emph{potentially} share atoms.
Moreover, when investigating the pairwise dependence and quantifying the extent to which species are shared between two random probabilities, being able to distinguish between shared atoms (with positive probability) and idiosyncratic ones is often helpful.
Definition~\ref{def: mSSP1} does not convey this information explicitly, which motivates the need for an alternative representation.
Consider the bivariate mSSP $(P_{1}, P_{2})$, which is obtained marginalizing a $J$-variate mSSP from Definition~\ref{def: mSSP1}.
Then, the following representation is equivalent.
\begin{proposition}\label{prop: characterization bivariate mssp}
    $(P_{1}, P_{2}) \sim \mSSP(\L_{\bpi}, P_{0})$ iff $\exists \, (\tP_{1}, \tP_{2})$ s.t.\ $ (\tP_{1}, \tP_{2}) \eqd (P_{1}, P_{2})$, where
    \begin{equation}\label{eq:3components}
        \tP_{j} \eqas \sum_{h \ge 1} \pi^{(1,2)}_{j,h} \delta_{\theta_{0,h}} + \sum_{h^{\prime} \geq 1 } \pi^{(j)}_{j,h^{\prime}} \delta_{\theta_{j,h^{\prime}}} + \pi^{(j)}_{j,0} P_{0}, \quad \quad \text{for } j=1,2.
    \end{equation}
    with $\bbP\big[ \pi^{(1,2)}_{1,h}\,\pi^{(1,2)}_{2,h}>0 \big]>0$, for every $h$ indexing the first sum, $\big(\pi_{j,h}^{(1,2)}, \pi_{j,h^{\prime}}^{(j)}: h \ge 1, h^{\prime} \ge 0 \big)$, for $j=1,2$, are two random probability vectors and are independent of the atoms, which are distributed as $\theta_{j,h} \simiid P_{0}$ for $j=0,1,2$ and for all $h$.
\end{proposition}

\noindent 
Recall that $\sum_{h=1}^{0} x_{h} = 0$, which allows the representation of $P_{j}$ as a sum of three components to simplify to two or one component.
Importantly, and in contrast to Definition~\ref{def: mSSP1}, we now require $\bbP[\pi^{(1,2)}_{1,h}\,\pi^{(1,2)}_{2,h}>0]>0$, which is clearly equivalent to $\bbP[\pi^{(1,2)}_{1,h}>0, \pi^{(1,2)}_{2,h}>0]>0$, for any $h$ in the first sum.
This condition implies that the atoms in the first sum in \eqref{eq:3components} are common to both random elements with positive probability, thus singling out atoms shared with positive probability.
Conversely, the terms in the second summation are almost surely specific to each $P_{j}$, meaning that the corresponding atoms cannot be shared across the two populations.
Note that the non-atomicity of $P_{0}$ implies $\bbP(\theta_{j,h}=\theta_{k,\ell})=0$ for all $(j,h) \ne (k,\ell)$.
Finally, $\pi^{(j)}_{j,0}$ is the cumulative frequency of almost surely non-shared species that would yield only singletons under an infinite sample from $P_{j}$.
The advantage of the representation in \eqref{eq:3components}, compared to Definition~\ref{def: mSSP1}, lies in the ability to immediately distinguish between shared and non-shared species.
This representation naturally extends to general $J$-variate processes with $J \geq 2$, enabling us to distinguish between species shared across any subset of the $P_{j}$'s.
However, as $J$ grows, the notation becomes increasingly cumbersome, so we omit the explicit formulation here.
More importantly, representation \eqref{eq:3components} enables us to identify a notable subclass of mSSPs, which we term \emph{regular} and examine in the next Section.
Although this subclass encompasses $J$-variate processes $(P_{1}, \ldots, P_{J})$, with $J\geq2$, it suffices to adopt the representation in \eqref{eq:3components} for any pair $(P_{j}, P_{k})$ to fully characterize it.
Finally, we define the pair consisting of a collection of random variables $\bX$ and the mSSP from which these observations are drawn as follows.
\begin{definition}
    \label{def: mSSM}
    A partially exchangeable array $\bX = (X_{j, i}:i\in \bbN, j \in [J])$, for some $J \in \bbN$, follows a multivariate species sampling model (mSSM) if its de Finetti measure is an mSSP.
	That is, for every $j \in [J]$ and for every $i=1,2, \ldots$
	\begin{equation*}
		X_{j,i} \mid (P_{1}, \ldots, P_{J}) \simind P_{j}, \qquad 
		(P_{1}, \ldots, P_{J}) \sim \mSSP(\L_{\bpi}, P_{0}).
	\end{equation*}
\end{definition}

\subsection{Regular multivariate species sampling processes} \label{sec: regular}

A notable subclass of mSSPs, which we term \emph{regular}, arises by imposing a simple independence condition on the weights associated with the non-shared atoms.
First, consider a bivariate mSSP $(P_{1},P_{2})$ and define 
\[
\bpi^{(j)} = \left(\frac{\pi^{(j)}_{j,h^{\prime}}}{\sum_{\ell \geq 0 } \pi^{(j)}_{j,\ell}}\right)_{h^{\prime}\geq0} \quad \text{for } j=1,2
\]
where the weights $\pi^{(j)}_{j,h^{\prime}}$ are as defined in \eqref{eq:3components} and by convention $0/0 = 0$.
\begin{definition}
    \label{def: rmSSP}
    A bivariate mSSP $(P_{1},P_{2})$ is \emph{regular} (rmSSP) if $\bpi^{(1)}\perp \bpi^{(2)}$.
    A $J$-variate mSSP $(P_{1}, \ldots, P_{J})$, with $J>2$, is \emph{regular} if $(P_{j},P_{k})$ is an rmSSP for any $j,k\in[J]$, with $j\neq k$.
\end{definition}
\noindent Intuitively, regularity requires that the relative frequencies of non-shared species are independent for each pair $(P_{j}, P_{k})$ of populations.
Note that if either $P_{j}$, $P_{k}$, or both have no non-shared species in the representation \eqref{eq:3components}, that is, if $\sum_{h^{\prime} \geq 0} \pi^{(\ell)}_{\ell,h^{\prime}} \eqas 0$, then regularity is trivially satisfied.
From a statistical modeling standpoint, it is important to note that the independence condition required by Definition~\ref{def: rmSSP} is relatively mild and, in most applied contexts, reasonable.
This condition implies that the relative frequencies of non-shared species, i.e., $\bpi^{(j)}$, should not influence the sharing of information across groups, once the total frequency $\sum_{h \geq 0} \pi^{(j)}_{j,h}$ of these idiosyncratic species has been accounted for.
Given that these species are (almost surely) not shared among groups, this assumption seems very reasonable.

In the following, we focus specifically on rmSSPs and their use within BNP models, while leaving a broader probabilistic study of general mSSPs for future work.
Special emphasis on the regular subclass is warranted for two main reasons.
First, rmSSPs differ from non-regular mSSPs due to the distinctive dependence structure that enables a remarkable characterization in terms of correlation between the $P_{j}$'s.
This result does not extend to general mSSPs, highlighting a fundamental difference between regular and non-regular mSSPs.
Second, the regular subclass is particularly relevant in statistics, as it encompasses all mSSPs currently employed in BNP.
Some of these will be illustrated in the examples below.

Henceforth, $\DP(\alpha,P_{0})$ denotes the law of a Dirichlet process with concentration parameter $\alpha$ and base measure $P_{0}$ \citep{ferguson1973bayesian}, and $\GEM(\alpha)$ denotes a Griffiths-Engen-McCloskey distribution \citep{sethuraman1994constructive}.
$\PYP(\sigma, \alpha, P_{0})$ stands for the law of a Pitman--Yor process with discount parameter $\sigma$, concentration parameter $\alpha$ and base measure $P_{0}$ \citep{pitman1997two}.
$\mathrm{CRM}(\rho, c, P_{0})$ and $\mathrm{NRMI}(\rho, c, P_{0})$ indicate, respectively, the laws of a completely random measure and a normalized completely random measure with intensity $\rho$, total mass parameter $c$, and base measure $P_{0}$ \citep{regazzini2003distributional, james2009posterior}.
$\GN(\gamma, P_{0})$ denotes the law of a Gnedin process with parameter $\gamma$, and base measure $P_{0}$ \citep{gnedin2010species}, $\DM_{M}(\tau, P_{0})$ is the law of a symmetric Dirichlet multinomial process with $M$ number of categories, concentration parameter $\tau$, and base measure $P_{0}$ \citep{richardson1997bayesian}.
Finally, $\SSP(\L_{\bpi}, P_{0})$ denotes the law of a (univariate) species sampling process with weights distribution defined by $\L_{\bpi}$ and base measure $P_{0}$ \citep{pitman1996some}.
Recall that a concise account of SSPs, including their associated exchangeable partition probability function (EPPF), prediction rule, and notable special cases of SSPs, is provided in Section~\ref{sec: SSP supp} of the Supplement.
In the following, $P_{0}$ always indicates a generic non-atomic deterministic probability measure.

In order to highlight the comprehensive nature of mSSPs and its subclass given by rmSSPs, we now show that several popular classes of dependent models are subclasses of rmSSPs and hence, a fortiori, mSSP.
\begin{example}[Hierarchical processes]\label{ex: hierarchical}
	Assume the distribution of $(P_{1}, \ldots, P_{J})$ coincides with any of the hierarchical specifications listed in Table \ref{tab: HSSP}.
	\begin{table}[H]
		\centering
		\caption{Hierarchical processes \citep{teh2006hierarchical, camerlenghi2019distribution, bassetti2020hierarchical}
			\label{tab: HSSP}}
		\resizebox{\textwidth}{!}{
			\begin{tabular}{lc}
				\toprule
				Hierarchical Dirichlet process (HDP)& 
				$ P_{j}\mid Q\simiid \DP(\alpha,Q), \quad Q\sim \DP(\alpha_{0},P_{0})$ \\
				\midrule
				Hierarchical Pitman--Yor process (HPY) & 
				$P_{j}\mid Q\simiid \PYP( \sigma,\alpha, Q), \quad Q\sim \PYP(\sigma_{0},\alpha_{0}, P_{0})$\\
				\midrule
				Hierarchical normalized completely random measure (HNRMI)&
				$P_{j}\mid Q\simiid \mathrm{NRMI}\left( \rho, c, Q\right), \quad Q\sim \mathrm{NRMI}\left( \rho_{0}, c_{0}, P_{0}\right)$
				\\
				\midrule
				Hierarchical Dirichlet multinomial (HDM):&
				$ P_{j}\mid Q\simiid \DM_{M}(\tau,Q), \quad Q\sim \DM_{M_{0}}(\tau_{0},P_{0})$\\
				\midrule
				Hierarchical Gnedin process (HGN)&
				$ P_{j}\mid Q\simiid \GN(\gamma,Q), \quad Q\sim \GN(\gamma_{0},P_{0})$\\
				\midrule
				Hierarchical species sampling process (HSSP)&
				$ P_{j}\mid Q\simiid \SSP(\L_{\bpi,j},Q), \quad Q\sim \SSP(\L_{\bpi,0},P_{0})$\\
				\bottomrule
			\end{tabular}
		} 
	\end{table}
	\vspace{-0.5\baselineskip}
	\noindent 
	Since $\sum_{h^{\prime} \geq 0} \pi^{(j)}_{j,h^{\prime}} \eqas 0$, for any $j\in[J]$, $(P_{1}, \ldots, P_{J})$ is trivially an rmSSP.
\end{example}

\begin{remark}
SSPs are defined in terms of a non-atomic base measure $P_{0}$ \citep{pitman1996some}.
Hence, writing $\SSP(\L_{\bpi,j}, Q)$ in Table \ref{tab: HSSP} represents an abuse of notation.
However, since the extension to the case where the ``base measure'' $Q$ can be an atomic discrete random measure is immediate, we will also employ it in the sequel.
Furthermore, after marginalizing out $Q$, we are back to the definition of mSSP with a non-atomic base measure.
\end{remark}

\begin{example}[Nested processes]\label{ex: nested}
	Assume the distribution of $(P_{1}, \ldots, P_{J})$ corresponds to any of the nested constructions listed in Table \ref{tab: NSSP}.
	\begin{table}[H]
		\centering
		\caption{Nested processes \citep{rodriguez2008nested, zuanetti2018clustering}
			\label{tab: NSSP}}
		\resizebox{\textwidth}{!}{
			\begin{tabular}{lc}
				\toprule
				Nested Dirichlet process (NDP) & 
				$P_{j}\mid Q\simiid Q, \quad Q\sim \DP(\alpha, \DP(\beta,P_{0}))$ \\
				\midrule
				Nested Pitman--Yor process (NPY) \qquad \qquad \qquad \qquad & 
				$P_{j}\mid Q\simiid Q, \quad Q\sim \PYP(\sigma_{\alpha}, \alpha, \PYP(\sigma_{\beta},\beta,P_{0}))$ \qquad \qquad \\
				\midrule
				Nested Dirichlet multinomial (NDM)&
				$P_{j} \mid Q \simiid Q, \quad Q\sim \DM_{M_{\alpha}}(\tau_{\alpha}, \DM_{M_{\beta}}(\tau_{\beta},P_{0}))$\\
				\midrule
				Nested Gnedin process (NGN)&
				$P_{j}\mid Q\simiid Q, \quad Q\sim \GN(\gamma_{\alpha}, \GN(\gamma_{\beta},P_{0}))$\\
				\midrule
				Nested species sampling process (NSSP)&
				$P_{j}\mid Q\simiid Q, \quad Q\sim \SSP(\L_{\bpi,0}, \SSP(\L_{\bpi},P_{0}))$\\
				\bottomrule
			\end{tabular}
		}
	\end{table}
	\vspace{-0.5\baselineskip}
	\noindent Since $\sum_{h^{\prime} \geq 0} \pi^{(j)}_{j,h^{\prime}} \eqas 0$, for any $j\in[J]$, $(P_{1}, \ldots, P_{J})$ is an rmSSP.
\end{example}

\begin{example}[Additive processes]\label{ex: +SSP}
	Assume the distribution of $(P_{1}, \ldots, P_{J})$ coincides with any of the additive specifications listed in Table \ref{tab: +SSP}.
	\begin{table}[H]
		\centering
		\caption{Additive processes \citep{muller2004method}
			\label{tab: +SSP}}
		\resizebox{\textwidth}{!}{
			\begin{tabular}{lc}
				\toprule
				Additive Dirichlet process (+DP)& 
				$P_{j} = \epsilon_{j} \, Q_{0} + (1-\epsilon_{j}) Q_{j}, \quad Q_{j}\simind \DP(\alpha_{j}, P_{0})$, $j=0,1,\ldots,J$ \\
				\midrule
				Additive Pitman--Yor process (+PY) & 
				$P_{j} = \epsilon_{j} \, Q_{0} + (1-\epsilon_{j}) Q_{j}, \quad Q_{j}\simind \PYP(\sigma_{j},\alpha_{j}, P_{0})$, $j=0,1,\ldots,J$\\
				\midrule
				Additive Dirichlet multinomial (+DM)&
				$P_{j} = \epsilon_{j} \, Q_{0} + (1-\epsilon_{j}) Q_{j}, \quad Q_{j}\simind\DM_{M_{j}}(\tau_{j}, P_{0})$, $j=0,1,\ldots,J$\\
				\midrule
				Additive Gnedin process (+GN)&
				$P_{j} = \epsilon_{j} \, Q_{0} + (1-\epsilon_{j}) Q_{j}, \quad Q_{j}\simind \GN(\gamma_{j}, P_{0})$, $j=0,1,\ldots,J$\\
				\midrule
				Additive species sampling process (+SSP)&
				$P_{j} = \epsilon_{j} \, Q_{0} + (1-\epsilon_{j}) Q_{j}, \quad Q_{j}\simind \SSP(\L_{\bpi,j}, P_{0})$, $j=0,1,\ldots,J$\\
				\bottomrule
			\end{tabular}
		}
	\end{table}
	\vspace{-0.5\baselineskip}
\noindent In this case, the idiosyncratic components are non-zero with positive probability.
However, it is simple to see that for any pair $(P_{j}, P_{k})$ the required independence condition holds, i.e., $\bpi^{(j)}\perp \bpi^{(k)}$.
Thus, $(P_{1}, \ldots, P_{J})$ is an rmSSP.
\end{example}

\begin{example}\citep[Completely random vectors,][]{catalano2021measuring}.
	If $(P_{1}, \ldots, P_{J})$ is distributed according to any of the following: 
	\begin{itemize}
		\item GM-dependent DP \citep[GM-DP,][]{lijoi2014bayesian}:\\ 
		\centerline{$P_{j} = \dfrac{\mu_{0} + \mu_{j} }{\mu_{0}(\bbX) +\mu_{j}(\bbX)}, \quad \mu_{0} \sim \mathrm{CRM}((1-z)\frac{\exp\{-s\}}{s},c, P_{0}), \quad \mu_{j} \simind \mathrm{CRM}(z\frac{\exp\{-s\}}{s},c, P_{0})$; }
		\item GM-dependent $\sigma$-stable \citep[GM-$\sigma$,][]{lijoi2014bayesian}:\\
		\centerline{$P_{j} = \dfrac{\mu_{0} + \mu_{j} }{\mu_{0}(\bbX) +\mu_{j}(\bbX)}, \quad \mu_{0} \sim \mathrm{CRM}((1-z)\frac{\sigma s^{-1-\sigma}}{\Gamma(1-\sigma)},c, P_{0}), \quad \mu_{j} \simind \mathrm{CRM}(z\frac{\sigma s^{-1-\sigma}}{\Gamma(1-\sigma)},c, P_{0})$ }
	\end{itemize}
	then, for any pair $(P_{j},P_{k})$ we have $\bpi^{(j)}\perp \bpi^{(k)}$ and, hence, $(P_{1}, \ldots, P_{J})$ is an rmSSP.
	
	\noindent Moreover, if $(P_{1}, \ldots, P_{J})$ is distributed according to a normalized compound random measures vector \citep[][]{griffin2017compound}, then $\sum_{h^{\prime} \geq 0 } \pi^{(j)}_{j,h^{\prime}} \eqas 0$ for any $j\in[J]$ and, hence, $(P_{1}, \ldots, P_{J})$ is an rmSSP.
\end{example}

\begin{example}[Hidden hierarchical DP]\label{ex: HHDP}
	$(P_{1}, \ldots, P_{J})$ is distributed as a hidden hierarchical Dirichlet process \citep[HHDP,][]{james2008discussion,lijoi2023flexible} if\\
	\centerline{$P_{j}\mid Q\simiid Q, \quad Q\mid Q_{0}\sim \DP(\alpha, \DP(\beta,Q_{0})) \quad Q_{0}\sim \DP(\beta_{0},P_{0})$}\\
	then $\sum_{h^{\prime} \geq 0 } \pi^{(j)}_{j,h^{\prime}} \eqas 0$ holds for any $j\in[J]$, and $(P_{1}, \ldots, P_{J})$ is an rmSSP.
\end{example}
\begin{example}[Semi-hierarchical DP]\label{ex: semi-HDP}
	$(P_{1}, \ldots, P_{J})$ is distributed according to a semi-hierarchical Dirichlet process \citep[semi-HDP,][]{beraha2021semi} if \\
	\centerline{
		$P_{j}\mid Q\simiid Q, \quad Q\mid Q_{0}\sim \DP(\alpha, \DP(\beta,\kappa P_{0} + (1-\kappa)Q_{0})), \quad Q_{0}\sim \DP(\beta_{0},P_{0}).$} Also in this case, 
	we have $\sum_{h^{\prime} \geq 0 } \pi^{(j)}_{j,h^{\prime}} \eqas 0$ for any $j\in[J]$ and $(P_{1}, \ldots, P_{J})$ is an rmSSP.
\end{example} 

\begin{example}[Processes based on stick-breaking constructions]\label{ex: SB}
	Assume the distribution of $(P_{1}, \ldots, P_{J})$ coincides with any of the following: 
	\begin{itemize}
		\item Nested common atoms process \citep[nCAM,][]{denti2023common}, which is given by\\
		\centerline{$P_{j}\mid Q\simiid Q, \, Q = \sum\limits_{s\geq1}\pi_{s}\delta_{G_{s}}, \, G_{s}= \sum\limits_{t \geq 1}\omega_{t,s}\delta_{\theta_{t}}, \, (\pi_{s})_{s\geq1}\sim\GEM(\alpha), \, (\omega_{t,s})_{t \geq 1}\simiid\GEM(\beta)$;}
		\item Tree stick-breaking process with covariates \citep[treeSB,][]{horiguchi2025tree}, which corresponds to\\
		$P_{j} \sim \mathrm{treeSB}(P_{0}, \{F_{j,\epsilon}\}, \mathscr{T})$.
	\end{itemize}
	Since $\sum_{h^{\prime} \geq 0 } \pi^{(j)}_{j,h^{\prime}} \eqas 0$ for any $j\in[J]$, one has that $(P_{1}, \ldots, P_{J})$ is an rmSSP.
\end{example}

\begin{example}[Vectors of normalized independent finite point processes]
	\label{ex: Vec-NIFPP}
    If $(P_{1}, \ldots, P_{J})$ is a vector of finite Dirichlet processes \citep[Vec-FDP,][]{colombi2025hierarchical}, i.e., $(P_{1}, \ldots, P_{J})\sim \mathrm{Vec}-\mathrm{FDP}(\Lambda, \gamma, P_{0})$, then $\sum_{h^{\prime} \geq 0 } \pi^{(j)}_{j,h^{\prime}} \eqas 0$ for any $j\in[J]$ and, thus, $(P_{1}, \ldots, P_{J})$ is an rmSSP.
\end{example}
\begin{example}[Independent processes]\label{ex: IND}
	If $(P_{1}, \ldots, P_{J})$ are independent SSPs, one trivially has $\bpi^{(j)}\perp \bpi^{(k)}$ for any $j\ne k$.
    Thus, $(P_{1}, \ldots, P_{J})$ is an rmSSP.
\end{example}

\section{Dependence structure and moments of mSSPs}\label{sec: depend}
\subsection{Correlation and dependence}\label{subsec: corr}
While mSSPs generalize SSPs, their essence lies in their multivariate nature and investigating the dependence between elements of the vector $(P_{1}, \ldots, P_{J})$ is a crucial task, even more so since this aspect is obviously absent in univariate SSPs.
Here we provide a solid foundation for the use of correlation as a measure of dependence for mSSPs: we derive interpretable expressions for the correlation between pairs of random probability measures in terms of observable variables, prove that the correlation equals one if and only if the data are fully exchangeable, and, furthermore, show how zero correlation characterizes independence for \emph{regular} mSSP.
First, we compute the marginal expected value and variance of the $P_{j}$'s, which will turn out to be helpful in the sequel.
\begin{proposition}
	\label{prop: mean var PA}
	If $(P_{1}, \ldots, P_{J})$ is an mSSP and $X_{j,i}\mid(P_{1}, \ldots, P_{J})\simind P_{j}$, for $i=1,2,\ldots$ and $j=1,\ldots, J$, then 
	\[
	\bbE[P_{j}(A)] = P_{0}(A),
	\qquad
	\Var[P_{j}(A)]= \bbP(X_{j,1}=X_{j,2}) P_{0}(A)\big[1-P_{0}(A)\big].
	\]
\end{proposition}
\noindent Note that by marginal exchangeability, the tie probability $\bbP(X_{j, i}=X_{j,m})$ between observations from the same population $j$ does not depend on the indices $(i,m)$.
Moreover, using representation \eqref{eq:3components}, one can rewrite the tie probability of an mSSP as 
\begin{equation}
	\label{eq:tie_within}
	\bbP(X_{j,1}=X_{j,2}) = \sum_{h \ge 1} \bbE\left[ \left(\pi^{(j,k)}_{j,h}\right)^{2}\right] + \sum_{h^{\prime}\geq 1} \bbE\left[ \left(\pi^{(j)}_{j,h^{\prime}}\right)^{2} \right].
\end{equation}
The link between tie probability and variance of a single homogeneous NRMI was first noted in \textcite{james2006conjugacy}.
Here we have established it for general mSSPs: since mSSPs do not require specifying a law for the weights, this means that the link between variance and tie probability is structural.
However, this represents only our starting point in uncovering the crucial role played by tie probabilities for mSSPs.\\ 
The following simple, yet important, step consists of looking not only at the tie probabilities within each population but also across.
Also in this case, the tie probability across populations $j$ and $k$, $\bbP(X_{j, i}=X_{k,m})$, does not depend on the indices $(i,m)$ and equals
\begin{equation}
	\label{eq:tie_across}
	\bbP(X_{j,1}=X_{k,1}) = \sum_{h \ge 1} \bbE\left[ \pi^{(j,k)}_{j,h}\, \pi^{(j,k)}_{k,h}\right].
\end{equation}
In Section~\ref{sec: pEPPF}, we will also express the tie probability in terms of the more general law of the partition induced at the level of the observables.
We are now ready to compute the correlation of mSSPs, a major highlight of this paper in terms of both understanding dependent models and methodological implications.
\begin{proposition}
    \label{prop: CorrProb}
    Let $(P_{1}, \ldots, P_{J})$ be an mSSP and $X_{j,i}\mid(P_{1}, \ldots, P_{J}) \simind P_{j}$, for $i=1,2,\ldots$ and $j=1,\ldots, J$.
    Then we have
	\begin{align*}
		\Cor[P_{j}(A),P_{k}(A)]
		&= \frac{\bbP(X_{j,1}=X_{k,1})}
		{ \sqrt{\bbP(X_{j,1}=X_{j,2})} \sqrt{\bbP(X_{k,1}=X_{k,2})} } \qquad \forall j\neq k\in[J],
	\end{align*}
    for any measurable set $A$ such that $\Var[P_{\ell}(A)]>0$, for $\ell=j,k$.
\end{proposition}
\noindent This key result has several intertwined ramifications.
First, it solves the open problem of identifying the reason for the correlation not to depend on the evaluation set $A$, which was observed on a case by case basis in most currently employed models: by Proposition~\ref{prop: CorrProb} for any mSSP, the pairwise correlation between its elements is expressed exclusively in terms of tie probabilities within and across groups; hence, by \eqref{eq:tie_within}--\eqref{eq:tie_across} it depends only on the weights of the mSSP and it cannot depend on the set $A$, which is characterized in terms of locations.
Second, Proposition~\ref{prop: CorrProb} implies that correlation between random probabilities is a consequence uniquely of ties between the observable species.
Thus, the dependence boils down to ties across populations and the learning mechanism runs exclusively through the ties.
Further insights on correlation as a measure of global dependence are collected in the following corollary.
\begin{corollary}
    \label{cor: CorrProb}
    Let $(P_{1}, \ldots, P_{J})$ be an mSSP and $X_{j,i}\mid(P_{1}, \ldots, P_{J}) \simind P_{j},$ for $i=1,2,\ldots$ and $j=1,\ldots, J$.
    Then, for any $j\neq k\in[J]$ and any measurable set $A$ such that $\Var[P_{\ell}(A)]>0$, for $\ell=j,k$, we have
	\begin{enumerate}[(c-i)]
		\item $\Cor[P_{j}(A),P_{k}(A)] \ge 0$;
		\item $\Cor[P_{j}(A),P_{k}(A)] = 0$ iff $\bbP(X_{j,1}=X_{k,1})=0$ iff $\bbE(\pi_{j,h}\pi_{k,h}) = 0$, for all $h$;
		\item If $\bbP(X_{j,1}=X_{j,2})=\bbP(X_{k,1}=X_{k,2})>0$ (e.g., if $P_{j}$ and $P_{k}$ are equal in distribution), then $\Cor[P_{j}(A),P_{k}(A)] = \bbP(X_{j,1}=X_{k,1})/\bbP(X_{j,1}=X_{j,2})$.
	\end{enumerate}
\end{corollary}

\begin{remark}
    The third statement is particularly appealing from an intuitive standpoint: in the common situation of equal marginals, one can think of correlation as the ratio of the probabilities of, respectively, ``tie across groups'' and ``tie within a group''.
    Implicitly, this also ensures that the tie probability ``across groups'' is always smaller than, or equal to, the one ``within a group'', which is a reasonable and appealing feature.
    See also \textcite{durante2025partially} for a discussion in the context of multilayer networks, where this ordering can be recast as a desirable generalized homophily property.
\end{remark}

\noindent Corollary~\ref{cor: CorrProb} links correlation among pairs of $P_{j}$'s with properties of observable quantities, providing both an intuitive and rigorous foundation to the use of correlation as a measure of dependence.
However, this leaves an important question unanswered: how well does the correlation capture dependence among different processes?
Are the extreme situations of full exchangeability, i.e., maximal dependence, and independence recovered when the correlation equals one and zero, respectively?
The next proposition shows that a correlation equal to one implies maximal dependence, that is, full exchangeability, of the observables.

\begin{theorem}
\label{thm: Corr1}
    Let $(P_{1}, \ldots, P_{J})$ be an mSSP and $X_{j,i}\mid(P_{1}, \ldots, P_{J})\simind P_{j},$ for $i=1,2,\ldots$ and $j=1,\ldots, J$.
    Let $A$ be a measurable set such that $0 < P_{0}(A) < 1$.
    If, for $j\neq k\in[J]$, at least one of the following holds 
    \begin{itemize}
        \item[{(p-i)}] $(P_{j},P_{k})$ is proper;
        \item[{(p-ii)}] $\bbP(X_{j,1} = X_{j,2}) = \bbP(X_{k,1} = X_{k,2})>0$,
    \end{itemize}
    then 
	\[
	\Cor[P_{j} (A), P_{k} (A)]= 1 \quad\text{if and only if}\quad P_{j}\eqas P_{k}
	\]
	and if and only if $\bX = (X_{\ell,i},\,i\geq1, \, \ell\in\{j,k\})$ is exchangeable.
\end{theorem}
\noindent Corollary~\ref{cor: CorrProb} and Theorem~\ref{thm: Corr1} jointly provide a straightforward interpretation of what happens when the probability of a tie across groups approaches the probability of a tie within: the correlation increases towards one and the dependence among the observations shifts from partial exchangeability towards the extreme of full exchangeability.

The other extreme case, namely independence, is harder to recover from a situation of zero correlation.
However, if we restrict attention to rmSSPs, we are able to show that it is impossible to construct zero-correlated rmSSPs whose components are not pairwise independent.
This yields the desired characterization, but also highlights the natural role played by rmSSPs within the general class of mSSP.
\begin{theorem}\label{thm: CorInd}
	Let $(P_{1}, \ldots, P_{J})$ be an rmSSP.
    Then, for any $j\neq k\in[J]$ and any measurable set $A$ such that $\Var[P_{\ell}(A)]>0$, for $\ell=j,k$, we have
	\[
	   \Cor[P_{j}(A), P_{k}(A)]= 0 \quad\text{if and only if} \quad P_{j} \perp P_{k}.
	\]
\end{theorem}
\noindent Hence, within the class of rmSSP, on one hand, correlation equal to zero implies independence among the $P_{j}$'s and across groups of observations, and, on the other hand (under assumption (p-i) or (p-ii)),  correlation equal to one implies exchangeability.

\begin{remark}
	Not all rmSSPs can achieve correlation exactly equal to zero or one, at least in their standard definitions.
	To fix ideas, consider rmSSPs lacking idiosyncratic and improper components (i.e., $\sum_{h \ge 1} \pi^{(i,j)}_{j,h} \eqas 1$); one remarkable instance is given by hierarchical constructions.
    For these processes, the situation of independence across groups is to be interpreted as a limiting case.
	For example, $J$ independent DPs arise from the HDP only if we let $\alpha_{0}\to\infty$.
    To make things concrete and highlight how much literature we cover with mSSPs, Table~\ref{tab: corr} presents the correlation, probability of ties, and the values of hyperparameters to attain independence and exchangeability for a wide variety of models.
\end{remark}

\noindent From an inferential perspective, the dependence among the latent $(P_{1}, \ldots, P_{J})$ plays a key instrumental role, since it induces dependence among the observations.
We have already recovered the extreme cases of exchangeability and independence of the observables as corresponding to, respectively, correlation one and zero of pairs of $P_{j}$'s.
Nonetheless, the following results, which hold for the entire class of mSSPs, highlight how the correlation among observables coincides with the tie probability.
Further, we stress the implications on the induced dependence among the data $\bX$.
\begin{proposition}
    \label{prop: CorrObs}
    Let $(P_{1}, \ldots, P_{J})\sim \mSSP(\L_{\bpi},P_{0})$ and $X_{j,i}\mid(P_{1}, \ldots, P_{J})\simind P_{j},$ for $i=1,2,\ldots$ and $j=1,\ldots, J$.
    Assume $\bbX = \bbR$ and $P_{0}$ has finite second moment.
    Then, for any $j, k\in[J]$ and any $i,m$, 
	\[
	\Cor(X_{j,i},X_{k,m}) = \bbP(X_{j,i}=X_{k,m}).
	\]
\end{proposition}

\begin{table}[!ht]
    \centering
    \caption{Correlation, tie probabilities and extreme cases.
    From left to right: type of mSSP (notation defined in Examples of Section~\ref{sec: regular}); pairwise correlation; probability of tie across and within groups; values to which the hyperparameters should converge for the correlation to converge, respectively, to $0$ and $1$ (while $P(\text{Ties Within})$ does not converge to $0$ or $1$).	
		\label{tab: corr}
	} 
	\resizebox{\textwidth}{!}{
		\begin{tabular}{lccccc}
			\toprule
			Process & Correlation & $P(\text{Ties Across})$ & $P(\text{Ties Within})$ & Indep.\ & Exch.\\
			\midrule
			HDP & $\dfrac{1+\alpha}{1+\alpha+\alpha_{0}}$ & $\dfrac{1}{1+\alpha_{0}}$ & $\dfrac{1+\alpha+\alpha_{0}}{(1+\alpha)\,(1+\alpha_{0})}$& $\alpha_{0}\rightarrow +\infty$ & $\alpha\rightarrow +\infty$
			\rule{0pt}{5ex} \\ 
			HPY
			& $\dfrac{(1+\alpha)(1-\sigma_{0})}{(1 - \sigma\sigma_{0}) +\alpha (1-\sigma_{0})+\alpha_{0}(1-\sigma)}$
			& $\dfrac{1-\sigma_{0}}{1+\alpha_{0}}$
			& $\dfrac{(1 - \sigma\sigma_{0}) +\alpha (1-\sigma_{0})+\alpha_{0}(1-\sigma)}{(1+\alpha)\,(1+\alpha_{0})}$
			& $\displaystyle {\genfrac{}{}{0pt}{}{\alpha_{0}\rightarrow +\infty}{\text{or }\sigma_{0}\rightarrow 1}}$ 
			& $\displaystyle {\genfrac{}{}{0pt}{}{\alpha\rightarrow +\infty}{\text{or }\sigma\rightarrow 1}}$ \rule{0pt}{5ex} \\
			HDM 
			& $\dfrac{(1+\tau_{0})(1+\tau\,M)}{(1 + \tau\,M)(1+\tau_{0}\,M_{0}) - \tau\tau_{0}(M-1)(M_{0}-1)}$
			& $\dfrac{1+\tau_{0}}{1+\tau_{0}\,M_{0}}$ 
			&$\dfrac{(1 + \tau\,M)(1+\tau_{0}\,M_{0}) - \tau\tau_{0}(M-1)(M_{0}-1)}{(1 + \tau\,M)(1+\tau_{0}\,M_{0})}$
			&$M_{0}\rightarrow +\infty$
			&$M\rightarrow +\infty$\rule{0pt}{5ex} \\
			HGN 
			& $\dfrac{\gamma_{0} (\gamma+1)}{(\gamma + \gamma_{0})}$
			& $\dfrac{2\gamma_{0}}{\gamma_{0}+1}$
			& $\dfrac{2(\gamma + \gamma_{0})}{(\gamma + 1)(\gamma_{0}+1)}$
			& $\gamma_{0}\rightarrow 0$
			& $\gamma\rightarrow 0$\rule{0pt}{5ex} \\
			HSSP 
			& $\dfrac{\EPPF_{1,0}^{(2)}(2)}{\EPPF_{1,1}^{(2)}(2) + \EPPF_{2,1}^{(2)}(1,1) \EPPF_{1,0}^{(2)}(2)}\;^{\star}$
			& $\EPPF_{1,0}^{(2)}(2)$
			& $\EPPF_{1,1}^{(2)}(2) + \EPPF_{2,1}^{(2)}(1,1) {\EPPF_{1,0}^{(2)}(2)}$
			& $\EPPF_{1,0}^{(2)}(2) =0$
			& $\EPPF_{1,1}^{(2)}(2)=0$ \rule{0pt}{5ex} \\
			NDP
			& $\dfrac{1}{1+\alpha}$ 
			& $\dfrac{1}{(1+\alpha)(1+\beta)}$ 
			& $\dfrac{1}{1+\beta}$ 
			& $\alpha\rightarrow +\infty$ 
			& $\alpha \rightarrow 0$ \rule{0pt}{5ex} \\
			NPY
			& $\dfrac{1-\sigma_{\alpha}}{1+\alpha}$ 
			& $\dfrac{(1-\sigma_{\alpha})(1-\sigma_{\beta})}{(1+\alpha)(1+\beta)}$ 
			& $\dfrac{1-\sigma_{\beta}}{1+\beta}$ 
			& $\displaystyle {\genfrac{}{}{0pt}{}{\alpha\rightarrow +\infty }{\text{or }\sigma_{\alpha}\rightarrow 1}} $ 
			& $\displaystyle {\genfrac{}{}{0pt}{}{(\alpha,\sigma_{\alpha})\rightarrow}{ \rightarrow (0,0)}}$ \rule{0pt}{5ex}\\
			NDM
			& $\dfrac{1+\tau_{\alpha}}{1+\tau_{\alpha}M_{\alpha}}$
			& $\dfrac{(1+\tau_{\alpha})(1+\tau_{\beta})}{(1+\tau_{\alpha}M_{\alpha})(1+\tau_{\beta}M_{\beta})}$
			&$\dfrac{1+\tau_{\beta}}{1+\tau_{\beta}M_{\beta}}$
			&$M_{\alpha}\rightarrow +\infty$
			&$M_{\alpha}\rightarrow 1$ \rule{0pt}{5ex} \\
			NGN
			&$\dfrac{2\gamma_{\alpha}}{\gamma_{\alpha} + 1}$ 
			&$\dfrac{4\gamma_{\alpha}\gamma_{\beta}}{(\gamma_{\alpha} + 1)(\gamma_{\beta} + 1)}$
			&$\dfrac{2\gamma_{\beta}}{\gamma_{\beta} + 1}$
			&$\gamma_{\alpha} \rightarrow 0$
			&$\gamma_{\alpha} \rightarrow 1$ \rule{0pt}{5ex} \\
			NSSP
			&$\EPPF^{(2)}_{1,0}(2)$
			&$\EPPF^{(2)}_{1,0}(2)\EPPF^{(2)}_{1,1}(2)$
			&$\EPPF^{(2)}_{1,1}(2)$
			&$\EPPF^{(2)}_{1,0}(2) = 0$
			&$\EPPF^{(2)}_{1,0}(2) = 1$
			\rule{0pt}{5ex}\\
			+DP 
			& $\dfrac{\dfrac{\epsilon_{j}\epsilon_{k}}{1+\alpha_{0}}}{\sqrt{\left(\dfrac{\epsilon_{j}^{2}}{1+\alpha_{0}} + \dfrac{(1-\epsilon_{j})^{2}}{1+\alpha_{j}}\right)\left(\dfrac{\epsilon_{k}^{2}}{1+\alpha_{0}} + \dfrac{(1-\epsilon_{k})^{2}}{1+\alpha_{k}}\right)}}$ 
			& $ \dfrac{\epsilon_{j} \epsilon_{k}}{1+\alpha_{0}}$ 
			& $\dfrac{\epsilon_{j}^{2}}{1+\alpha_{0}} + \dfrac{(1-\epsilon_{j})^{2}}{1+\alpha_{j}}$ 
			& $\epsilon_{j}=0 \text{ or } \epsilon_{k} = 0$
			&$\epsilon_{j} = \epsilon_{k} =1$\rule{0pt}{5ex}\\
			+PY
			& $\dfrac{\dfrac{\epsilon_{j}\epsilon_{k}\, (1-\sigma_{0})}{1+\alpha_{0}}}{\sqrt{\left(\dfrac{\epsilon_{j}^{2}\, (1-\sigma_{0})}{1+\alpha_{0}} + \dfrac{(1-\epsilon_{j})^{2}\, (1-\sigma_{j})}{1+\alpha_{j}}\right)\left(\dfrac{\epsilon_{k}^{2}\, (1-\sigma_{0})}{1+\alpha_{0}} + \dfrac{(1-\epsilon_{k})^{2}\, (1-\sigma_{k})}{1+\alpha_{k}}\right)}}$ 
			& $ \dfrac{\epsilon_{j}\epsilon_{k}\, (1-\sigma_{0})}{1+\alpha_{0}}$ 
			& $\dfrac{\epsilon_{j}^{2}\, (1-\sigma_{0})}{1+\alpha_{0}} + \dfrac{(1-\epsilon_{j})^{2}\, (1-\sigma_{j})}{1+\alpha_{j}}$ 
			& $\epsilon_{j}=0 \text{ or } \epsilon_{k} = 0$
			&$\epsilon_{j} = \epsilon_{k} =1$\rule{0pt}{5ex}\\
			+DM
			&$\dfrac{\dfrac{\epsilon_{j}\epsilon_{k}(1+\tau_{0})}{1+\tau_{0}\,M_{0}}}{\sqrt{\left(\dfrac{\epsilon_{j}^{2}(1+\tau_{0})}{1+\tau_{0}\,M_{0}} + 
					\dfrac{(1-\epsilon_{j})^{2}(1+\tau_{j})}{1+\tau_{j}\,M_{j}}\right)\left(\dfrac{\epsilon_{k}^{2}(1+\tau_{0})}{1+\tau_{0}\,M_{0}} + 
					\dfrac{(1-\epsilon_{k})^{2}(1+\tau_{k})}{1+\tau_{k}\,M_{k}}\right)}}$ 
			&$\dfrac{\epsilon_{j}\epsilon_{k}(1+\tau_{0})}{1+\tau_{0}\,M_{0}}$
			&$\dfrac{\epsilon_{j}^{2}(1+\tau_{0})}{1+\tau_{0}\,M_{0}} + 
			\dfrac{(1-\epsilon_{j})^{2}(1+\tau_{j})}{1+\tau_{j}\,M_{j}}$
			& $\epsilon_{j}=0 \text{ or } \epsilon_{k} = 0$
			&$\epsilon_{j} = \epsilon_{k} =1$\rule{0pt}{5ex}\\
			+GN
			& $\dfrac{\dfrac{\epsilon_{j}\epsilon_{k}\,2\gamma_{0}}{\gamma_{0}+1}}{\sqrt{\left(\dfrac{\epsilon_{j}^{2}\,2\gamma_{0}}{\gamma_{0}+1}+
					\dfrac{(1-\epsilon_{j})^{2}\,2\gamma_{j}}{\gamma_{j}+1}\right)\left(\dfrac{\epsilon_{k}^{2}\,2\gamma_{0}}{\gamma_{0}+1}+
					\dfrac{(1-\epsilon_{k})^{2}\,2\gamma_{k}}{\gamma_{k}+1}\right)}}$ 
			& $\dfrac{\epsilon_{j}\epsilon_{k}\,2\gamma_{0}}{\gamma_{0}+1}$
			& $\dfrac{\epsilon_{j}^{2}\,2\gamma_{0}}{\gamma_{0}+1}+
			\dfrac{(1-\epsilon_{j})^{2}\,2\gamma_{j}}{\gamma_{j}+1}$ 
			& $\epsilon_{j}=0 \text{ or } \epsilon_{k} = 0$
			&$\epsilon_{j} = \epsilon_{k} =1$\rule{0pt}{5ex}\\
			+SSP
			& $\dfrac{\epsilon_{j}\epsilon_{k}\EPPF^{(2)}_{1,0}(2)}{\sqrt{(\epsilon_{j}^{2} \EPPF^{(2)}_{1,0}(2)+
					(1-\epsilon_{j})^{2} \EPPF^{(2)}_{1,1}(2))(\epsilon_{k}^{2} \EPPF^{(2)}_{1,0}(2)+
					(1-\epsilon_{k})^{2} \EPPF^{(2)}_{1,1}(2))}}$
			& $\epsilon_{j}\epsilon_{k}\EPPF^{(2)}_{1,0}(2)$ 
			& $\epsilon_{j}^{2} \EPPF^{(2)}_{1,0}(2)+
			(1-\epsilon_{j})^{2} \EPPF^{(2)}_{1,1}(2)$
			& $\epsilon_{j}=0 \text{ or } \epsilon_{k} = 0$
			&$\epsilon_{j} = \epsilon_{k} =1$\rule{0pt}{5ex}\\
			GM-DP & $\dfrac{(1-z)c}{1+c}\,{}_{3}F_{2}(a,1,1;b,b;1)^{\star\star}$& $\dfrac{(1-z)c}{(1+c)^{2}}\,{}_{3}F_{2}(a,1,1;b,b;1)$ & $\dfrac{1}{1+c}$ & $z=1$ & $z=0$ \rule{0pt}{5ex} \\
			GM-$\sigma$ & $(1-z) \I(c,z)^{\star\star\star}$ & $(1-z)(1-\sigma) \I(c,z)$ & $1-\sigma$& $z=1$ & $z=0$ \rule{0pt}{5ex}\\
			HHDP & $1 - \dfrac{\alpha \beta_{0}}{(1+\alpha)(\beta_{0} +\beta +1)}$ & $\dfrac{1}{\beta_{0}+1}+\dfrac{\beta_{0}}{(1+\alpha)(1+\beta)(1+\beta_{0})}$ &$\dfrac{1+\beta+\beta_{0}}{(1+\beta)\,(1+\beta_{0})}$ & $\displaystyle{ \genfrac{}{}{0pt}{}{(\alpha, \beta_{0}) \rightarrow} {\rightarrow(+ \infty, + \infty)}}$&$\alpha\rightarrow 0$ \rule{0pt}{5ex}\\
			nCAM
			&$1-\dfrac{\beta\alpha}{(2\beta+1)(1+\alpha)}$
			&$\dfrac{1}{1+\alpha}\left(\dfrac{1}{1+\beta}+ \dfrac{\alpha}{2\beta+1}\right)$
			&$\dfrac{1}{1+\beta}$
			&None
			&$\alpha \rightarrow 0$\rule{0pt}{5ex}\\
			\bottomrule
			\multicolumn{6}{l}{{\small
					${}^{\star}$ $\EPPF_{\cdot,1}^{\cdot}$ and $\EPPF_{\cdot,0}^{\cdot}$ are 
					induced by $\L_{\bpi,1}=\ldots=\L_{\bpi,J}$ and $\L_{\bpi,0}$, respectively.
                    $\qquad$ ${}^{\star\star}$ ${}_{3}F_{2}$ is the generalized hypergeometric function and $a = \alpha(1-z)+2, b = \alpha+2$ $\qquad$ 
					${}^{\star\star\star}$ $\I(c,z) = \dfrac{1}{\sigma} \int_{0}^{1} \dfrac{\omega^{1/\sigma - 1}}{[1+z(1-\omega^{1/\sigma})^{\sigma} - z(1-\omega)]}\d \omega$}} 
	\end{tabular}}
\end{table}
\begin{corollary}\label{cor: CorrX}
	Let $(P_{1}, \ldots, P_{J})\sim \mSSP(\L_{\bpi},P_{0})$ and $X_{j,i}\mid(P_{1}, \ldots, P_{J})\simind P_{j},$ for $i=1,2,\ldots$ and $j=1,\ldots, J$.
    Assume $\bbX = \bbR$ and $P_{0}$ has finite second moment.
    Then, for any $j,k\in[J]$ and any $i\neq m$, 
	\begin{enumerate}[(c-i)]
		\item $\Cor(X_{j,i},X_{k,m}) \ge 0$;
		\item $\Cor(X_{j,i},X_{k,m}) = 0$ if and only if $\bbP(X_{j,i} = X_{k,m})=0$ if and only if $X_{j,i} \perp X_{k,m}$;
		\item $\Cor(X_{j,i},X_{k,m}) = 0$ if and only if $\bbE(\pi_{j,h}\pi_{k,h}) = 0$, for all $h$.
	\end{enumerate}
\end{corollary}

\noindent Note that Proposition~\ref{prop: CorrObs} and Corollary~\ref{cor: CorrX} hold true both within (i.e., $j=k$) and across (i.e., $j\neq k$) groups, and thus, also for (univariate) SSP.

\subsection{Higher moments of mSSPs}
We now derive both marginal and mixed moments of any order.
These can also be seen as generalizations, to all SSPs and mSSPs, of the powerful results on joint moments of normalized completely random measures \citep{james2006conjugacy} and of hierarchical normalized completely random measures \citep{camerlenghi2019distribution}, which leverage the Laplace functional characterization of completely random measures.
Here we show that moments can be computed in the class of mSSPs even for elements unrelated to completely random measures and/or to hierarchical processes.
The following two propositions provide the expressions for the marginal moments.
\begin{proposition}
	\label{prop: MomSSP}
	Let $(P_{1}, \ldots, P_{J})$ be an mSSP, $X_{j,i}\mid(P_{1}, \ldots, P_{J})\simind P_{j},$ for $i=1,2,\ldots$ and $j=1,\ldots, J$.
    Then, for every $q \in \bbN$ and measurable set $A$,
	\[
	\bbE[P_{j}(A)^{q}] = \bbE\bigg[P_{0}(A)^{K^{(j)}_{1:q}} \bigg],
	\]
	where $K^{(j)}_{1:q}$ is the random number of unique species in a sample of size $q$ from $P_{j}$.
\end{proposition}

\newpage
\begin{proposition}
	\label{prop: JointMomSSP}
	Let $(P_{1}, \ldots, P_{J})$ be an mSSP and $\{A_{1},\ldots, A_{h}\}$ be pairwise disjoint measurable sets.
    Then, for any sequence $q_{1},q_{2},\ldots,q_{h}$, with $q_{i} \in \bbN$ for $i=1 \ldots, h$, we have
   \begin{align*}
		&\bbE[P_{j}(A_{1})^{q_{1}} \cdots P_{j}(A_{h})^{q_{h}}] = \bbE \left[P_{0}(A_{1})^{K^{(j)}_{1:q_{1}}}\, P_{0}(A_{2})^{K^{(j)}_{q_{1}+1:q_{1}+q_{2}}} \cdots P_{0}(A_{h})^{K^{(j)}_{q_{h-1}+1: q_{1} + \cdots + q_{h}}} \mid E_{\neq} \right] \bbP(E_{\neq}),
	\end{align*}
	where $K^{(j)}_{a:b}$ is the random number of species in the ``block of observations'' from the $a$-th to the $b$-th observation, in a sample of size $q_{1} + \cdots + q_{h}$ from $P_{j}$, and $E_{\neq}$ is the event that no shared species are recorded across the different blocks of observations.
\end{proposition}
\noindent The following two theorems provide the expressions for the mixed moments.
\vspace{-0.5\baselineskip}
\begin{theorem}
	\label{thm: MommSSP}
	Let $(P_{1}, \ldots, P_{J})$ be an mSSP, $X_{j,i}\mid(P_{1}, \ldots, P_{J})\simind P_{j},$ for $i=1,2,\ldots$ and $j=1,\ldots, J$.
    Then, for any measurable set $A$ and sequence $q_{1},q_{2},\ldots,q_{J}$, with $q_{i} \in \bbN$ for $i=1 \ldots, J$, we have
	\[
    \bbE[P_{1}(A)^{q_{1}} \cdots P_{J}(A)^{q_{J}}] = \bbE\big[ P_{0}(A)^{K_{q_{1}, \ldots, q_{J}}}\big],
	\]
	where $K_{q_{1},\ldots,q_{J}}$ is the overall number of species observed in a sample that contains $q_{j}$ observations from $P_{j}$, for $j=1,\ldots, J$.
\end{theorem}

\begin{theorem}
	\label{thm: JointMommSSP}
	Let $(P_{1}, \ldots, P_{J})$ be an mSSP and $\{A_{1},\ldots, A_{J}\}$ be
	pairwise disjoint measurable sets.
    Then, for any sequence $q_{1},q_{2},\ldots,q_{J}$, with $q_{i} \in \bbN$ for $i=1 \ldots, J$, we have
	\begin{align*} 
		&\bbE[P_{1}(A_{1})^{q_{1}} \cdots P_{J}(A_{J})^{q_{J}}]= \bbE\left[ P_{0}(A_{1})^{K^{(1)}_{1:q_{1}}} \cdots P_{0}(A_{J})^{K^{(J)}_{1:q_{J}}}\mid E_{\neq}\right] \bbP(E_{\neq}),
	\end{align*}
	where $K^{(j)}_{1:q_{j}}$ is the number of observed species from population $j$, in a sample which contains $q_{j}$ observations from $P_{j}$, for $j=1,\ldots,J$ and $E_{\neq}$ is the event that no species is shared across the $J$ groups.
\end{theorem}
\noindent 
Importantly, these results showcase that higher-order moments of an mSSP evaluated on (measurable) sets can be meaningfully interpreted in terms of simple observable quantities like the random number of observed species within and across groups.
Hence, interpretability is not unique to correlation.

\section{Partially exchangeable partition 
function}
\label{sec: pEPPF}
In the exchangeable case, the random partition induced by a discrete nonparametric prior is characterized by the exchangeable partition probability function (EPPF), a concept introduced in \textcite{pitman1995exchangeable}.
The EPPF plays a fundamental role across several domains, including combinatorics, stochastic process theory, population genetics, Bayesian statistics and machine learning; see \textcite{pitman2006combinatorial} and references therein.
It also underpins models in ecology and natural language processing, where exchangeable clustering structures naturally arise.
The EPPF associated with the DP corresponds to Ewens' sampling formula \citep{antoniak1974mixtures, ewens1972sampling}; see \textcite{crane2016ubiquitous, tavare2021magical} for accounts of its widespread applications.

\noindent Moving from exchangeability to partial exchangeability, the counterpart of the EPPF is the partially exchangeable partition probability function (pEPPF), which characterizes the random partition induced by a partially exchangeable array.
This concept is different from the one introduced in \textcite{pitman1995exchangeable}, which is unrelated to the original definition of partial exchangeability in the sense of de Finetti that we adopt here.
The notion of pEPPF first appeared in \textcite{leisen2011vectors} and \textcite{lijoi2014bayesian} for specific instances of dependent processes, and it started being leveraged in a systematic way for other subclasses of dependent processes only recently in, e.g., \textcite{camerlenghi2017bayesian, camerlenghi2019distribution, camerlenghi2019latent, beraha2021semi, lijoi2023flexible, denti2023common}.
Its absence from the classical probabilistic literature may stem from the fact that, unlike in the exchangeable setting where the EPPF can often be defined directly without invoking an associated exchangeable sequence, the partially exchangeable case lacks a similarly tractable direct construction.
Instead, the pEPPF arises naturally by marginalizing a partially exchangeable array of random elements, which represents the canonical approach to deriving the corresponding random partition within the BNP framework.
In multi-population species sampling problems, where the values sampled from $P_{0}$ serve solely as species labels with no numerical meaning, the pEPPF uniquely determines the marginal likelihood of the observations.
Similarly, in the context of model-based clustering or latent multi-level modeling, the pEPPF encapsulates the underlying clustering mechanism.
Importantly, from a computational perspective, the pEPPF also plays a pivotal role as it provides the key ingredient for deriving marginal posterior sampling schemes.

Let us now formally introduce the pEPPF induced by any vector $(P_{1}, \ldots, P_{J})$ of random probability measures with possibly discrete components.
A sample $(X_{j, i}: i \in [I_{j}], j \in [J])$, where $I_{j}$ denotes the sample size of group $j$ and $n = \sum_{j=1}^{J} I_{j}$ is the total sample size, induces a random partition of the integers $[n]$ based on the ties among the observations.
To see this, let the integers label the observations according to their \emph{order of arrival by group}, that is, observations are indexed first by group $j = 1, \ldots, J$, and then by within-group order of arrival.
Specifically, observation $X_{j,i}$ is associated with the label $\sum_{k = 1}^{j-1} I_{k} + i$, for any $i=1,\ldots,I_{j}$.
The resulting random partition can be usefully characterized by the corresponding pEPPF.
To this end, let $D$ be the number of distinct values among the $n=\sum_{j=1}^{J} I_{j}$ observations in the sample $(X_{j,i}: i \in [I_{j}], \, j \in [J])$.
For each group $j$, define the vector of frequency counts $\bn_{j}=(n_{j,1},\ldots,n_{j, D})$, where $n_{j,d}$ indicates the number of observations in the $j$th group that coincide with the $d$th distinct value, indexed according to the order of arrival by groups.
Clearly, $n_{j,d} \ge 0$ and by construction $\sum_{j=1}^{J} n_{j,d} \ge 1$, since each distinct value must appear in at least one group.
The count $n_{j,d}=0$ indicates that the $d$th distinct value does not occur in group $j$, while it is shared between groups $k$ and $l$ if and only if $n_{k,d} \, n_{l,d} \ge 1$.
The law of the resulting random partition is characterized by its pEPPF, defined as 
\begin{equation}\label{eq: pEPPF}
	\pEPPF_{D}^{(n)}(\bn_{1}, \ldots, \bn_{J}) = \bbE \bigg[ \int_{\bbX^{D}_{\star}} \prod_{d=1}^{D} P_{1}(\d x_{d})^{n_{1,d}} \ldots P_{J}(\d x_{d})^{n_{J,d}} \bigg],
\end{equation}
under the constraint that $\sum_{d=1}^{D} n_{j,d} = I_{j}$ for each $j = 1, \ldots, J$, and where $\bbX$ denotes the space in which the $X_{j,i}$'s take values, while $\bbX^{D}_{\star}$ denotes the subset of $\bbX^{D}$ consisting of vectors with all distinct entries.
We stress that the expectation in \eqref{eq: pEPPF} is taken with respect to the joint distribution of the vector of random probability measures $(P_{1}, \ldots, P_{J})$, that is, the de Finetti measure associated with the partially exchangeable array.
An important special case is immediately recovered when $J=1$, namely the single population setting: indeed, the pEPPF in \eqref{eq: pEPPF} reduces to a standard EPPF.
Moreover, if $J=2$, the probability of a tie across groups coincides with $\pEPPF_{1}^{(2)}(1,1)$.

Clearly, if $(P_{1}, \ldots, P_{J})$ is an mSSP, it induces a pEPPF as defined by \eqref{eq: pEPPF}.
But what about the converse?
Given a pEPPF, does there exist an mSSP, up to the choice of an independent, non-atomic base measure $P_{0}$, that generates it?
This amounts to asking whether every pair of pEPPF and independent non-atomic $P_{0}$ determines a unique mSSP.
The next result provides an affirmative answer and, as a by-product, yields an intuitive generative construction.

\begin{theorem}\label{thm: pEPPFchar}
	Consider an arbitrary pEPPF as in \eqref{eq: pEPPF} and let $P_{0}$ be a non-atomic (deterministic) probability measure.
    Consider the partially exchangeable array $\bX=(X_{j,i}:\: i\in \bbN,\: j\in[J])$ such that for any non-negative integers $I_{1},\ldots,I_{J}$ the variables $(X_{j, i}: i \in [I_{j}], \, j \in [J])$ follow the generative scheme:
	\begin{enumerate}
		\item sample a random partition $\Pi_{n}$ from the given pEPPF;
		\item conditionally on the partition $\Pi_{n}$, sample from $P_{0}$ the i.i.d.\ unique values associated with each partition set.
	\end{enumerate}
    Then, the de Finetti measure associated with $\bX$ is the law of an mSSP.
\end{theorem}

\begin{remark}
At first glance, the previous result may seem surprising: any pEPPF in \eqref{eq: pEPPF}, regardless of whether $(P_{1}, \ldots, P_{J})$ generating it is an mSSP or not, identifies an mSSP by pairing it with an independent non-atomic base measure $P_{0}$.
The core idea behind such a fundamental result is that the distribution of the weights of the directing probability measure is what characterizes the partition induced by the ties of any (infinite) partially exchangeable array.
Hence, if the dependent process $(P_{1}, \ldots, P_{J})$ inducing the pEPPF is not an mSSP (for instance, due to dependence between weights and locations), one can still identify an mSSP, say $(P_{1}^{\star},\ldots, P_{J}^{\star})$, based on the same pEPPF.
The two vectors $(P_{1}, \ldots, P_{J})$ and $(P_{1}^{\star},\ldots, P_{J}^{\star})$ will generally have different distributions, although they share the same pEPPF.
In other words, the pEPPF characterizes the multivariate partition structure, but not the law of the partially exchangeable array.
Recovering the latter requires specifying a mechanism for atom assignment, with the independence required by mSSPs representing a simple and tractable choice.
Crucially, this implies that the random partition structure of any partially exchangeable array can be studied via the pEPPF of an mSSP, making mSSPs the natural framework for analyzing and understanding the discrete structure of arbitrary dependent vectors and partially exchangeable partition models.
\end{remark}

The notion of pEPPF also enables us to restate the pairwise correlation results from Section~\ref{subsec: corr} within its more general structure.
For instance, one can express 
\[
\Cor[P_{j}(A),P_{k}(A) ] = \frac{\pEPPF_{1}^{(2)} (1,1)}{\sqrt{\EPPF_{j,1}^{(2)} (2)} \sqrt{\EPPF_{k,1}^{(2)}(2)}},
\]
where $\EPPF_{j}$ denotes the marginal EPPF corresponding to $P_{j}$.
Similarly, the correlation between observations satisfies $\Cor(X_{j,i}, X_{k,m}) = \pEPPF_{1}^{(2)} (1,1)$.

Finally, we record an alternative representation of the pEPPF in terms of the weights associated with a proper mSSP.
\begin{proposition}
	\label{prop: peppf from weight}
	Let $(P_{1}, \ldots, P_{J})$ be a proper mSSP.
	Then
	\begin{equation*}
		\pEPPF_{D}^{(n)}(\bn_{1},\ldots,\bn_{J}) = \bbE \bigg[ \sum_{h_{1} \ne \cdots \ne h_{D}} \prod_{j=1}^{J} \prod_{d=1}^{D} \pi_{j,h_{d}}^{n_{j,d}} \bigg].
	\end{equation*}
\end{proposition}

\section{Predictive structure and inference}\label{sec: pred}
In the exchangeable case, the predictive distribution of an SSP admits a simple and elegant representation, with weights expressed as ratios of the associated EPPF \citep{pitman1996some}.
The sequential mechanism that generates these prediction rules is known as the generalized Chinese restaurant process: observations correspond to customers entering a restaurant, each choosing to sit either at an already occupied table or at a new one.
Each table serves a unique dish drawn independently from $P_{0}$.
The predictive distribution is given by
\begin{equation*}
	\bbP \big(X_{n+1}=x \mid \bX \big)=\begin{cases}
		\quad \frac{\EPPF_{K}^{(n+1)}(n_{1},\ldots,n_{k}+1, \ldots, n_{K})}{\EPPF_{K}^{(n)}(n_{1},\ldots,n_{k}, \ldots, n_{K})} &\qquad \text{if } x=X_{k}^{\star} \quad \text{ and } k=1, \ldots, K\\[5pt]
		\quad \frac{\EPPF_{K+1}^{(n+1)}(n_{1},\ldots,n_{k}, \ldots, n_{K},1)}{\EPPF_{K}^{(n)}(n_{1},\ldots,n_{k}, \ldots, n_{K})} &\qquad \text{if } x=X_{K+1}^{\star},
	\end{cases}
\end{equation*}
where $(X_{k}^{\star} : k = 1, \ldots, K)$ are the $K$ distinct values observed among $X_{1}, \ldots, X_{n}$, appearing with frequencies $(n_{1}, \ldots, n_{K})$ and drawn i.i.d.\ from $P_{0}$.
These probabilities follow from conditioning on the observed partition: the EPPF in the numerator is updated either by increasing the count of an existing cluster or by adding a new singleton cluster.
See Section~\ref{sec: SSP supp} of the Supplementary Material for further details.

In the partially exchangeable framework, the pEPPF associated with an mSSP naturally gives rise to a multivariate generative mechanism, where predictive distributions are again expressed as ratios of pEPPFs.
We refer to this construction as the multivariate generalized Chinese restaurant process (mgCRP).
It differs both from the classical generalized Chinese restaurant process and from common multi-population extensions typically modeled as restaurant franchises.
Unlike the latter, we do not introduce multiple restaurants.
Unlike the former, although we retain a single restaurant in which each table serves a unique dish and a customer at a new table receives a previously unserved dish, the allocation mechanism is more intricate.
Specifically, the probability that a customer sits at a given table depends not only on the current seating configuration, but also on the group of the incoming customer and the group membership of those already seated.
Remarkably, these allocation probabilities can still be expressed as ratios of pEPPFs, a fact that is both natural and striking.
The resulting mgCRP is formalized in the next proposition.

\begin{proposition}\label{prop: urn from peppf}
    Let $\bX$ be a partially exchangeable array directed by a de Finetti measure given by the law of an mSSP $(P_{1}, \ldots, P_{J})$.
    For any $j\in[J]$, the corresponding predictive distributions are characterized by an mgCRP of the form
	\begin{align*}
		X_{j, I_{j}+1}\mid (\bX_{j, 1:I_{j}})_{j=1}^{J} = 
		\begin{cases} \quad
			X^{\star}_{l} &\text{w.p.\ } \frac{\pEPPF^{(n+1)}_{D}(\bn_{1}, \ldots, [n_{j,1},\ldots,n_{j,l}+1,\ldots, n_{j,D}], \ldots, \bn_{J})} {\pEPPF^{(n)}_{D}(\bn_{1}, \ldots,[n_{j,1},\ldots,n_{j,l},\ldots, n_{j,D}], \ldots, \bn_{J})} \\[5pt]
			\quad X^{\star}_{new} &\text{w.p.\ } \frac{\pEPPF^{(n+1)}_{D+1}( [\bn_{1},0], [n_{j,1}, \ldots, n_{j,D},1], \ldots, [\bn_{J},0])}{\pEPPF^{(n)}_{D}(\bn_{1}, \ldots, [n_{j,1},\ldots, n_{j,D}], \ldots, \bn_{J})}
		\end{cases}
	\end{align*} 
	where $(X_{1}^{\star}, \ldots, X_{D}^{\star})$ are the $D$ unique values in $(\bX_{j, 1:I_{j}})_{j=1}^{J}$ listed in order of arrival by group, $n = \sum_{j} I_{j}$, and $X^{\star}_{new}$ represents a new species sampled independently from $P_{0}$.
\end{proposition}

\noindent Although the predictive scheme in Proposition~\ref{prop: urn from peppf} follows naturally from the structure of the pEPPF and stands out for its theoretical elegance, its computational feasibility depends heavily on the ability to evaluate ratios of pEPPFs.
Unlike the case of univariate SSPs, where such ratios are sometimes available in closed form, the multivariate setting rarely admits simple analytic expressions.
Exceptions are limited to trivial cases that reduce to univariate specifications, such as independent or almost surely identical Gibbs-type priors.
Section~\ref{sec: SSP supp} details explicit predictive schemes for specific univariate SSPs within the Gibbs-type family \citep{gnedin2006exchangeable,lijoi2007bayesian}, arguably the most tractable generalization of the DP \citep{deblasi2015gibbs}.
Nevertheless, it is important to note that mSSPs used in Bayesian modeling give rise to pEPPFs that are obtained as mixtures of EPPFs.
This includes models such as the HSSP, NSSP, +SSP, and various combinations thereof.
Thus, implementable predictive sampling schemes can typically be derived through data augmentation strategies that exploit the tractability of the underlying EPPFs.
These techniques leverage latent variables that simplify ratios of pEPPFs to ratios of products of EPPFs, greatly simplifying computations.
A prominent example is the Chinese restaurant franchise representation for the HDP \citep{teh2006hierarchical}.
The following examples present such augmented formulations of the pEPPF, which enable tractable predictive schemes and facilitate the design of marginal Gibbs samplers, for three large classes of regular mSSPs, namely HSSP, NSSP, and +SSP.
Let $\pEPPF^{(n)}_{D,\mathrm{aug}}(\bn_{1}, \ldots, \bn_{J}, \bm{\ell}, \bq)$ denote the augmented pEPPF, in the sense that the original pEPPF can be recovered by summing over all possible values of the latent variables $\bm{\ell}$ and $\bq$, i.e., $\pEPPF^{(n)}_{D}(\bn_{1}, \ldots, \bn_{J})= \sum_{\bm{\ell}, \bq}\pEPPF^{(n)}_{D,\mathrm{aug}}(\bn_{1}, \ldots, \bn_{J}, \, \bm{\ell}, \bq)$.
\begin{repexample}{ex: hierarchical}
	If $(P_{1}, \ldots, P_{J})$ is an HSSP, then
	\begin{align}\label{eq: aug pEPPF hier}
		\resizebox{0.83\textwidth}{!}{
			$
			\pEPPF^{(n)}_{D,\mathrm{aug}}(\bn_{1}, \ldots, \bn_{J}, \, \bm{\ell}, \bq) = \EPPF_{D,0}^{(\ell_{\cdot,\cdot})}(\ell_{\cdot,1}, \ldots, \ell_{\cdot,D}) \prod_{j=1}^{J}
			\EPPF_{\ell_{j,\cdot},j}^{(I_{j})}(q_{j,1}, \ldots, q_{j,\ell_{j, \cdot}}),
			$
		}
	\end{align}
	where for $j=1,\ldots,J$, $\EPPF_{\ell_{j,\cdot},j}^{(I_{j})}(q_{j,1}, \ldots, q_{j,\ell_{j, \cdot}})$ denotes the EPPF induced by ${\L_{\bpi,j}}$, which characterizes a latent partition of the $I_{j}$ observations of group $j$ into $\ell_{j,\cdot}$ blocks of cardinalities $q_{j,1}, \ldots, q_{j,\ell_{j, \cdot}}$.
	Conditionally on these partitions, all the $\ell_{\cdot, \cdot} = \sum_{j=1}^{J} \ell_{j, \cdot}$ blocks (we use ‘$\cdot$’ to indicate summation over that index) are grouped into a coarser partition of $D$ blocks, each corresponding to a distinct observed species.
    The distribution of this coarser partition is characterized by the $\EPPF_{D,0}^{(\ell_{\cdot,\cdot})}(\ell_{\cdot,1}, \ldots, \ell_{\cdot, D})$ induced by $Q$.
\end{repexample}

\begin{repexample}{ex: nested}
	If $(P_{1}, \ldots, P_{J})$ is an NSSP, then 
	\begin{equation}\label{eq: pEPPF aug NSSP}
		\resizebox{0.83\textwidth}{!}{
			$
			\pEPPF^{(n)}_{D,\mathrm{aug}}(\bn_{1}, \ldots, \bn_{J}, \, \bm{\ell}, \, \bq) =
			\EPPF_{R,0}^{(J)}(\ell_{1},\ldots,\ell_{R}) \prod_{r=1}^{R} \EPPF_{D_{r}}^{(I^{\star}_{r})}(q_{r,1},\ldots,q_{r, D_{r}}),
			$
		}
	\end{equation}
	where $\EPPF_{R,0}^{(J)}(\ell_{1},\ldots,\ell_{R})$ denotes the EPPF induced by $\L_{\bpi,0}$ that controls the clustering of the group labels $j=1,\ldots,J$ into $R$ blocks (obtained from the ties among the $P_{j}$'s).
	Let $P_{r}^{\star} \simiid \SSP(\L_{\bpi},P_{0})$ $r=1,\ldots,R$ be the unique values of $(P_{j})_{j=1}^{J}$ in order of arrival and let $I_{r}^{\star}=\sum_{j: P_{j}=P_{r}^{\star}}I_{j}$ be the number of observations from the $\ell_{r}$ groups assigned to $P_{r}^{\star}$.
	Conditionally on this clustering of the groups, for $r=1,\ldots, R$, the $\EPPF_{D_{r}}^{(I^{\star}_{r})}(q_{r, 1},\ldots,q_{r, D_{r}})$ describes the distribution induced by $P_{r}^{\star}$ characterizing the partition of the $I_{r}^{\star}$ observations assigned to $P_{r}^{\star}$ into $D_{r}$ distinct species.
	Since the $P_{r}^{\star}$'s do not share atoms a.s., it follows that the total number of distinct species is given by $D=\sum_{r=1}^{R} D_{r}$.
\end{repexample}

\begin{repexample}{ex: +SSP}
	If $(P_{1}, \ldots, P_{J})$ is a +SSP, then 
	\begin{equation}\label{eq: pEPPF aug +SSP}
		\pEPPF^{(n)}_{D,\mathrm{aug}}(\bn_{1}, \ldots, \bn_{J}, \, \bm{\ell}, \, \bq) =
		\prod_{j=1}^{J} \epsilon_{j} ^{\ell_{0,j}} (1-\epsilon_{j})^{\ell_{j}} \prod_{j=0}^{J}
		\EPPF^{(\ell_{j})}_{D_{j},j}(q_{j,1},\ldots,q_{j,D_{j}}),
	\end{equation}
	where
	$\ell_{0,j}$ and $\ell_{j}=I_{j}-\ell_{0,j}$ denote, for each $j\in[J]$, the number of observations assigned to the shared SSP $Q_{0}$ and to the idiosyncratic SSP $Q_{j}$, respectively.
    Here, $\epsilon_{j}^{\ell_{0,j}} (1 -\epsilon_{j})^{\ell_{j}}$ is the probability of the i.i.d.\ latent assignment of the $I_{j}$ observations via $\mathrm{Bern}(\epsilon_{j})$.
    Moreover, $\ell_{0} = \sum_{j=1}^{J}\ell_{0,j}$ is the total number of observations assigned to the shared SSP $Q_{0}$.
    Conditionally on these latent assignments, for $j=0,\ldots,J$, $\EPPF^{(\ell_{j})}_{D_{j},j}(q_{j,1},\ldots,q_{j,D_{j}})$ is the EPPF induced by $Q_{j}$ that governs the clustering of the $\ell_{j}$ observations assigned to $Q_{j}$ into $D_{j}$ unique species.
	Since the $Q_{j}$'s do not share species a.s., the total number of distinct species across all groups is given by $D=\sum_{j=0}^{J} D_{j}$.
\end{repexample}

The hierarchical representations of the pEPPF derived in \eqref{eq: aug pEPPF hier}, \eqref{eq: pEPPF aug NSSP}, and \eqref{eq: pEPPF aug +SSP}, expressed as products of simple EPPFs in an augmented space, allow for simplifying the ratio in Proposition~\ref{prop: urn from peppf} into a product of tractable predictive expressions, analogous to those of the Chinese restaurant process.

\section{Multi-armed bandits for species discovery}\label{sec: illustration}

Among the numerous application areas of dependent processes, which include density regression, spatio-temporal analysis, functional data, survival analysis, topic modeling, hierarchical and multi-level clustering, and ANOVA-type models, here we focus on a multi-armed bandit problem connected to species sampling.
Embedding the analysis within the framework of mSSMs offers a principled way to compare different classes of models.
Our structural results for mSSMs allow model parameters to be chosen so that relevant prior quantities coincide, enabling fair performance comparisons among competing approaches.
A systematic investigation of which models are preferable in specific settings is beyond the scope of this paper, but the general strategy is clear: calibrate prior parameters to match the aspects that are most critical for the task at hand, and evaluate inferential performance accordingly.
This will be the focus of future work.

\noindent A Bayesian nonparametric approach to species sampling problems in the single population case, i.e., $J=1$, was introduced by \textcite{lijoi2007bayesian}, where Bayesian analogs of the classical Turing and Good-Toulmin estimators \citep{good1953population,good1956number} were derived.
In this setting, a random probability measure $P$ models the species proportions in the population, and, given an observed sample, the main goal is to estimate the probability of discovering a new species either at the next step or after an additional $m$ unobserved draws.
Following \textcite{lijoi2007bayesian}, there has been a rich literature exploring alternative prior specifications, estimation of diverse functionals and quantities of interest, and a wide range of applications.
For detailed reviews, see \textcite{deblasi2015gibbs} and \textcite{balocchi2026bayesian}, and references therein.

The multi-sample setup with $J$ populations modeled through a vector of dependent random probability measures $(P_{1}, \ldots, P_{J})$ was first studied in \textcite{camerlenghi2017bayesian}.
A sequential perspective was adopted in \textcite{battiston2018multi, camerlenghi2020nonparametric}, to design an optimal sequential sampling strategy to maximize the diversity of the observed species.
This involves deciding, at each step, which population to sample from, while sequentially incorporating information from previously observed species across populations.
The problem naturally fits within the framework of multi-armed bandits, where each arm represents a population and a unit reward is earned upon discovering a new species.
Such problems arise in ecology and biology, where sampling from diverse environments aims to uncover new species, and in genomics, where the objective is often to detect as many genetic variants as possible \citep[see, e.g.,][]{lijoi2008bayesian, masoero2022more}.

\subsection{Real data}
\begin{figure}
    \begin{subfigure}{0.6\linewidth}
		\centering
		\includegraphics[width=\textwidth]{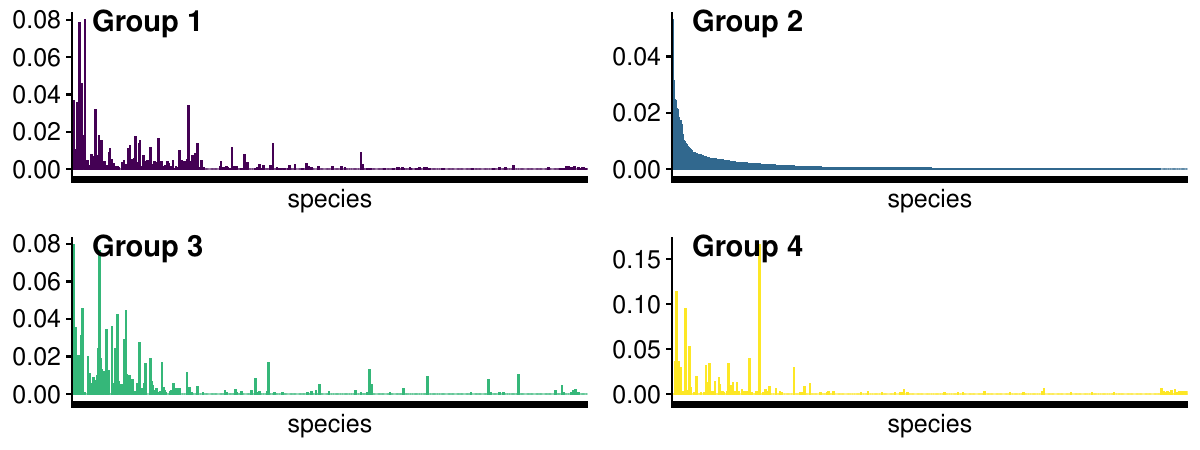}
	\end{subfigure}
 \hfill
 \begin{subfigure}{0.38\linewidth}
		\centering
		\includegraphics[ width=\linewidth]{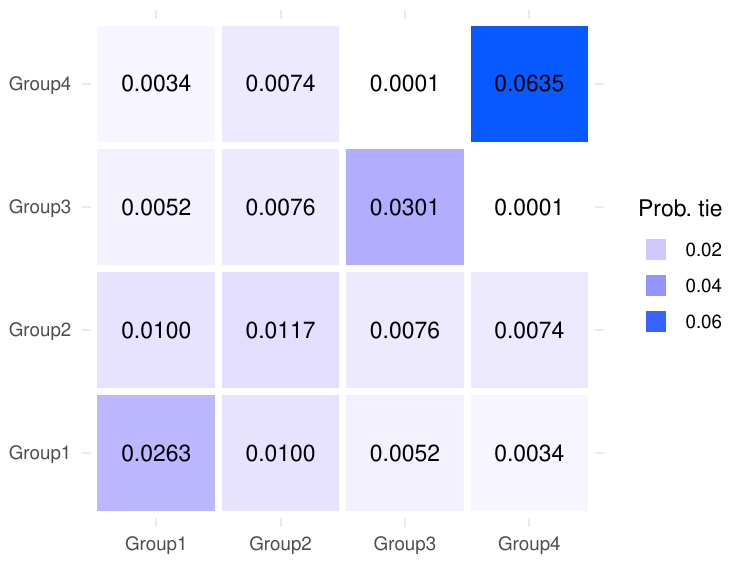}
	\end{subfigure}
    \caption{
    Empirical distribution functions of the four groups (left) and empirical tie probabilities (right), computed from the full dataset.
    In the left panel, species within each group are ordered according to their frequency ranking in Group $2$.}
    \label{fig: data}
\end{figure}

Here we consider a multi-armed bandit problem of tree-species discovery, using the dataset of South American tree species publicly available in the supplementary materials of \textcite{condit2002beta}.
The dataset records $41,688$ trees observed across $100$ plots, comprising $802$ distinct species.
In accordance with \textcite{battiston2018multi}, we aggregated the $100$ plots into four larger groups based on spatial location, joining columns in the dataset whose codes begin with BCI, P, S, and C.
These four groups define the $J=4$ alternative arms.
The empirical distributions and empirical tie probabilities (i.e., relative frequencies) for each group are shown in Figure~\ref{fig: data}.
Further details on the dataset can be found in \textcite{pyke2001floristic, condit2002beta, battiston2018multi}.
\noindent The four arms corresponding to the trees' populations in the four different regions are modeled as a vector of dependent random probabilities $(P_{1}, \ldots, P_{4})$, each representing a population whose species composition is initially unknown, both in terms of presence of a species and relative abundance.
Species may be shared across arms, possibly with different frequencies, making the rmSSP framework a natural modeling choice.

\noindent Within this setup, we compare the performance of six rmSSPs in maximizing the number of distinct observed species in an additional sample.
For each model, the sampling strategy to achieve this goal consists of selecting at step $n+1$ the arm with the highest estimated probability of discovering a new species, based on observations collected up to step $n$.
Specifically, at each step we choose the arm $j$ that maximizes $\bbP(X_{j, I_{j}+1} \notin \bX_{\text{past}} \mid \bX_{\text{past}})$, where $\bX_{\text{past}}$ denotes the previously observed species across all sites.
We also contrast these model-based approaches with a simple baseline that selects an arm uniformly at random at each step, which we refer to as the uniform model.

\noindent The six rmSSP models we compare are: independent DP and PYP, additive DP and PYP, and hierarchical DP and PYP.
We assign hyperpriors to the concentration parameter in each DP-based rmSSP and to both concentration and discount parameters in each PYP-based rmSSP.
To ensure a fair comparison, these hyperpriors are chosen so that the prior mean and variance of the tie probabilities, within groups for all models and across groups for those that borrow information, match across all six specifications.
This calibration reflects two considerations.
First, in a species sampling problem, where our objective is to maximize the species diversity, it is sensible to set the probability of not discovering a new species, i.e., the probability of ties, equal across all models.
Second, by the results of the previous section, the tie probabilities effectively capture the dependence structure and information-sharing behaviour for any rmSSP, regardless of the specific application at hand.
Full details on model definitions, sampling algorithms, and hyperprior settings are provided in Section~\ref{sec: algo supp} of the Supplementary Material.

\noindent Figure~\ref{fig: real} showcases the average cumulative number of species discovered by the two hierarchical, the two additive, and the two independent rmSSP models, each also compared to the uniform model, as a function of the number of additional samples.
Table~\ref{tab: real} reports the average number of new species discovered per sampling step.
All results are averages based on $20$ runs.
In each run, we begin with an initial sample of $30$ observations per arm (drawn without replacement from the full dataset), then sample $300$ further observations sequentially according to each strategy and record the species discoveries.
\begin{figure}[tb]
	\begin{subfigure}{0.32\linewidth}
		\centering
		\includegraphics[ width=\linewidth]{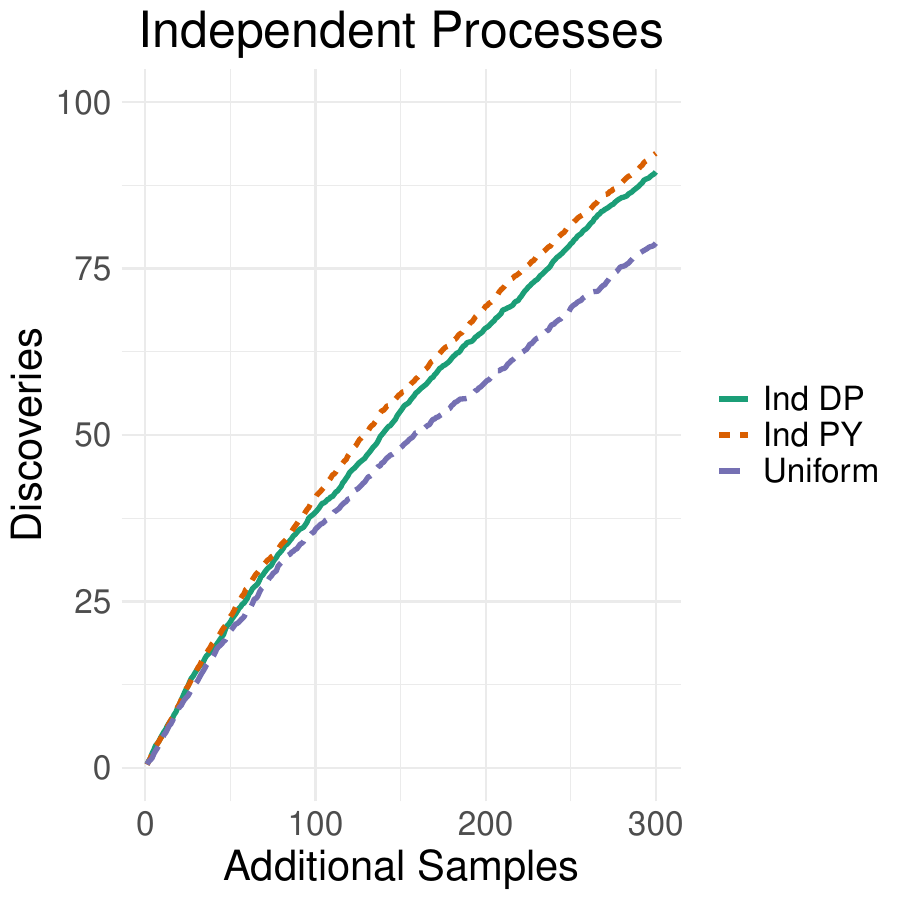}
	\end{subfigure}
 \hfill
 \begin{subfigure}{0.32\linewidth}
		\centering
		\includegraphics[ width=\linewidth]{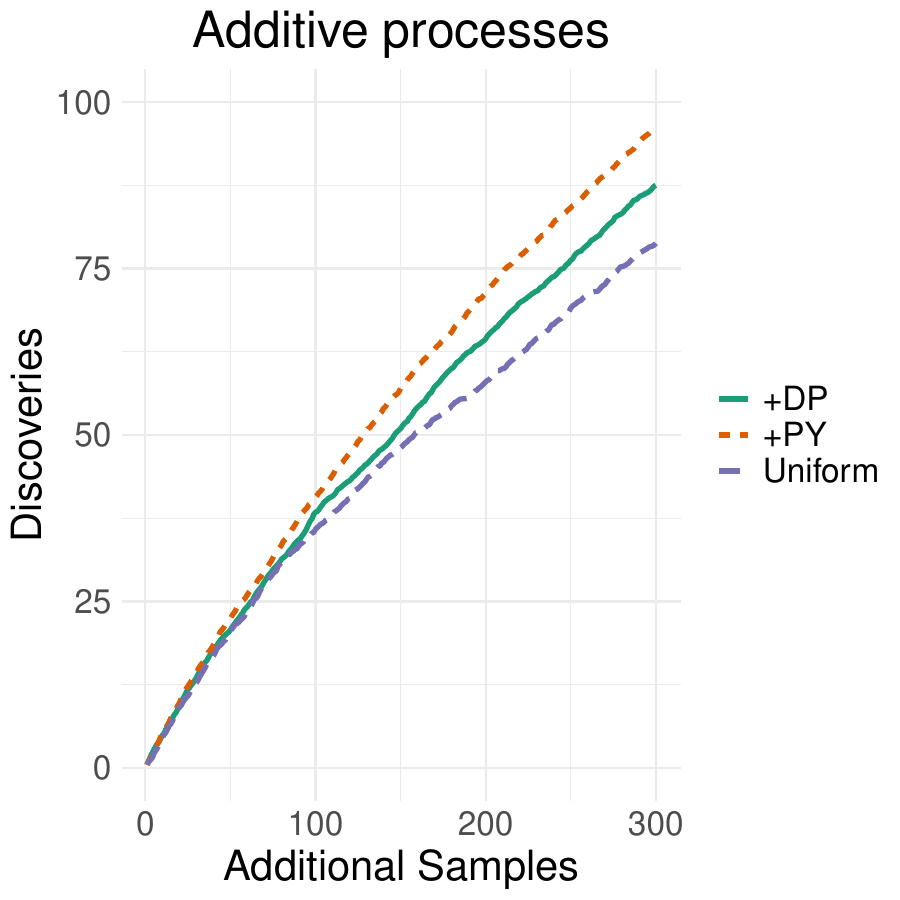}
	\end{subfigure}
 \hfill
		\begin{subfigure}{0.32\linewidth}
		\centering
		\includegraphics[ width=\linewidth]{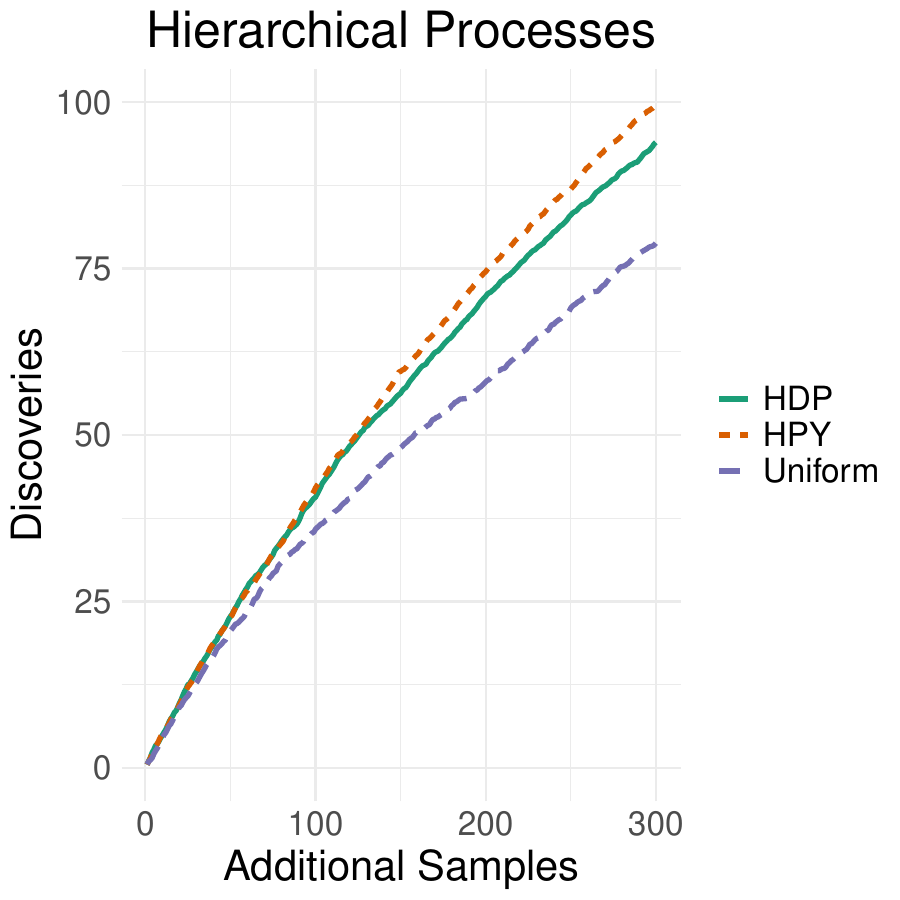}
	\end{subfigure}
 \caption{\label{fig: real} 
 	Tree species data: cumulative number of species discovered as a function of the additional sample size for each rmSSP model and the uniform baseline.}
\end{figure}

\begin{table}[tb]
    \centering
    \renewcommand{\arraystretch}{1} 
    \resizebox{\textwidth}{!}{
    \begin{tabular}{*{8}{>{\centering\arraybackslash}m{2cm}}}
        &\multicolumn{7}{c}{\textbf{Tree data}} \\ \cmidrule{2-8}
        &Uniform & DP & PY & +DP & +PY & HDP & HPY \\ \midrule
        Avg.\ num.\ &0.2627 & 0.2983 & 0.3077 & 0.2917 & 0.3202 & 0.3132 & 0.3315 \\ \bottomrule
    \end{tabular}}
    \caption{\label{tab: real} Tree species data: Average number of new species discovered per sampling step for each rmSSP model and the uniform baseline.}
\end{table}
\noindent Several noteworthy insights emerge from this experiment:
(a) All rmSSP models are clearly superior to the uniform baseline.
(b) The two PYP-based rmSSPs consistently outperform their DP-based counterparts, thanks to the extra flexibility provided by the discount parameter, which governs the rate at which new species appear.
(c) With the exception of the +DP, all models that borrow information across populations yield higher discovery rates than the independent specifications.
The +DP's weak performance stems from its underlying assumption that shared-species frequencies are proportional across populations, which seems inappropriate in this setting, and its lower flexibility compared to PYP, which prevents it from compensating for this misspecification.
To the best of our knowledge, this limitation of the +DP has not been previously noted in the literature.

\subsection{Synthetic data}

The \emph{tree dataset} exhibits high probabilities of ties across samples (see the right panel of Figure~\ref{fig: data} and recall that such probabilities are bounded above by the probability of a tie within a sample) and, thus, distributions in different samples are highly similar.
This makes it quite apparent that borrowing information across groups is advantageous.
Therefore, one could argue that the setting considered is overly favourable to rmSSPs relative to the independent models.
To assess whether, and under which conditions, borrowing information may become detrimental, we repeat the analysis on a simulated dataset.
\begin{figure}[htb]
  \centering
    \includegraphics[width=0.4\linewidth]{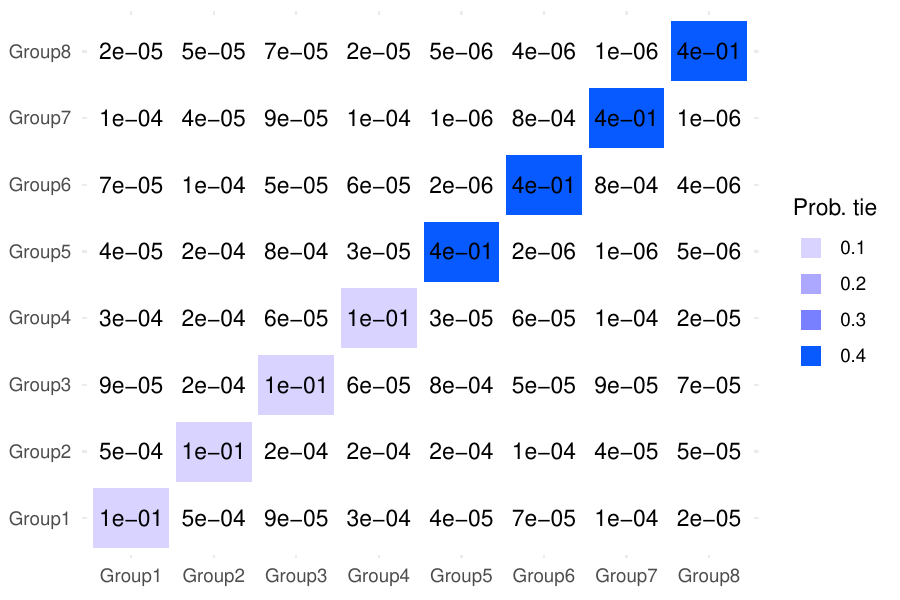}
  \caption{\label{fig: ptiesim} Probabilities of ties based on the \emph{true} distributions in the simulated scenario.}
\end{figure}
In the simulation experiment, we consider eight populations.
The true distribution of each arm is supported on a subset of $2,500$ species randomly drawn from a total of $3,000$, allowing for partial overlap of the supports across arms.
Each arm follows a Zipf distribution, where the probability assigned to the $k$th most frequent species in population $j$ is proportional to $k^{-s_{j}}$.
We set $s_{j} = 1.3$ for $j = 1,2,3,4$ and $s_{j} = 2$ for $j = 5,6,7,8$ \citep[cf.][]{battiston2018multi}.
However, before assigning the Zipf probability mass function, the $2,500$ selected species in each population are randomly permuted, leading to markedly different probability mass functions and low probabilities of ties across populations.
See Figure~\ref{fig: ptiesim}.
This scenario represents a worst-case setting for non-independent rmSSPs: although some species are shared across populations, borrowing information is undesirable.

\begin{figure}[htb]
 \begin{subfigure}{0.32\linewidth}
    \centering
    \includegraphics[ width=\linewidth]{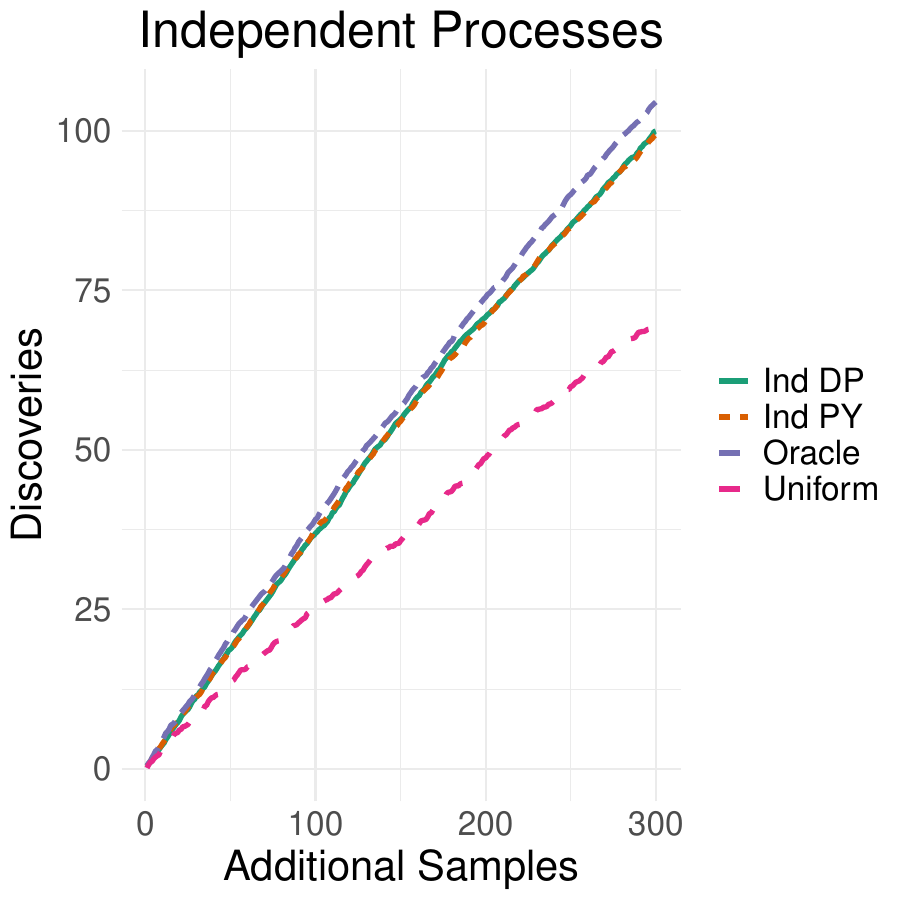}
 \end{subfigure}
 \hfill
 \begin{subfigure}{0.32\linewidth}
    \centering
    \includegraphics[ width=\linewidth]{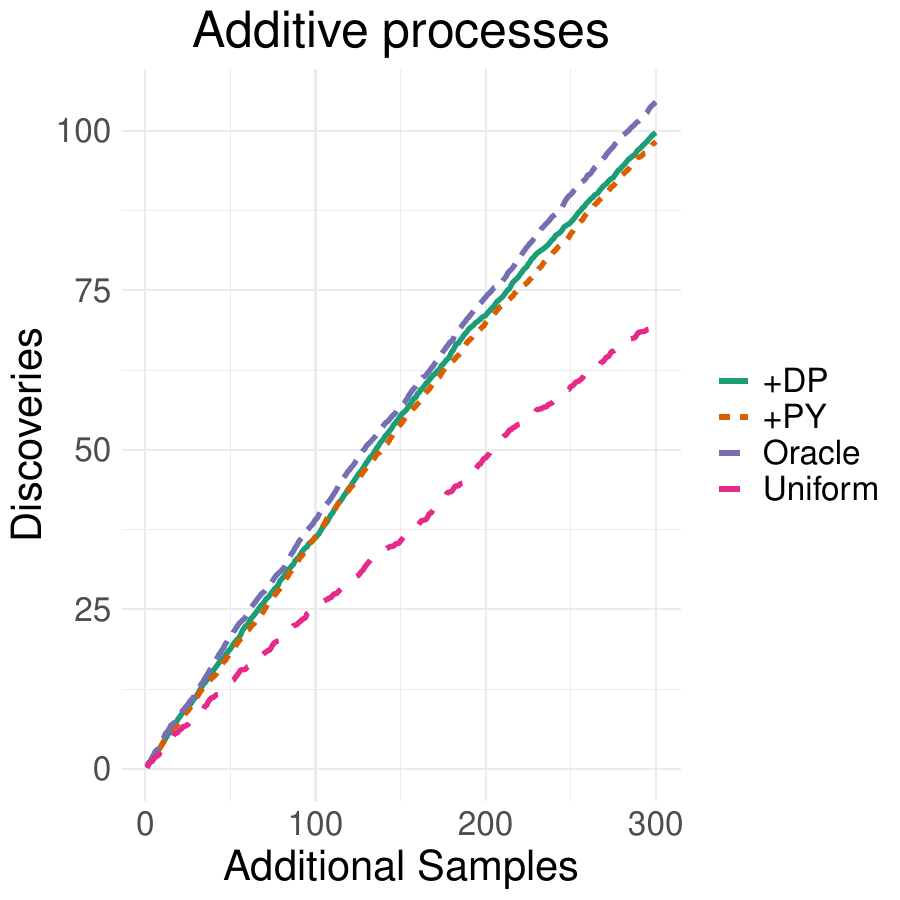}
 \end{subfigure}
 \hfill
 \begin{subfigure}{0.32\linewidth}
    \centering
    \includegraphics[ width=\linewidth]{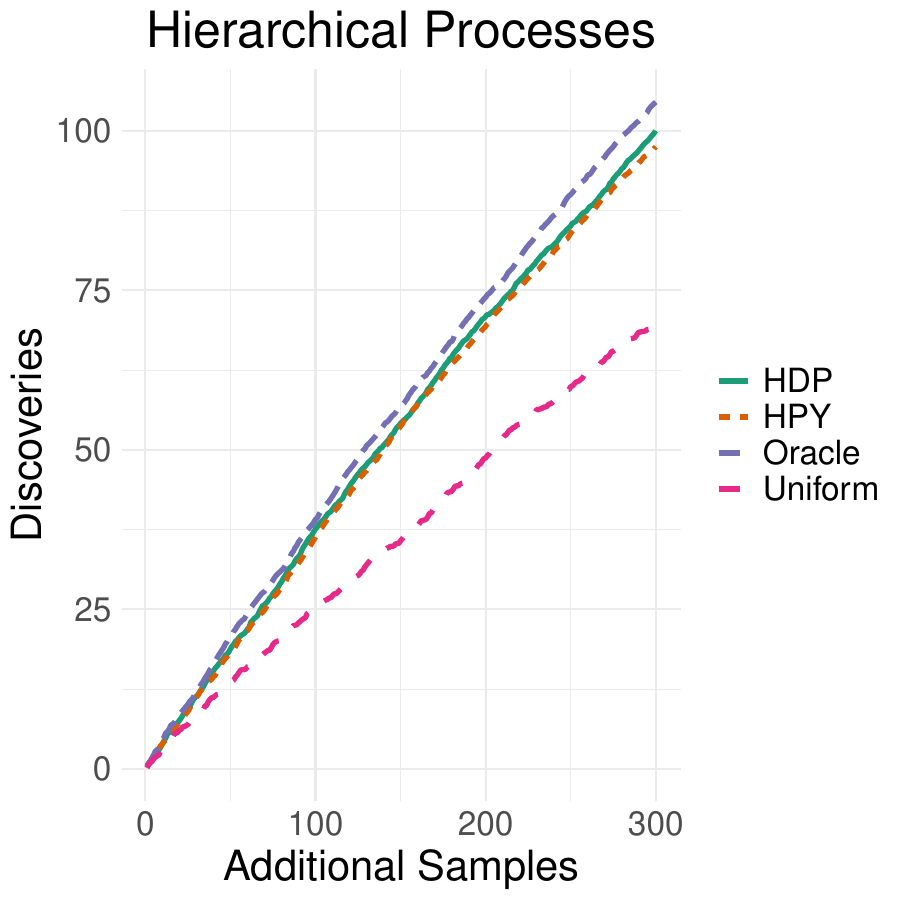}
 \end{subfigure}
\caption{\label{fig: sim low} {\small Simulated scenario with low probability of ties across populations:
Number of species discovered as a function of the additional sample size in the rmSSPs, the uniform model, and the oracle model.}
}
\end{figure}

\noindent Figure~\ref{fig: sim low} and Table~\ref{tab: sim} report averages over $20$ runs.
Comparisons are made against both the uniform model and the oracle model, which selects the arm with the highest \emph{true} frequency of unobserved species.
The results show that even in this scenario, the hierarchical and additive rmSSPs perform on par with the independent models, and close to the oracle in terms of species discovery.
This finding is reassuring, as it indicates that borrowing information, while unnecessary here, does not degrade performance.
\begin{table}[htb]
    \centering
    \renewcommand{\arraystretch}{1}
    \resizebox{0.9\textwidth}{!}{
    \begin{tabular}{*{9}{>{\centering\arraybackslash}m{2cm}}}
     &\multicolumn{8}{c}{\textbf{Simulated Scenario with low prob.\ of ties}} \\ \cmidrule{2-9}
     &Uniform & DP & PY & +DP & +PY & HDP & HPY & Oracle\\ \midrule
     Avg.\ num.\ &0.2340 & 0.3332 & 0.3312 & 0.3323 & 0.3277 & 0.3333 & 0.3250 & 0.3483 \\ \midrule
     RMSE & NA & 0.1563 & 0.0743 & 0.1621 & 0.0655 & 0.1929 & 0.0655 & 0 \\ \bottomrule
    \end{tabular}} 
    \caption{\label{tab: sim}
    Simulated scenario with low probability of ties across populations: Average number of species discovered per sampling step (Avg.\ num.) and root mean squared error (RMSE) of the estimated discovery probabilities in each population.}
\end{table}

\section{Conclusion}\label{sec: conclusion}
We introduced the class of mSSPs, a general framework extending Pitman's classical theory of species sampling models to the partially exchangeable setting.
Our contribution is twofold.\\
First, the mSSP framework provides a unifying perspective that encompasses most existing dependent nonparametric priors.
It fundamentally advances the understanding of their behaviour by revealing that borrowing of information across groups is fully determined by ties within and across groups.
These insights lead to principled strategies for prior specification, model calibration, and fair comparison across different subclasses of mSSPs.
A systematic empirical comparison of competing models will be pursued in future work.\\
Second, our approach is constructive.
It provides a modular recipe for building new models by combining EPPFs into structured dependence mechanisms.
This allows for the design of both new classes of mSSPs and novel models satisfying alternative probabilistic symmetries beyond partial exchangeability.
A first contribution along this path can be found in \textcite{fasano2025probabilistic}.\\
In addition to the systematic performance comparison of existing mSSPs and the development of new models, two further research directions emerge.
Our finding that borrowing of information is entirely governed by ties suggests that standard dependence measures may not be well-suited to random discrete structures.
This calls for a new theoretical framework based on the role of ties in generating, interpreting, and quantifying dependence.
A second, more probabilistic direction is to develop a standalone framework for pEPPFs, decoupled from partially exchangeable arrays, that incorporates any sequential generative construction.

\if0\blind 

\else
\section*{Acknowledgements} We are grateful to Jim Griffin, Fabrizio Leisen, Steven MacEachern, Li Ma, Peter M\"uller, Long Nguyen, Peter Orbanz, Riccardo Passeggeri, Judith Rousseau, and Aad van der Vaart for their valuable feedback and suggestions following presentations of this work at various seminars and conferences.
B.~Franzolini is supported by the National Recovery and Resilience Plan of Italy (PE1 FAIR - CUP B43C22000800006).
A.~Lijoi, I.~Pr\"unster and G.~Rebaudo are partially supported by the European Union - NextGenerationEU PRIN-PNRR (project P2022H5WZ9).
\fi

\printbibliography

\newpage

\begin{center}
{\LARGE{Supplementary Materials for\\} 
\bf Multivariate Species Sampling Models
}

\if0\blind
\else
{\large
\vspace*{0.5cm}
Beatrice Franzolini\\
DEMS Department, University of Milano-Bicocca
\vspace*{0.2cm}\\
Antonio Lijoi, Igor Pr\"unster\\
Bocconi Institute for Data Science and Analytics, Bocconi University\\
\vspace*{0.2cm}
Giovanni Rebaudo\\
ESOMAS Department, University of Torino and Collegio Carlo Alberto
}
\fi
\end{center}

\newrefsection

\setcounter{equation}{0}
\setcounter{definition}{0}
\setcounter{remark}{0}
\setcounter{page}{1}
\setcounter{table}{1}
\setcounter{figure}{0}
\setcounter{section}{0}
\numberwithin{table}{section}
\renewcommand{\theequation}{S.\arabic{equation}}
\renewcommand{\theremark}{S.\arabic{remark}}
\renewcommand{\thesubsection}{S.\arabic{section}.\arabic{subsection}}
\renewcommand{\thesection}{S.\arabic{section}}
\renewcommand{\thepage}{S.\arabic{page}}
\renewcommand{\thetable}{S.\arabic{table}}
\renewcommand{\thefigure}{S.\arabic{figure}}
\renewcommand{\thedefinition}{S.\arabic{definition}}

\vspace{0cm}

\section{Some basics on (univariate) species sampling}\label{sec: SSP supp}

In classical species sampling problems, a random sample $(X_{1},\ldots, X_{n})$ is extracted from an unknown and typically discrete distribution and each observed value corresponds to the species of a drawn individual.
Denoting with $P$ the unknown distribution of species in the population, we have
\[
X_{i}\mid P \simiid P \qquad \text{ for } i=1,\ldots,n.
\] 
To develop a Bayesian model for species sampling problems, a prior must be defined for the unknown distribution $P$.
In the univariate setting, the problem can be tackled relying on the large class of priors provided by species sampling processes (SSP), introduced by \textcite{pitman1996some} as a generalization of the Dirichlet process of \textcite{ferguson1973bayesian}.
\begin{definition}[SSP]\label{def: SSP}
	A random probability measure $P$ is a \emph{species sampling process} (SSP) if 
	\[
	P \eqas \sum_{h \geq 1} \pi_{h} \delta_{\theta_{h}} + \bigg(1-\sum_{h \ge 1} \pi_{h}\bigg) P_{0},
	\]
	where the atoms $(\theta_{h})_{h\geq1}$ are i.i.d.\ from the non-atomic probability measure $P_{0}$ and are independent of the random sub-probability vector of the weights $\bpi=(\pi_{h})_{h \ge 1}$.
	Moreover, if $\sum_{h \geq 1} \pi_{h}\eqas 1$, $P$ is said to be \emph{proper}.
\end{definition}
The corresponding model is defined once the observations are sampled independently from $P$ given $P$.
\begin{definition}[SSM]
	An infinite sequence of random variables $X_{1}, X_{2}, \ldots$ follows a \emph{species sampling model} (SSM) if it is exchangeable with an SSP directing measure.
	That is 
	\begin{align*}
			X_{i} \mid P \simiid P \quad \text{ for }i=1,2,\ldots, \quad 
			P &\sim \SSP(\L_{\bpi}, P_{0}).
	\end{align*}
\end{definition}
Any sample $(X_{1}, \ldots, X_{n})$ arising from a $P \sim \SSP(\L_{\bpi}, P_{0})$ induces a random partition of the labels of the observations in the sample, i.e., of $[n]=\{1, \ldots, n\}$.
More precisely, two observation labels $i$ and $l$ belong to the same block of the partition of $[n]$ (i.e., $X_{i}$ and $X_{l}$ are clustered together) if and only if $X_{i}=X_{l}$.
The discrete part of the SSP entails that two observations are clustered together with positive probability since, unless $\sum_{h \geq 1} \pi_{h} \eqas0$, $\bbP(X_{i}=X_{l})>0$.
The law of such a random partition (denoted $\Pi_{n}$) of $[n]$ is characterized by the exchangeable partition probability function (EPPF) \citep{pitman1996some}.

More precisely, let $\{C_{1}, \ldots, C_{K} \}$ be an arbitrary partition of $[n]$ for a given $n \in \bbN$ and $n_{k}=|C_{k}|$ for $k\in[K]$, then 
\begin{equation*}
	\bbP(\Pi_{n} = \{C_{1}, \ldots, C_{K}\})= \EPPF_{K}^{(n)}(n_{1},\ldots,n_{K}).
\end{equation*}
In words, $\EPPF_{K}^{(n)}(n_{1},\ldots,n_{K})$ can be interpreted as the probability of observing a particular (unordered) partition of $n$ observations into $K$ subsets of cardinalities $\{n_{1},\ldots,n_{K}\}$.
Note that the EPPF is defined on the space of the compositions of $n$, which can be interpreted as the space of the frequency vector of the partition in a given arbitrary order (e.g., the order of arrival).
Let $P=\sum_{h\ge1} \pi_{h} \delta_{\theta_{h}}$ be a proper SSP.
Then the induced EPPF can be computed as
\begin{align}
	\EPPF_{K}^{(n)}(n_{1},\ldots,n_{K}) = \bbE \bigg[ \sum_{h_{1} \ne \cdots \ne h_{K}} \prod_{k=1}^{K} \pi_{h_{k}}^{n_{k}} \bigg].
	\label{eq: EPPF from weights}
\end{align}
The EPPF characterizes the SSM \citep{pitman1996some}.
For any $n \in \bbN$, if $(X_{1}, \ldots, X_{n})$ arises from an SSM, its law can be obtained hierarchically as 
\begin{enumerate}
	\item sample the random partition $\Pi_{n}$ from the induced EPPF obtained as in \eqref{eq: EPPF from weights};
	\item sample i.i.d.\ the unique values associated with each set in the partition from $P_{0}$.
\end{enumerate}
The EPPF and the SSP can also be characterized by a specific sequence of predictive distributions \citep{pitman1996some} also known as the \emph{generalized Chinese restaurant process} (gCRP).
In the culinary metaphor, we can think of observations corresponding to customers in a restaurant, who arrive sequentially and sit at an already occupied table or a new table and each table serves a different dish (i.i.d.\ sampled from $P_{0}$).

It is theoretically straightforward to derive the predictive distribution associated with any SSP via ratios of EPPFs as an application of the definition of conditional probability, leading to
\begin{equation} \label{eq: CRP SSP supp}
	\bbP \big(X_{n+1}=x \mid \bX \big)=\begin{cases}
		\frac{\EPPF_{K}^{(n+1)}(n_{1},\ldots,n_{k}+1, \ldots, n_{K})}{\EPPF_{K}^{(n)}(n_{1},\ldots,n_{k}, \ldots, n_{K})} &\text{if } x=X_{k}^{\star} \quad \text{ and } k=1, \ldots, K\\
		\frac{\EPPF_{K+1}^{(n+1)}(n_{1},\ldots,n_{k}, \ldots, n_{K},1)}{\EPPF_{K}^{(n)}(n_{1},\ldots,n_{k}, \ldots, n_{K})} &\text{if } x=X_{K+1}^{\star},
	\end{cases}
\end{equation}
where $(X_{k}^{\star}: k=1,\ldots,K)$ denote the $K$ unique values of $\bX = (X_{1}, \ldots, X_{n})$ that were recorded with frequency $n_{1},\ldots,n_{K}$ and are i.i.d.\ sampled from $P_{0}$.
See \textcite{pitman1996some, pitman2006combinatorial, lee2013defining, ghosal2017fundamentals} for details and proofs about different characterizations of (univariate) SSM.

Although the analytical expression of the gCRP is available from the EPPF as shown in \eqref{eq: CRP SSP supp}, such an expression does not reduce to simple and tractable quantities in general.
However, a notable exception is the subclass of Gibbs-type priors \citep{gnedin2006exchangeable, deblasi2015gibbs}, which, thanks to the product partition form of the EPPF, allows the ratio of EPPFs in the gCRP to boil down to a simple ratio of constants for several notable examples, as in the well-known Chinese restaurant process (CRP) \citep{blackwell1973ferguson} induced by the Dirichlet process (DP) \citep{ferguson1973bayesian}.
The class of Gibbs-type priors is the most natural tractable generalization of the DP \citep{deblasi2015gibbs} and it includes the symmetric finite Dirichlet prior \citep{green2001modelling}, the Pitman--Yor process (PYP) \citep{pitman1997two}, the normalized inverse Gaussian (NIG) \citep{lijoi2005hierarchical}, the normalized generalized gamma process (NGGP) \citep{lijoi2007controlling}, mixtures of finite symmetric Dirichlet priors \citep{nobile1994bayesian, richardson1997bayesian, 
nobile2007bayesian, miller2018mixture} and the mixture of DP (MDP) models \citep{antoniak1974mixtures}.
In the following sections, we recall the analytical expression of the three different characterizations (i.e., SSP, EPPF, and gCRP) of some relevant and tractable examples of Gibbs-type priors commonly used in Bayesian analysis.

\subsection{Pitman--Yor process (PYP)}
We say that $P\sim \SSP(\L_{\bpi}, P_{0})$ follows a Pitman--Yor process, i.e., $P \sim \PYP(\sigma, \alpha, P_{0})$, with $P_{0}$ a non-atomic measure if it is a proper SSP with $\L_{\bpi} \sim \GEM(\sigma,\alpha)$, where the two-parameter GEM distribution, named after Griffiths, Engen, and McCloskey, can be thought as arising from the stick-breaking construction where the $\pi_{i}$'s are such that $\pi_{i}=v_{i} \prod_{\ell=1}^{i-1} (1-v_{\ell} )$, with $v_{i}\sim \Beta (1-\sigma, \alpha + i\sigma)$, $i\ge 1$, $\sigma\in[0,1)$ and $\alpha>-\sigma$.

The following EPPF characterizes the PYP
\begin{equation}\label{eq: EPPF PYP supp}
	\EPPF_{K}^{(n)}(n_{1}, \ldots, n_{K}) = \dfrac{\prod_{k=1}^{K-1} (\alpha +k\, \sigma)}{(\alpha+1)_{n-1}} \prod_{k=1}^{K} (1-\sigma)_{n_{k}-1},
\end{equation}
where $(x)_{n} = x(x+1) \cdots (x+n-1)$ is the $n$th ascending factorial.

Denoting with $X_{1}, X_{2}, \ldots$ an SSM from $P \sim \PYP(\sigma, \alpha, P_{0})$, we can derive the well-known gCRP of the PYP from the EPPF in \eqref{eq: EPPF PYP supp} applying the definition of conditional probability.
\begin{equation*}
	\bbP \big(X_{n+1}=x \mid \bX \big)=\begin{cases}
		\frac{n_{k}- \sigma }{\alpha + n} &\text{if } x=X_{k}^{\star} \quad \text{ and } k=1, \ldots, K\\
		\frac{\alpha + \sigma K}{\alpha + n} &\text{if } x=X_{K+1}^{\star}.
	\end{cases}
\end{equation*}

\subsection{Dirichlet process (DP)}
If we consider $P \sim \PYP(\sigma, \alpha, P_{0})$ as in the previous section and we restrict $\sigma=0$ and $\alpha>0$ we obtain the relevant special case of the Dirichlet process, i.e., $P \sim \DP(\alpha, P_{0})$.
Thus we can specialize the distribution of the weights to $\GEM(\alpha)$, the induced EPPF in \eqref{eq: EPPF PYP supp} simplifies to
\begin{equation*}
	\EPPF_{K}^{(n)}(n_{1}, \ldots, n_{K} ) = \dfrac{\alpha^{K} \Gamma(\alpha)}{\Gamma(\alpha+n)} \prod_{k=1}^{K} (n_{k}-1)!,
\end{equation*}
and the corresponding CRP
\begin{equation*}
	\bbP \big(X_{n+1}=x \mid \bX \big)=\begin{cases}
		\frac{n_{k} }{\alpha + n} &\text{if } x=X_{k}^{\star} \quad \text{ and } k=1, \ldots, K\\
		\frac{\alpha }{\alpha + n} &\text{if } x=X_{K+1}^{\star}.
	\end{cases}
\end{equation*}

\subsection{Dirichlet multinomial (DM)}
Here we consider an SSP $P$ with a fixed known number $M$ of species in the population (with $M \in \bbN$) that follow a finite-dimensional (symmetric) Dirichlet multinomial (DM).
That is, for a fixed $M \in \bbN$,
\begin{equation*}
	P = \sum_{h=1}^{M} \pi_{h} \delta_{\theta_{h}},
\end{equation*}
where $(\pi_{1}, \ldots, \pi_{M}) \sim \Dir(\tau, \ldots, \tau)$ $\perp$ $\theta_{h} \simiid P_{0}$.
We write $P \sim \DM_{M}(\tau,P_{0})$.

Then we can derive the induced EPPF as 
\begin{equation}\label{eq: EPPF Dir supp}
	\EPPF_{K}^{(n)}(n_{1},\ldots,n_{K}) = \frac{M!}{(M-K)!} \frac{\Gamma( \tau \, M) }{\Gamma(n + \tau \, M)\Gamma(\tau )^{K}} \prod_{k=1}^{K} \Gamma(n_{k} + \tau )
\end{equation}
and the corresponding gCRP
\begin{equation*}
	\bbP \big(X_{n+1}=x \mid \bX \big) = \begin{cases}
		\frac{n_{k} +\tau}{\tau M + n} &\text{if } x=X_{k}^{\star} \quad \text{ and } k=1, \ldots, K\\
		\frac{\tau (M-K)}{\tau M + n} &\text{if } x=X_{K+1}^{\star} \text{ and } K < M.
	\end{cases}
\end{equation*}

\subsection{Gnedin process (GN)}
Allowing for an unknown $M$ in a finite-dimensional symmetric Dirichlet multinomial process, the model becomes a mixture of symmetric Dirichlet models.
A relevant example is the {\em Gnedin process} (with discount parameter equal to $-1$ and $\zeta=0$ in the original notation of \cite{gnedin2010species}).
The corresponding EPPF is
\begin{equation*}
	\EPPF_{K}^{(n)}(n_{1},\ldots,n_{K}) = \sum_{m=K}^{\infty} \EPPF_{K}^{(n)}(n_{1},\ldots,n_{K}\mid M=m)\, \bbP(M=m),
\end{equation*}
where $\EPPF_{K}^{(n)}(n_{1},\ldots,n_{K}\mid M=m)$ is the EPPF of the $M$-symmetric Dirichlet prior in \eqref{eq: EPPF Dir supp}, with $\tau=1$ and $ \bbP(M=m) = \frac{\gamma (1-\gamma)_{m-1}}{m!}$, $\gamma \in (0,1)$.

The corresponding gCRP boils down to the following simple tractable expression
\begin{equation*}
	\bbP \big(X_{n+1}=x \mid \bX \big) = \begin{cases}
		\frac{(n_{k}+1)(n-K+\gamma)}{n(n+\gamma)} &\text{if } x=X_{k}^{\star} \quad \text{ and } k=1, \ldots, K\\
		\frac{K^{2}- K \gamma}{n(n+\gamma)} &\text{if } x=X_{K+1}^{\star}.
	\end{cases}
\end{equation*}
We denote the corresponding SSP $P \sim \GN(\gamma, P_{0})$.

\section{Proofs} 
\subsection{Proof of Proposition~\ref{prop: marginal mSSP}}
The proof immediately follows from the Definition~\ref{def: mSSP1} of mSSP.
More precisely, restricting to $j\in \J$ yields the same expression with the same atom sequence $(\theta_{h})_{h \ge 1}$ and the sub-array $(\bpi_{j})_{j\in \J}$, which remains independent of $(\theta_{h})_{h \ge 1}$.
Thus, $(P_{j})_{j \in \J}$ is an mSSP with parameters $(\L_{\bpi_{\J}}, P_{0})$.
Finally, if $\sum_{h\ge1}\pi_{j,h}=1$ a.s.\ for all $j\in[J]$, the same immediately holds for all $j\in \J$, entailing that also the marginal distributions are proper mSSP.

\subsection{Proof of Proposition~\ref{prop: characterization bivariate mssp}}
\begin{proof}
	To prove the statement, we want to show the non-trivial implication of the if and only if, i.e., if $(P_{1}, P_{2})$ is an mSSP with non-atomic base measure $P_{0}$, there exist $(\tP_{1}, \tP_{2})$ defined as in \eqref{eq:3components}, that is
	\begin{equation*}
		\tP_{j} \eqas \sum_{h \ge 1} \pi^{(1,2)}_{j,h} \delta_{\theta_{0,h}} + \sum_{h^{\prime} \geq 1 } \pi^{(j)}_{j,h^{\prime}} \delta_{\theta_{j,h^{\prime}}} + \pi^{(j)}_{j,0} P_{0}, \quad \quad \text{for } j=1,2.
	\end{equation*}
	where, by definition of random probability vector, $\sum_{h\geq1}\pi^{(1,2)}_{j,h} + \sum_{h^{\prime}\geq 0 } \pi^{(j)}_{j,h^{\prime}} \eqas 1 $, $\pi_{j,h}^{(1,2)}, \pi_{j,h^{\prime}}^{(j)} \overset{\mathrm{a.s.}}{\ge} 0$ for $j=1,2$, $h \ge 1$ and $h^{\prime} \ge 0$, the atoms are jointly independent of the weights and such that $\theta_{j,h} \simiid P_{0}$, for $j=0,1,2, \,\ h \ge 1$ and \mbox{$\bbP[\pi^{(1,2)}_{1,h}>0, \pi^{(1,2)}_{2,h}>0]>0$} for all $h$ indexing the first sum.
	
	From the Definition~\ref{def: mSSP1} of mSSP, we write for $j=1,2$
	\[
	P_{j} \eqas \sum_{h \ge 1} \pi_{j,h} \delta_{\theta_{h}} + \bigg(1-\sum_{h \ge 1} \pi_{j,h}\bigg) P_{0} = \sum_{h \in \H} \pi_{j,h} \delta_{\theta_{h}} + \bigg(1-\sum_{h \in \H} \pi_{j,h}\bigg) P_{0},
	\]
	where we denote by $\H \coloneqq \{h \ge 1\}$ the index set of the first sum.
	Note that the cardinality of $\H$ can be in $\{0\} \cup \bbN \cup \{\infty\}$, and recall that we use the convention that, for any $(x_{h})_{h}$, $\sum_{h=1}^{0} x_{h} = \sum_{h \in \emptyset} x_{h} = 0$.
	We define $\pi_{j,0}^{(j)} \coloneqq 1-\sum_{h \in \H} \pi_{j,h}$ and we partition $\H$ into $\{\H_{0}, \bar{\H}_{0}\}$,
	where 
	\[
	\H_{0} \coloneqq \{h\in \H :\:\: \bbP(\pi_{1,h}>0, \pi_{2,h}>0)>0\} = \{h \in \H :\:\: \bbP(\pi_{1,h}\pi_{2,h}>0)>0\}
	\]
	is the set of shared atoms and $\bar{\H}_{0} = \H \setminus \H_{0}$.

     Let us define, $\pi_{j,h}^{(1,2)} \coloneqq \pi_{j,h}$ for $h \in \H_{0}$ and $j=1,2$, and $\pi_{j,h}^{(j)} \coloneqq \pi_{j,h}$ for $h \in \bar{\H}_{0}$ and $j=1,2$.
     
    Note that the subset of atoms $\theta_{h}$ that are indexed by $\H_{0}$, i.e, $(\theta_{h})_{h\in \H_{0}}$, are actually shared across the two groups with positive probability (i.e., they are atoms associated to weights that are jointly positive in both $P_{1}$ and $P_{2}$ with positive probability).
    Conversely, the subset of atoms $\theta_{h}$ that are indexed by $\bar{\H}_{0}$, i.e, $(\theta_{h})_{h\in \bar{\H}_{0}}$, while not necessary separable into two independent sequences of atoms each corresponding to one of the two processes, are almost surely not shared across $j$ by definition of $\bar{\H}_{0}$.
    Indeed, for every $h \in \bar{\H}_{0}$ almost surely at most one between $\pi_{1,h}^{(1)}=\pi_{1,h}$ and $\pi_{2,h}^{(2)}=\pi_{2,h}$ can be different from $0$.
    
    Now, let $\theta_{0,h}:=\theta_{h}$, for $h \in \H_{0} $, and, thus, $\theta_{0,h}\simiid P_{0}$, for $h \in \H_{0}$ and define two mutually independent sequences $(\theta_{1,h})_{h \in \bar{\H}_{0}}$ and $(\theta_{2,h})_{h \in \bar{\H}_{0}}$, i.i.d.\ from $P_{0}$ and mutually independent of the weights $(\pi_{j,h}^{(1,2)})_{j,h}$ and $(\pi_{j,h}^{(j)})_{j,h}$ and of the shared atoms $(\theta_{0,h})_{h \in \H_{0}}$.
    Note that such sequences of atoms $(\theta_{j,h})_{j,h}$, for $j=1,2$ are iid from $P_{0}$ across $j=1,2$ and $h\in \bar{\H}_{0}$ by construction and, thus, are not shared almost surely across $j$ by the non-atomicity of $P_{0}$.
    Thus, $(\tP_{1}, \tP_{2})$ will have the same law of $(P_{1},P_{2}) \sim \mSSP(\L_{\bpi},P_{0})$, even though the non-shared atoms are not a.s.\ equal across the two representations (since they contribute to at most one component almost surely).

    Indeed, note that $\theta_{j,h} \simiid P_{0}$, for $j=0, 1,2$ and each $h$, and are independent of the entire weight array and, for $j=1,2$,
	\[
	\sum_{h\in \H}\pi_{j,h} \delta_{\theta_{h}} = \sum_{h \in \H_{0}}\pi_{j,h} \delta_{\theta_{h}} + \sum_{h \in \bar{\H}_{0}}\pi_{j,h} \delta_{\theta_{h}},
	\]
	\[
	\sum_{h \in \H_{0}}\pi_{j,h} \delta_{\theta_{h}} \eqd \sum_{h \in \H_{0}} \pi_{j,h}^{(1,2)}\,\delta_{\theta_{0,h}} \text{ and } \sum_{h \in \bar{\H}_{0}} \pi_{j,h} \delta_{\theta_{h}} \eqd \sum_{h \in \bar{\H}_{0}} \pi_{j,h}^{(j)}\,\delta_{\theta_{j,h}}.
	\]
    and 
    \[
        \bigg(\sum_{h \in \H_{0}}\pi_{j,h} \delta_{\theta_{h}}, \sum_{h \in \bar{\H}_{0}} \pi_{j,h} \delta_{\theta_{h}} \bigg) \eqd \bigg(\sum_{h \in \H_{0}} \pi_{j,h}^{(1,2)}\,\delta_{\theta_{0,h}}, \sum_{h \in \bar{\H}_{0}} \pi_{j,h}^{(j)}\,\delta_{\theta_{j,h}} \bigg).
\]
    Thus, $\sum_{h \in \H_{0}} \pi^{(1,2)}_{j,h} + \sum_{h \in \bar{\H}_{0}} \pi^{(j)}_{j,h} + \pi_{j,0}^{(j)} \eqas 1$.
	To conclude the proof, we just relabel the indices in both $\H_{0}$ and $\bar{\H}_{0}$ (not affecting the law of the sums) such that they are ordered integers starting from $1$ with no gaps and remap the elements in the corresponding sums accordingly.
\end{proof}

\subsection{Proof of Proposition~\ref{prop: mean var PA}}
\begin{proof}
	By the law of iterated expectations, the first and second moments of $P_{j}(A)$ are equal to
	\begin{align*}
		\bbE[P_{j}(A)] &= \bbP(X_{j,i} \in A)=P_{0}(A) \\ 
        \bbE[P_{j}(A)^{2}] &= \bbP(X_{j,i} \in A, X_{j,l} \in A), \text{ with } i \neq l.
	\end{align*}
	Disintegrating with respect to $\{X_{j, i}= X_{j,l}\}$ and its complement to exploit that the unique values are i.i.d.\ sampled from the non-atomic probability measure $P_{0}$ leads to 
	\begin{align*}
        &\bbP(X_{j,i} \in A, X_{j,l} \in A) \\
        &\quad=\,\bbP(X_{j,i} = X_{j,l}) \bbP(X_{j,i} \in A, X_{j,l} \in A \mid X_{j,i} = X_{j,l})\\
        &\qquad+ \bbP(X_{j,i} \ne X_{j,l}) \bbP(X_{j,i} \in A, X_{j,l} \in A \mid X_{j,i} \ne X_{j,l})\\
        &\quad=\, \bbP(X_{j,i} = X_{j,l}) P_{0}(A) + \bbP(X_{j,i} \ne X_{j,l}) P_{0}(A)^{2}.
	\end{align*}
	Finally,
	$\Var[P_{j}(A)] = \bbE[P_{j}(A)^{2}] - \bbE[P_{j}(A)]^{2} = \bbP(X_{j,i}= X_{j,l}) P_{0}(A) [1-P_{0}(A)]$.
\end{proof}

\subsection{Proof of Proposition~\ref{prop: CorrProb}}
\begin{proof}
	For $j\neq k$, by the law of iterated expectations, we get
	\begin{align*}
		\bbE[P_{j}(A) P_{k}(A)] = \bbP(X_{j,i} \in A, X_{k,m} \in A).
	\end{align*}
	 Disintegrating with respect to $\{X_{j, i}= X_{j,l}\}$ and its complement to exploit that the unique values are i.i.d.\ sampled from the non-atomic probability measure $P_{0}$ leads to 
	\begin{align*}
		&\bbP(X_{j,i} \in A, X_{k,m} \in A) &\\
        & \quad =\bbP(X_{j,i} = X_{k,m}) \bbP(X_{j,i}\in A, X_{k,m} \in A \mid X_{j,i} = X_{k,m})\\
		& \qquad + \bbP(X_{j,i} \ne X_{k,m}) \bbP(X_{j,i} \in A, X_{k,m} \in A \mid X_{j,i} \ne X_{k,m})\\
		& \quad = \bbP(X_{j,i} = X_{k,m}) P_{0}(A) + \bbP(X_{j,i} \ne X_{k,m}) P_{0}(A)^{2}.
	\end{align*}
	Thus, $\Cov[P_{j}(A),P_{k}(A)]
	=\bbP(X_{j,i} = X_{k,m}) P_{0}(A) [1-P_{0}(A)]$.
	The correlation is obtained using Proposition~\ref{prop: mean var PA}.
\end{proof}

\subsection{Proof of Corollary~\ref{cor: CorrProb}}
The proof follows trivially from Proposition~\ref{prop: CorrProb}.

\subsection{Proof of Theorem~\ref{thm: Corr1}}

\begin{proof}
First note that the correlation in the statement is well defined since by assumption $A$ is a measurable set such that $0< P_{0}(A) <1$ and if at least one of (p-i) or (p-ii) holds, then 
$\bbP(X_{j,1}=X_{j,2})>0$ and $\bbP(X_{k,1}=X_{k,2})>0$.
Thus, 
\[
    0<\Var[P_{\ell}(A)] = \bbP(X_{\ell,1}=X_{\ell,2}) P_{0}(A)\big[1-P_{0}(A)\big] < \infty, \text{ for }\ell=j,k.
\]
Note that $P_{j} \eqas P_{k}$ if and only if (iff) $\bX$ is exchangeable by the de Finetti theorem.
Note also that if $P_{j} \eqas P_{k}$, then $P_{j}(A) \eqas P_{k}(A)$ and $\Cor[P_{j}(A), P_{k}(A)]=1$ trivially.
Therefore, below we prove the other implication: under (p-i) or (p-ii) if $\Cor[P_{j}(A), P_{k}(A)]=1$, then $P_{j} \eqas P_{k}$.
\\
\\
\textbf{Case (p-i)}: Assume $(P_{j},P_{k})$ is a proper mSSP.

By the definition of proper mSSP, we know that 
\[
	P_{j} \eqas \sum_{h \ge 1} \pi_{j,h} \delta_{\theta_{h}} \quad \text{ and } \quad
	P_{k} \eqas \sum_{h \ge 1} \pi_{k,h} \delta_{\theta_{h}}.
\]
where $(\pi_{\ell,h})_{h\geq 1}$ is a probability vector, for $\ell = j,k$ and $\theta_{h}\simiid P_{0}$.
    
Moreover, by Cauchy--Schwarz inequality, we have a.s.\ that
\begin{align*}
	\sqrt{\sum_{h \ge 1} \pi_{j,h}^{2}}\sqrt{\sum_{h \ge 1} \pi_{k,h}^{2}} \geq \sum_{h \ge 1} \pi_{j,h} \pi_{k,h}.
\end{align*}
Assume by contradiction that the event $\{\pi_{j,h} \neq \pi_{k,h}$ for at least one $h\}$ has positive probability.
This implies that with positive probability, the above inequality is strict and thus, with positive probability, we have
\[
	\bbP(X_{j,1}=X_{k,1} \mid P_{j}, P_{k}) < \sqrt{\bbP(X_{j,1}=X_{j,2} \mid P_{j}, P_{k}) } \sqrt{\bbP(X_{k,1}=X_{k,2} \mid P_{j}, P_{k}) }
\]
which implies (by monotonicity of the expectation and Cauchy--Schwarz) 

\begin{align*}
    \bbP(X_{j,1}=X_{k,1}) &< 
    \bbE\bigg[\sqrt{\bbP(X_{j,1}=X_{j,2} \mid P_{j}, P_{k}) } \sqrt{\bbP(X_{k,1}=X_{k,2} \mid P_{j}, P_{k}) }\bigg]\\
    & \le \sqrt{\bbE\bigg[\bbP(X_{j,1}=X_{j,2} \mid P_{j}, P_{k}) \bigg]} \sqrt{\bbE\bigg[\bbP(X_{k,1}=X_{k,2} \mid P_{j}, P_{k}) \bigg]}\\
    &=
    \sqrt{\bbP(X_{j,1}=X_{j,2} )\bbP(X_{k,1}=X_{k,2})}.
\end{align*}
Thus,
\[
    \Cor[P_{j}(A),P_{k}(A)] < 1,
\]
completing the contradiction.
Therefore, we have $\pi_{j,h} = \pi_{k,h}$ a.s., for all $h$, and thus $P_{j}\eqas P_{k}$.
\\
\\
\textbf{Case (p-ii)}: Assume $\bbP(X_{j,1}=X_{j,2}) = \bbP(X_{k,1}=X_{k,2})=:\rho$.

 By Proposition~\ref{prop: CorrProb}, we have
\[
    \Cor[P_{j}(A),P_{k}(A)] = \frac{\bbP(X_{j,1}=X_{k,1})}{\rho} \qquad \forall j\neq k\in[J]
\]
Thus, $\Cor[P_{j}(A),P_{k}(A)] = 1$ entails $\bbP(X_{j,1}=X_{k,1})=\rho$.

Recall that,
\[
    \bbP(X_{\ell,1}=X_{\ell,2})=\bbE \bigg[\sum_{h\ge 1} \pi_{\ell,h}^{2} \bigg] \text{ for }\ell =j,k \text{ and } \bbP(X_{j,1}=X_{k,1})=\bbE \bigg[\sum_{h\ge 1} \pi_{j,h}\pi_{k,h} \bigg].
\]
Therefore,
\begin{align*}
    \bbE\bigg[\sum_{h\ge 1} \big(\pi_{j,h} -\pi_{k,h} \big)^{2} \bigg] = \bbE\bigg[\sum_{h\ge 1} \pi_{j,h}^{2} \bigg] +\bbE\bigg[\sum_{h\ge 1}\pi_{k,h}^{2} \bigg] - 2 \bbE\bigg[\sum_{h\ge 1} \pi_{j,h} \pi_{k,h} \bigg] = \rho + \rho - 2\rho= 0.
\end{align*}
Thus, by the non-negativity of the square,
\[
    \pi_{j,h} \eqas \pi_{k,h} \quad \forall h.
\]
This entails that also $\big(1-\sum_{h\ge 1} \pi_{j,h} \big) \eqas \big(1-\sum_{h\ge 1} \pi_{k,h} \big)$.
Recalling that the atoms $\theta_{h}$ are shared across $k$ and $j$, we can conclude that $P_{k} \eqas P_{j}$.
\end{proof}

\subsection{Proof of Theorem~\ref{thm: CorInd}}
\begin{proof}
Clearly $P_{j} \perp P_{k}$ entails $\Cor[P_{j} (A), P_{k} (A)]= 0$.
We want to show that $\Cor[P_{j}(A), P_{k} (A)]= 0$ entails $P_{j} \perp P_{k}$.

Let us consider the representation of $(P_{j}, P_{k}) \sim \mSSP(\L_{\bpi}, P_{0})$ provided in Proposition~\ref{prop: characterization bivariate mssp} and let us rewrite explicitly each measure as a mixture of two components, i.e.,
\[
    P_{j} \eqas \omega^{(j,k)}_{j}\sum_{h \ge 1} \bar\pi^{(j,k)}_{j,h} \delta_{\theta_{h}} + \left[1-\omega^{(j,k)}_{j}\right]\bigg[\bar\pi^{(j)}_{j,0}P_{0} +\sum_{h^{\prime}\geq 1}\bar\pi^{(j)}_{j,h^{\prime}} \delta_{\theta_{j,h^{\prime}}} \bigg]
\]
and 
\[
    P_{k} \eqas \omega^{(j,k)}_{k} \sum_{h \ge 1} \bar\pi^{(j,k)}_{k,h} \delta_{\theta_{h}} + \left[1-\omega^{(j,k)}_{k}\right]\bigg[\bar\pi^{(k)}_{k,0}P_{0} +\sum_{h^{\prime}\geq 1}\bar\pi^{(k)}_{k,h^{\prime}} \delta_{\theta_{k,h^{\prime}}} \bigg],
\]
where for $\ell =j, k$
\[
    \bar\pi^{(j,k)}_{\ell,h} = \frac{\pi^{(j,k)}_{\ell,h}}{\sum_{s \geq 1 } \pi^{(j,k)}_{\ell,s}}, \qquad \bar\pi^{(\ell)}_{\ell,h^{\prime}} = \frac{\pi^{(\ell)}_{\ell,h^{\prime}}}{\sum_{s \geq 0 } \pi^{(\ell)}_{\ell,s}}, \quad \text{and} \quad \omega^{(j,k)}_{\ell} = \sum\limits_{h \geq 1} \pi^{(j,k)}_{\ell,h}.
\]
Recall the conventions $\sum_{h \in \emptyset} x_{h} =0$ and $0/0=0$ and that by Proposition~\ref{prop: CorrProb} and Corollary~\ref{cor: CorrProb}, 
\[
    \Cor[P_{j}(A), P_{k}(A) ]= 0 \text{ iff } \bbP(X_{j,1}=X_{k,1})=0.
\]
Note that, by Proposition~\ref{prop: characterization bivariate mssp}, we can show that $\omega^{(j,k)}_{j} \eqas \omega^{(j,k)}_{k}\eqas 0$.
Indeed, if the shared index set, denoted by $\H_{0}$ in the proof of Proposition~\ref{prop: characterization bivariate mssp}, is empty, then $\omega^{(j,k)}_{j}\eqas\omega^{(j,k)}_{k}\eqas 0$ by construction.
Instead, if we assume by contradiction that (w.l.o.g.) $\bbP\left[\omega^{(j,k)}_{j}>0\right]>0$ then by Proposition~\ref{prop: characterization bivariate mssp}, we have that
\[
    \bbP(X_{j,1}=X_{k,1}) = \sum_{h \ge 1} \bbE\left[ \omega^{(j,k)}_{j} \bar\pi^{(j,k)}_{j,h}\, \omega^{(j,k)}_{k} \bar\pi^{(j,k)}_{k,h}\right] \ge \bbE\left[ \omega^{(j,k)}_{j} \bar\pi^{(j,k)}_{j,h}\, \omega^{(j,k)}_{k} \bar\pi^{(j,k)}_{k,h}\right],
\]
for any fixed $h$ in the shared index, e.g., $h=1$, where $\bbP\left[\omega^{(j,k)}_{j}>0\right]>0$  implies 
\[
    \bbP\left[\omega^{(j,k)}_{j} \bar\pi^{(j,k)}_{j,1} \omega^{(j,k)}_{k} \bar\pi^{(j,k)}_{k,1} >0 \right]>0
\]
that entails
\[
    \bbP(X_{j,1}=X_{k,1}) \ge \bbE\left[ \omega^{(j,k)}_{j} \bar\pi^{(j,k)}_{j,1}\, \omega^{(j,k)}_{k} \bar\pi^{(j,k)}_{k,1}\right] >0
\]
that contradicts $\bbP(X_{j,1}=X_{k,1})=0$.
    
Since $\omega^{(j,k)}_{j}\eqas\omega^{(j,k)}_{k}\eqas 0$, we can rewrite 
\begin{equation}\label{eq:temp}
    P_{j} \eqas \bar\pi^{(j)}_{j,0}P_{0} +\sum_{h^{\prime}\geq 1}\bar\pi^{(j)}_{j,h^{\prime}} \delta_{\theta_{j,h^{\prime}}}
\end{equation}
and 
\begin{equation}\label{eq:temp2}
	P_{k} \eqas \bar\pi^{(k)}_{k,0}P_{0} +\sum_{h^{\prime}\geq 1}\bar\pi^{(k)}_{k,h^{\prime}} \delta_{\theta_{k,h^{\prime}}}.
\end{equation}
Therefore $P_{j} \perp P_{k}$ since by Proposition~\ref{prop: characterization bivariate mssp} and Definition~\ref{def: rmSSP} of rmSSP both the weights and the atoms in \eqref{eq:temp} and \eqref{eq:temp2} are independent across measures.
\end{proof}

\subsection{Proof of Proposition~\ref{prop: CorrObs}}
\begin{proof}
Define the random variable $Z$, so that $Z=1$, if $X_{j,i}=X_{k,m}$, and $Z=0$, otherwise.
By the law of total covariance
\begin{align*}
	\Cov(X_{j,i},X_{k,m})&= \bbE\left[\Cov(X_{j,i},X_{k,m}\mid Z)\right] + \Cov\left(\bbE\left[X_{j,i}\mid Z\right],\bbE\left[X_{k,m}\mid Z\right]\right)\\
	&=\bbE\left[\Cov(X_{j,i},X_{k,m}\mid Z)\right] + 0\\
		&=\Cov(X_{j,i},X_{k,m}\mid Z=1) \bbP(X_{j,i}=X_{k,m})\\
		&=\bbP(X_{j,i}=X_{k,m}) \Var(X^{\star}),
\end{align*}
where $X^{\star}\sim P_{0}$ and $\Cov(X_{j,i}, X_{k,m}\mid Z=1) = \Var(X^{\star})$ is obtained since the conditioning to $Z = 1$ implies that both observations are equal to the same atom, which is itself sampled from $P_{0}$.
Note also that we set
$\Cov(\bbE[X_{j,i}\mid Z],\bbE[X_{k,m}\mid Z])=0$ since $\bbE(X_{j,i}\mid Z=1)=\bbE(X_{j,i}\mid Z=0)=\bbE(X^{\star})$ and exploit that $\Cov(X_{j,i},X_{k,m}\mid Z=0)=0$ since given $Z=0$ $X_{j,i}$ and $X_{k,m}$ are i.i.d.\ sampled from $P_{0}$ (either they are equal to distinct atoms $\theta_{h}$s or arise from the non-atomic part of the mSSP).
    
The statement follows by dividing the expression of the covariance by $\sqrt{\Var(X_{j,i})\Var(X_{k,m})} = \Var(X_{j,i}) = \Var(X^{\star})$.

Finally, note that the facts that $P_{0}$ has finite second moment and it is a non-atomic probability measures entail that $P_{0}$ has finite and positive variance and thus the correlation is well defined since marginally the observations are distributed according to $P_{0}$.
\end{proof}

\subsection{Proof of Corollary~\ref{cor: CorrX}}
The proof follows trivially from Proposition~\ref{prop: CorrObs}.

\begin{remark}
Note that $X_{j,i} \perp X_{k,m}$ is a marginal independence statement between two arbitrary random variables in groups $j$ and $k$, that does not imply independence of larger samples from the two groups without further assumptions such as that $(P_{j},P_{k})$ is a rmSSP.
\end{remark}

\subsection{Proof of Proposition~\ref{prop: MomSSP}}
\begin{proof} Define $\bX_{j,1:q} = (X_{j,1}, \ldots, X_{j,q})$,
\begin{align*}
	\bbE[P_{j}(A)^{q}] = \bbP(\bX_{j,1:q} \in A^{q}).
\end{align*}
Let $\Pi_{q}^{(j)}$ be the random partition of the labels $[q]:=\{1, \ldots, q\}$ of the data units in $\bX_{j,1:q}$, which is induced by the ties among $\bX_{j,1:q}$ and takes values in the set $\P(q)$, that is the set of all possible partitions of $[q]$.
Disintegrate with respect to the random partition $\Pi_{q}^{(j)}$ to recover independence and aggregate by symmetry induced by exchangeability.
Let also $K^{(j)}_{1:q}$ denote the number of blocks in $\Pi_{q}^{(j)}$.
\begin{align*}
    \bbP\left(\bX_{j,1:q} \in A^{q} \right) &= \sum_{ \pi \in \P(q) } \bbP \big[\bX_{j,1:q} \in A^{q}\mid \Pi_{q}^{(j)}= \pi \big] \, \bbP\big[\Pi_{q}^{(j)} = \pi \big] \\
	&= \sum_{s=1}^{q} P_{0}(A)^{s} \sum_{ \pi \in \{ \P(q) : K^{(j)}_{1:q} = s \} } \bbP\big[\Pi_{q}^{(j)} = \pi \big]\\
	&=\sum_{s=1}^{q} P_{0}(A)^{s} \bbP(K^{(j)}_{1:q}=s) = \bbE\big[P_{0}(A)^{K^{(j)}_{1:q}}\big].
\end{align*}
\end{proof}

\subsection{Proof of Proposition~\ref{prop: JointMomSSP}}
\begin{proof}
	For notational convenience, we prove the proposition for $h=2$.
	The general case can be proven with the same argument.
	Notation is the same as in the proof of Proposition~\ref{prop: MomSSP}.
	\begin{align*}
		\bbE[P_{j}(A_{1})^{q_{1}} P_{j}(A_{2})^{q_{2}}] = \bbP \left(\bX_{j,1:q} \in A_{1}^{q_{1}} \times A_{2}^{q_{2}} \right),
	\end{align*}
	where $q=q_{1}+q_{2}$.
	Denote now with $\A_{q_{1},q_{2}}\subset  \P(q)
    $ the set of all possible partitions $\Pi_{q}^{(j)}$ induced by the ties in $\bX_{j,1:q}$ such that the elements in $\bX_{j,1:q_{1}}$ and in $\bX_{j,q_{1}+1:q}$ do not have ties.
	It follows that 
	\begin{align*}
		\bbP\left( \bX_{j,1:q} \in A_{1}^{q_{1}} \times A_{2}^{q_{2}} \right) =& \, \bbP \big[(\bX_{j,1:q} \in A_{1}^{q_{1}} \times A_{2}^{q_{2}}) \cap (\Pi_{q}^{(j)}\in \A_{q_{1},q_{2}}) \big]\\
		=& \, \sum\limits_{s_{1}=1}^{q_{1}} \sum\limits_{s_{2}=1}^{q_{2}} \bbP(\Pi_{q}^{(j)}\in \A_{q_{1},q_{2}}, K^{(j)}_{1:q_{1}}=s_{1}, K^{(j)}_{q_{1}+1:q} =s_{2}) \\
	 & \times \bbP\big[\bX_{j,1:q} \in A_{1}^{q_{1}} \times A_{2}^{q_{2}} \mid \Pi_{q}^{(j)}\in \A_{q_{1},q_{2}}, K^{(j)}_{1:q_{1}}=s_{1}, K^{(j)}_{q_{1}+1:q} =s_{2} \big]\\
		=& \, \sum\limits_{s_{1}=1}^{q_{1}} \sum\limits_{s_{2}=1}^{q_{2}} P_{0}(A_{1})^{s_{1}} P_{0}(A_{2})^{s_{2}} \bbP\big[\Pi_{q}^{(j)}\in \A_{q_{1},q_{2}}, K^{(j)}_{1:q_{1}}=s_{1}, K^{(j)}_{q_{1}+1:q} =s_{2} \big] \\
		=& \, \bbE \left[ P_{0}(A_{1})^{K^{(j)}_{1:q_{1}}} P_{0}(A_{2})^{K^{(j)}_{q_{1}+1:q}}\mid \Pi_{q}^{(j)}\in \A_{q_{1},q_{2}}\right] \bbP \big[\Pi_{q}^{(j)}\in \A_{q_{1},q_{2}} \big].
	\end{align*}
Note that we exploit 
\[
\{\bX_{j,1:q} \in A_{1}^{q_{1}}\times A_{2}^{q_{2}}\}\subseteq E_{\neq} \coloneqq  \big\{\Pi_{q}^{(j)}\in \A_{q_{1},q_{2}} \big\}
\]
because $A_{1}\cap A_{2}=\emptyset$.
Thus, $\bbP\big(\bX_{j,1:q} \in A_{1}^{q_{1}}\times A_{2}^{q_{2}}\big)=\bbP\big(E_{\neq} \big)
\bbP\big(\bX_{j,1:q}\in A_{1}^{q_{1}}\times A_{2}^{q_{2}}\mid E_{\neq}\big).$
\end{proof}

\subsection{Proof of Theorem~\ref{thm: MommSSP}}
\begin{proof}
	\begin{align*}
		\bbE[P_{1}(A)^{q_{1}} \cdots P_{J}(A)^{q_{J}}] = \bbP\left(\bX_{j,1:q_{j}} \in A^{q_{j}}: j=1,\ldots,J \right).
	\end{align*}
	Disintegrate with respect to the random partition $\Pi_{q}$ of the labels $[q]$ (with $q = q_{1}+\cdots+q_{J})$ induced by the ties of the observations $\bX_{1:q_{1},\ldots,1:q_{J}}$ (within and across groups) to recover independence and aggregate by symmetry induced by partial exchangeability.
    Let $\P({q_{1},\ldots, q_{J}})$ denote the space of possible partitions of such labels induced by the ties among $\bX_{1:q_{1},\ldots,1:q_{J}}$.
    Thus,
	\begin{align*}
	\bbP\left(\bX_{j,1:q_{j}} \in A^{q_{j}}: j=1,\ldots,J\right) &= \sum_{
     \pi \in \P({q_{1},\ldots, q_{J}}) } \bbP\big(\bX_{j,1:q_{j}} \in A^{q_{j}}: j=1,\ldots,J \mid \Pi_{q} = \pi \big) \bbP\big( \Pi_{q} = \pi \big)\\
    &=\sum_{s=1}^{q} P_{0}(A)^{s} \sum_{ \pi \in \{\P({q_{1},\ldots, q_{J}}) : K_{q_{1},\ldots,q_{J}} = s \}} \bbP \big( \Pi_{q} = \pi \big)\\
	&= \sum_{s=1}^{q} P_{0}(A)^{s} \bbP(K_{q_{1},\ldots,q_{J}} = s) = \bbE\big[P_{0}(A)^{K_{q_{1}, \ldots, q_{J}}} \big].
	\end{align*}
\end{proof}

\subsection{Proof of Theorem~\ref{thm: JointMommSSP}}
\begin{proof}
	First note that
	\begin{align*}
		\bbE\bigg[\prod_{j=1}^{J} P_{j}(A_{j})^{q_{j}} \bigg] = \bbP\bigg(\bX_{1:q_{1},\ldots,1:q_{J}} \in \bigtimes\limits_{j=1}^{J} A_{j}^{q_{j}}\bigg).
	\end{align*}
	Denote now by $\A_{q_{1},\ldots,q_{J}}\subset \P(q_{1}, \ldots, q_{J})$ the set of all possible partitions $\Pi_{q}$ of the labels of elements in $\bX_{1:q_{1}, \ldots, 1:q_{J}}$ such that the labels of elements in $\bX_{j,1:q_{j}}$ and in $\bX_{j^{\prime},1:q_{j^{\prime}}}$ do not belong to the same set, for any $j\neq j^{\prime}$ according to $\Pi_{q}$, with $q = q_{1} + \cdots + q_{J}$.
    Thus, $E_{\neq}\coloneqq \{\Pi_{q}\in \A_{q_{1},\ldots,q_{J}}\}$.
	\begin{align*}
    &\bbP\left[\bX_{1:q_{1},\ldots,1:q_{J}} \in \bigtimes\limits_{j=1}^{J} A_{j}^{q_{j}}\right] = \, \bbP\left[ \left(\bX_{1:q_{1},\ldots,1:q_{J}} \in \bigtimes\limits_{j=1}^{J} A_{j}^{q_{j}} \right) \cap \left(\Pi_{q}\in \A_{q_{1}, \ldots, q_{J}} \right)\right]\\
		&\quad = \, \bbP\left(\Pi_{q}\in \A_{q_{1},\ldots,q_{J}}\right) \bbP\left( \bX_{1:q_{1},\ldots,1:q_{J}} \in \bigtimes\limits_{j=1}^{J} A_{j}^{q_{j}} \mid \Pi_{q}\in \A_{q_{1},\ldots,q_{J}}\right) \\
		& \quad = \, \sum\limits_{s_{1}=1}^{q_{1}} \cdots \sum\limits_{s_{J}=1}^{q_{J}} \bbP \left[\Pi_{q}\in \A_{q_{1},\ldots,q_{J}}, K^{(1)}_{1:q_{1}}=s_{1},\ldots, K^{(J)}_{1:q_{J}} =s_{J}\right] \\
& \quad \quad \times \bbP\left[\bX_{1:q_{1}, \ldots, 1:q_{J}} \in \bigtimes\limits_{j=1}^{J} A_{j}^{q_{j}} \mid \Pi_{q}\in \A_{q_{1},\ldots,q_{J}}, K^{(1)}_{1:q_{1}}=s_{1},\ldots, K^{(J)}_{1:q_{J}} =s_{J}\right] \\
		&\quad = \, \sum\limits_{s_{1}=1}^{q_{1}} \cdots \sum\limits_{s_{J}=1}^{q_{J}} P_{0}(A_{1})^{s_{1}}\cdots P_{0}(A_{J})^{s_{J}}\\
	& \quad \quad	\times \bbP\left[\Pi_{q}\in \A_{q_{1},\ldots,q_{J}}, K^{(1)}_{1:q_{1}} = s_{1}, \ldots, K^{(J)}_{1:q_{J}} =s_{J} \right]\\
		&\quad = \, \bbE\left[P_{0}(A_{1})^{K^{(1)}_{1:q_{1}}} \cdots P_{0}(A_{J})^{K^{(J)}_{1:q_{J}}}\mid \Pi_{q}\in \A_{q_{1}, \ldots, q_{J}}\right] \bbP(\Pi_{q}\in \A_{q_{1}, \ldots, q_{J}}).
	\end{align*}
Note that we exploit that \[
\big\{\bX_{1:q_{1},\ldots,1:q_{J}}\in \bigtimes_{j=1}^{J} A_{j}^{q_{j}}\big\} \subseteq E_{\neq},
\]
because the $A_{j}$'s are pairwise disjoint sets.
\end{proof}

\subsection{Proof of Theorem~\ref{thm: pEPPFchar}}
\begin{proof}
To prove the theorem, we first show that when the random array follows an mSSM, then it can be obtained by sampling first the partition from the corresponding pEPPF and then associating unique values sampled independently from $P_{0}$ to each partition set.
Formally, for any family of measurable sets $(A_{j,i}: i\in [I_{j}], j \in [J])$,
\begin{align*}
	& \bbP \left[(X_{j,i}: i \in [I_{j}], \, j \in [J]) \in (A_{j,i}: i \in [I_{j}], \, j \in [J]) \right] \\
	&\quad = \sum_{\displaystyle{ \genfrac{}{}{0pt}{}{\{C_{1}, \ldots, C_{D}\}}{\in \P(I_{1}, \ldots, I_{J})} }} \bbP[(X_{j,i}: i \in [I_{j}], \, j \in [J]) \in (A_{j,i}: i \in [I_{j}], \, j \in [J]) \mid \Pi_{n} = \{C_{1}, \ldots, C_{D}\} ]\\
    & \hspace*{4cm} \times \bbP(\Pi_{n} = \{C_{1}, \ldots, C_{D}\})\\
	& \quad = \sum_{\displaystyle{ \genfrac{}{}{0pt}{}{\{C_{1}, \ldots, C_{D}\}}{\in \P(I_{1}, \ldots, I_{J})} }} \bbP(\Pi_{n} = \{C_{1}, \ldots, C_{D}\}) \prod_{d=1}^{D} P_{0} \left(\bigcap_{(j,i): (\sum_{k=1}^{j-1} I_{k} + i) \, \in \, C_{d}} A_{j,i}\right),
\end{align*}
where $C_{1}, \ldots, C_{D}$ are partition blocks, whose elements are collected according to the order of arrival by groups.
Note that $\sum_{k=1}^{j-1} I_{k}+i$ is the vectorization map from $([J]\times \bbN)$ to $[n]$, with $n =\sum_{j=1}^{J} I_{j}$ used to identify which $(j,i)$ lie in each block $C_{d}$ in order of arrival.
Now, to complete the proof, what is left to show is that when the pEPPF is obtained from an arbitrary vector of random probability measures according to $\eqref{eq: pEPPF}$ in the main paper, then there always exists an mSSP that induces the same pEPPF.

To this aim, let us first consider a pair of dependent random probability measures (not necessarily mSSP) that admit the following representation 
\begin{equation}\label{eq:genericmeasuresSupp}
    P_{j} \eqas \sum_{h \ge 1} \pi_{j,h} \delta_{\theta_{h}} + \bigg(1-\sum_{h \ge 1} \pi_{j,h}\bigg) P_{0}, \quad \quad \text{for } j=1,2
\end{equation}
where $P_{0}$ is a non-atomic (deterministic) probability measure on a space $\bbX$, $\bpi_{j} = (\pi_{j,h})_{h\geq1}$ is a random sub-probability sequence, for $j=1,2$, and the sequence of $(\theta_{h})_{h\geq1}$, conditionally on $\bpi= (\bpi_{1},\bpi_{2})$, follows any joint distribution such that $\bbP(\theta_{h}=\theta_{\ell}\mid \bpi) = 0$, for any $h \neq \ell$.
Define $\omega_{j} \eqas \sum_{h \ge 1} \pi_{j,h}$, such that
 \[
    P_{j} = \omega_{j} \, \tP_{j} + (1-\omega_{j})P_{0} \quad \text{for } j = 1,2
\]
where $\tP_{j} \eqas \sum\limits_{h \ge 1} \frac{\pi_{j,h}}{\omega_{j}} \delta_{\theta_{h}} {=:} \sum\limits_{h \ge 1} \tpi_{j,h} \delta_{\theta_{h}}$.
Then, applying the binomial expansion
\begin{align*}
		&\pEPPF_{D}^{(n)}(\bn_{1},\bn_{2}) = \bbE \bigg[ \int_{\bbX^{D}_{\star}} \prod_{d=1}^{D} P_{1}(\d x_{d})^{n_{1,d}} P_{2}(\d x_{d})^{n_{2,d}} \bigg]\\
		& \quad = \bbE \bigg\{ \int_{\bbX^{D}_{\star}} \prod_{d=1}^{D} \prod_{j=1}^{2} \bigg[\omega_{j} \, \tP_{j}(\d x_{d}) + (1-\omega_{j})P_{0} (\d x_{d})\bigg]^{n_{j,d}} \bigg\}\\
        & \quad = \bbE \bigg\{ \int_{\bbX^{D}_{\star}} \prod_{d=1}^{D} \prod_{j=1}^{2} \bigg[\sum_{r_{j,d}=0}^{n_{j,d}} \binom{n_{j,d}}{r_{j,d}} \omega_{j}^{r_{j,d}} (1-\omega_{j})^{n_{j,d}-r_{j,d}} \tP_{j}(\d x_{d})^{r_{j,d}} P_{0}(\d x_{d})^{n_{j,d}-r_{j,d}} \bigg] \bigg\}\\
        & \quad = \bbE \bigg\{ \int_{\bbX^{D}_{\star}} \prod_{d=1}^{D} \sum_{r_{1,d}=0}^{n_{1,d}}\sum_{r_{2,d}=0}^{n_{2,d}} \bigg[\binom{n_{1,d}}{r_{1,d}} \binom{n_{2,d}}{r_{2,d}}\omega_{1}^{r_{1,d}} (1-\omega_{1})^{n_{1,d}-r_{1,d}}\omega_{2}^{r_{2,d}} (1-\omega_{2})^{n_{2,d}-r_{2,d}} \\
        & \qquad \qquad \times\tP_{1}(\d x_{d})^{r_{1,d}}\,\tP_{2}(\d x_{d})^{r_{2,d}}\, P_{0}(\d x_{d})^{n_{\cdot,d}-r_{1,d} - r_{2,d}} \bigg] \bigg\},
        \end{align*}
        where for each species $d$, $n_{\cdot,d} = \sum_{j=1}^{2} n_{j,d}$ is the overall frequency of species $d$ in the sample.
        That is, $n_{\cdot,d}$ is the cardinality of the block $C_{d}$.
        Recall that $\bbX_{\star}^{D} =\{(x_{1}, \ldots, x_{D}) \in \bbX^{D} : x_{d} \ne x_{d^{\prime}} \text{ for all } d \ne d^{\prime}\}$.

        Let $\D :=\{d \in [D]: n_{\cdot,d} = 1\}$ denote the singleton species and $\bar{\D} :=\{ d \in [D]: n_{\cdot,d} \ge 2\}$ the non-singleton species.
        Note that $\D$ and $\bar{\D}$ form a partition of $[D]$.
        
        For any non-singleton species $d \in \bar{\D}$, when integrating over $\bbX_{\star}^{D}$, any term of the sum in $(r_{1,d},\,r_{2,d})$ for which $n_{\cdot,d}-r_{1,d} - r_{2,d}>0$ yields zero contribution, by the non-atomicity of $P_{0}$.
        That is, only atomic assignments: $r_{j,d}=n_{j,d}$ for all $j=1,2$ have positive contributions.
        Thus, for any $d \in \bar{\D}$ the contribution simplify to
        \begin{align*}
         \prod_{j=1}^{2} \bigg[\sum_{r_{j,d}=0}^{n_{j,d}} \binom{n_{j,d}}{r_{j,d}} \omega_{j}^{r_{j,d}} (1-\omega_{j})^{n_{j,d}-r_{j,d}} \tP_{j}(\d x_{d})^{r_{j,d}} P_{0}(\d x_{d})^{n_{j,d}-r_{j,d}} \bigg] = \prod_{j=1}^{2}
         \omega_{j}^{n_{j,d}} \tP_{j}(\d x_{d})^{n_{j,d}}.
        \end{align*}
        For any singleton species $d \in \D$, denote by $k(d)$ the label of the sample where the singleton is observed, i.e., $k(d)\in \{1, 2\}$ is such that $n_{k(d),d}=1$.
        Note that exactly one $k(d) \in \{1, 2\}$ is such that $n_{k(d),d}=1$ for any $d \in \D$.
        Thus, for any $d \in \D$ the contribution simplify to
        \begin{align*}
            &\prod_{j=1}^{2} \bigg[\sum_{r_{j,d}=0}^{n_{j,d}} \binom{n_{j,d}}{r_{j,d}} \omega_{j}^{r_{j,d}} (1-\omega_{j})^{n_{j,d}-r_{j,d}} \tP_{j}(\d x_{d})^{r_{j,d}} P_{0}(\d x_{d})^{n_{j,d}-r_{j,d}} \bigg]\\
            & \quad = \omega_{k(d)} \tP_{k(d)}(\d x_{d}) + (1-\omega_{k(d)}) P_{0}(\d x_{d}).
        \end{align*}
        Therefore,
        \begin{align*}
        \pEPPF_{D}^{(n)}(\bn_{1},\bn_{2})= \bbE \bigg\{ \int_{\bbX^{D}_{\star}} \prod_{d\in \D} \bigg[ \omega_{k(d)} \tP_{k(d)}(\d x_{d}) + (1-\omega_{k(d)}) P_{0}(\d x_{d}) \bigg] \prod_{d\in \bar{\D}} \prod_{j=1}^{2} \omega_{j}^{n_{j,d}} \tP_{j}(\d x_{d})^{n_{j,d}} \bigg\}.
        \end{align*}
        Now, let $\B \subseteq \D$ denote the subset of singleton species arising from the discrete part of the mSSP and $\bar{\B}:= \D \setminus \B$, that is, the subset of singleton species arising from the non-atomic part of the mSSP.
        We can rewrite
        \begin{align*}
        &\pEPPF_{D}^{(n)}(\bn_{1},\bn_{2}) = \bbE \bigg\{ 
        \int_{\bbX^{D}_{\star}} \prod_{d\in \D} \bigg[
        \omega_{k(d)} \tP_{k(d)}(\d x_{d}) + (1-\omega_{k(d)}) P_{0}(\d x_{d}) \bigg] \prod_{d\in \bar{\D}} \prod_{j=1}^{2} \omega_{j}^{n_{j,d}} \tP_{j}(\d x_{d})^{n_{j,d}} \bigg\} \\
        &\quad = \bbE \bigg\{ \sum_{\B: \B \subseteq \D} \bigg[ \prod_{d\in \bar{\B}} (1- \omega_{k(d)}) \bigg] \sum_{\displaystyle{ \genfrac{}{}{0pt}{}{(h_{d})_{d \in \bar{\D} \cup \B}}{ \text{\scriptsize{pairwise distinct}}}} } \bigg[\prod_{d\in \bar{\D}} \prod_{j=1}^{2}\big(\omega_{j} \tpi_{j,h_{d}} \big)^{n_{j,d}} \bigg] \bigg[ \prod_{d \in \B} \omega_{k(d)} \tpi_{k(d),h_{d}} \bigg] \bigg\}.
        \end{align*}
        
        Importantly, the expression and derivation above for the pEPPF also hold when $(P_{1}, P_{2})$ is an mSSP.
        Crucially, this expression depends only on the law of the weights of the random probability measures in \eqref{eq:genericmeasuresSupp}.
        Since the class of mSSPs imposes no restriction on this law, one can always choose an mSSP whose pEPPF matches the expression for any prescribed weight distribution.
        More precisely, given a pEPPF derived from $P_{1}$ and $P_{2}$ in \eqref{eq:genericmeasuresSupp}, the corresponding mSSP may be constructed by adopting the same weight law of $P_{1}$ and $P_{2}$ and selecting an arbitrary non-atomic base measure.

        To extend the argument to arbitrary pairs of probabilities not necessarily satisfying \eqref{eq:genericmeasuresSupp}, note that they can always be decomposed into an atomic and a non-atomic component.
        More precisely, let $(G_{1}, G_{2})$ be an arbitrary pair of random probability measures.
        For each $j = 1, 2$ we can rewrite
        \[
            G_{j} \eqas \sum_{h \ge 1} \pi_{j,h} \delta_{\theta_{h}} + \bigg(1-\sum_{h \ge 1} \pi_{j,h}\bigg) G_{0,j}, \quad \quad \text{for }j=1,2
        \]
        where $G_{0,j}$ is a non-atomic random probability measure on a space $\bbX$, $\bpi_{j} = (\pi_{j,h})_{h\geq1}$ is a random sub-probability sequence, for $j=1,2$, and the sequence of $(\theta_{h})_{h\geq1}$, conditionally on $\bpi= (\bpi_{1},\bpi_{2})$, follows any joint distribution such that $\bbP(\theta_{h}=\theta_{\ell}\mid \bpi) = 0$, for any $h \neq \ell$.
        Note that the pair $(G_{1}, G_{2})$ induces the same pEPPF if we substitute the random non-atomic probability measures $G_{0,1}$ and $G_{0,2}$ with a common deterministic arbitrary probability measure $P_{0}$, since applying the same argument above, the pEPPF is a functional of the weights only.

        Thus, the pEPPF induced by an arbitrary pair of random probabilities $(G_{1}, G_{2})$ can be obtained as the pEPPF induced by a pair of random probabilities $(P_{1}, P_{2})$ satisfying \eqref{eq:genericmeasuresSupp} and thus by a pair of mSSP.
        Finally, note that extending the above argument to $J > 2$ is a matter only of notation.
\end{proof}

\subsection{Proof of Proposition~\ref{prop: peppf from weight}}
\begin{proof}
The proof follows by reasoning analogous to the derivation of the pEPPF for generic measures in the proof of Theorem~\ref{thm: pEPPFchar}.
However, since here we assume the mSSP is proper, its non-atomic components vanish, and the derivation simplifies to:
	\begin{align*}
		\pEPPF_{D}^{(n)}(\bn_{1},\ldots,\bn_{J}) &= \bbE \bigg[ \int_{\bbX^{D}_{\star}} \prod_{d=1}^{D} P_{1}(\d x_{d})^{n_{1,d}} \ldots P_{J}(\d x_{d})^{n_{J,d}} \bigg]\\
		& = \bbE \bigg\{ \int_{\bbX^{D}_{\star}} \prod_{d=1}^{D} \prod_{j=1}^{J} \bigg[\sum_{h\ge 1} \pi_{j,h} \delta_{\theta_{h}}(\d x_{d})\bigg]^{n_{j,d}} \bigg\}\\
        &= \bbE \bigg[ \sum_{h_{1} \ne \cdots \ne h_{D}} \prod_{j=1}^{J} \prod_{d=1}^{D} \pi_{j,h_{d}}^{n_{j,d}} \bigg].
	\end{align*}
\end{proof}

\subsection{Proof of Proposition~\ref{prop: urn from peppf}}
The proof follows trivially from the definition of conditional probability.

\section{Algorithms and models details for the multi-armed bandit illustration}\label{sec: algo supp}
The algorithms used for all six strategies considered in Section~\ref{sec: illustration} are Markov chain Monte Carlo marginal algorithms.
These algorithms are obtained using the augmented representation of the pEPPF described in Section~\ref{sec: pred} for the additive and hierarchical processes, and the sequential sampling schemes detailed in Section~\ref{sec: SSP supp} for the independent processes.
All models are generalized to accommodate random hyperparameters to achieve greater flexibility in the learning mechanisms, leading to the following specifications for the six strategies.
Throughout, $P_{0}$ is an arbitrary deterministic non-atomic probability measure only used to generate the species labels.
\begin{itemize}
	\item Independent Dirichlet process
	\begin{align*}
			X_{j,i} \mid (P_{1}, \ldots, P_{J}) &\simind P_{j} \quad \text{ for } j \in [J] \text{ and } i=1,2, \ldots\\
			P_{j} \mid \alpha_{j} &\simind \DP (\alpha_{j}, P_{0}) \quad \text{ for } j \in [J]\\
			\alpha_{j} &\simiid \Ga (0.75, 1).
	\end{align*}
	where $\Ga (a, b)$ denotes a gamma distribution with expected value equal to $a / b$.
	\item Independent Pitman--Yor process
	\begin{align*}
			X_{j,i} \mid (P_{1}, \ldots, P_{J}) &\simind P_{j} \quad \text{ for } j \in [J] \text{ and } i=1,2, \ldots\\
			P_{j} \mid \sigma_{j}, \alpha_{j} &\simind \PYP (\sigma_{j}, \alpha_{j}, P_{0}) \quad \text{ for } j \in [J]\\
			\sigma_{j} \simiid \Beta(1, 3)\qquad &
			\alpha_{j} \simiid \Ga (0.2, 1).
	\end{align*}
	where $\Beta$ denotes a beta distribution.
	\item Additive Dirichlet process
	\begin{align*}
			X_{j,i} \mid (P_{1}, \ldots, P_{J}) &\simind P_{j} \quad \text{ for } j \in [J] \text{ and } i=1,2, \ldots\\
			P_{j} = \epsilon_{j} \, Q_{0} + (1-\epsilon_{j}) Q_{j}, &\quad 
			\epsilon_{j} \simiid 0.15 \, \delta_{0} + 0.15 \, \delta_{1} + 0.7 \, \Unif(0,1) \quad \text{ for } j \in [J]\\
			Q_{j} \mid \alpha_{j} &\simind \DP (\alpha_{j}, P_{0}) \quad \text{ for } j=0,1, \ldots,J\\
			\alpha_{0} \simiid \Ga (0.5, 2),\qquad
			&\alpha_{j} \simiid \Ga (6, 2) \quad\text{ for } j \in [J].
	\end{align*}
	\item Additive Pitman--Yor process
	\begin{align*}
		X_{j,i} \mid (P_{1}, \ldots, P_{J}) &\simind P_{j} \quad \text{ for } j \in [J] \text{ and } i=1,2, \ldots\\
		P_{j} = \epsilon_{j} \, Q_{0} + (1-\epsilon_{j}) Q_{j}, &\quad \epsilon_{j} \simiid 0.15 \, \delta_{0} + 0.15 \, \delta_{1} + 0.7 \, \Unif(0,1) \quad \text{ for } j \in [J]\\
			Q_{j} \mid \sigma_{j}, \alpha_{j} &\simind \PYP(\sigma_{j}, \alpha_{j}, P_{0}) \quad \text{ for } j=0,1, \ldots,J\\
			\sigma_{0} \simiid \Beta(1, 3), \qquad
			&\sigma_{j} \simiid \Beta(1, 2)\quad\text{ for } j \in [J]\\
			\alpha_{0} \simiid \Ga (0.25, 4), \qquad
			&\alpha_{j} \simiid \Ga (2, 2) \quad\text{ for } j \in [J].
	\end{align*}
	\item Hierarchical Dirichlet process
	\begin{align*}
			X_{j,i} \mid (P_{1}, \ldots, P_{J}) &\simind P_{j} \quad \text{ for } j \in [J] \text{ and } i=1,2, \ldots\\
			P_{j}\mid Q, \alpha_{j} &\simind \DP(\alpha_{j},Q) \quad \text{ for } j \in [J] \\
            Q \mid \alpha_{0} &\sim \DP(\alpha_{0},P_{0})\\
			\alpha_{0} \simiid \Ga (1, 1/3), \qquad
			&\alpha_{j} \simiid \Ga (1, 1/2) \quad\text{ for } j \in [J].
	\end{align*}
	\item Hierarchical Pitman--Yor process
	\begin{align*}
			X_{j,i} \mid (P_{1}, \ldots, P_{J}) &\simind P_{j} \quad \text{ for } j \in [J] \text{ and } i=1,2, \ldots\\
			P_{j}\mid Q, \sigma_{j},\alpha_{j} &\simind \PYP(\sigma_{j},\alpha_{j},Q) \quad \text{ for } j \in [J] \\ 
			Q \mid \sigma_{0},\alpha_{0} &\sim \PYP(\sigma_{0},\alpha_{0},P_{0})\\
			\sigma_{0} \simiid \Beta(1, 2), \qquad
			&\sigma_{j} \simiid \Beta(1, 2)\quad\text{ for } j \in [J]\\
			\alpha_{0} \simiid \Ga (1, 1), \qquad
			&\alpha_{j} \simiid \Ga (1, 1) \quad\text{ for } j \in [J].
	\end{align*}
\end{itemize}
The choice of the parameters of the hyperpriors on discount and concentration parameters is performed as follows.
We use the values suggested in \textcite{battiston2018multi} for the hierarchical Pitman--Yor process, and then we fix the ones of the other strategies by considering the probabilities of ties as a function of the hyperparameters and approximately match their expected values and variances.
This selection procedure ensures a fair performance comparison, as the probabilities of ties provide an excellent summary of dependence for rmSSPs.
The resulting expected probability of ties and corresponding variances are reported in Table~\ref{tab:hyperparam}.

\begin{table}[ht]
	\centering
	\resizebox{\textwidth}{!}{
		\begin{tabular}{lcccc}
			\toprule
			\textbf{Model} & \textbf{$\bbE[\text{prob tie across}]$} & \textbf{$\mathbb{V}[\text{prob tie across}]$} & \textbf{$\bbE[\text{prob tie within}]$} & \textbf{$\mathbb{V}[\text{prob tie within}]$} \\
			\midrule
			Independent DP & 0 & 0 & 0.672 & 0.049\\
			Independent PYP & 0 & 0 & 0.669 & 0.047\\
			+DP & 0.388 & 0.092 & 0.666 & 0.052\\
			+PY & 0.400 & 0.064 & 0.628 & 0.038\\
			HDP & 0.389 & 0.056 & 0.671 & 0.041 \\
			HPY & 0.397 & 0.043 & 0.638 & 0.033 \\
			\bottomrule
	\end{tabular}
    }
\caption{ \label{tab:hyperparam} Expected probabilities of ties across and within as functions of the hyperparameters and corresponding variances.
Values are obtained via Monte Carlo approximation by simulating $2,000$ samples of the hyperparameters from the hyperpriors of each model.}
\end{table}
To sample the concentration parameters of the Dirichlet processes, we employed a Gibbs sampler via an augmented representation of the full-conditional of the concentration parameter, avoiding a Metropolis within the Gibbs step.
For the hyperparameters of the Pitman--Yor processes, we devised an adaptive Metropolis--Hastings, obtained via $10$ repeated steps within the main Gibbs algorithm.
At each of the $300$ sequential sampling steps of the multi-armed bandit problem, we perform $200$ iterations of the MCMC algorithm, leading to a total of $60,000$ iterations (not including the possible additional Metropolis--Hastings steps, when present) per strategy.
After observing a new data point in a sequential step, we initialize the MCMC for the next step with a warm start based on the last iteration of the MCMC output in the previous step.
For instance, we initialize the values of the hyperparameters with the last sampled value in the previous MCMC chain that targets their posterior distribution without conditioning on the new data point.

Moreover, for hierarchical processes, we perform $1,000$ iterations of the MCMC before estimating the probability of discovery at the first sequential sampling step to also achieve a warm start at the first sampling step.
Code for all six strategies is freely available at \if0\blind the GitHub page of the authors.\else \url{https://github.com/GiovanniRebaudo/MSSP}.\fi 

\printbibliography
\end{document}